\documentclass[12pt,a4paper,twoside]{article}
\usepackage[french] {babel}
\usepackage{amssymb,amsmath,amsfonts,graphics}
\usepackage{stmaryrd}
\usepackage[all]{xy}

\def\rondI{\newbox\boxx{\hbox{$1$}}\hskip-11pt\raisebox{0pt}{$\displaystyle{\textrm{\Large$\bigcirc$}}$}}
\def\rondII{\newbox\boxx{\hbox{$2$}}\hskip-11pt\raisebox{0pt}{$\displaystyle{\textrm{\Large$\bigcirc$}}$}}
\def\rondIII{\newbox\boxx{\hbox{$3$}}\hskip-11pt\raisebox{0pt}{$\displaystyle{\textrm{\Large$\bigcirc$}}$}}
\def\rondIV{\newbox\boxx{\hbox{$4$}}\hskip-11pt\raisebox{0pt}{$\displaystyle{\textrm{\Large$\bigcirc$}}$}}

\begin{document}

\newcommand{\hooklongrightarrow}{\lhook\joinrel\longrightarrow}

\begin{center}
\textbf{Images directes II: $F$-isocristaux convergents}

\vskip20mm

Jean-Yves ETESSE
  \footnote{(CNRS - Institut de Math\'ematique, Universit\'e de Rennes 1, Campus de Beaulieu - 35042 RENNES Cedex France)\\
E-mail : Jean-Yves.Etesse@univ-rennes1.fr}
\end{center}

 \vskip40mm
\noindent\textbf{Sommaire}\\

	\begin{enumerate}
		\item[0.] Introduction
		\item[1.] $F$-isocristaux convergents sur un sch\'ema affine et lisse
			\begin{enumerate}
			\item[1.1.] Notations
			\item[1.2.] Des \'equivalences de cat\'egories
			\end{enumerate}
		\item[2.] $F$-isocristaux convergents sur un sch\'ema lisse formellement\\ relevable
			\begin{enumerate}
			\item[2.1.] Espaces rigides associ\'es aux sch\'emas formels
			\item[2.2.] $F$-isocristaux convergents et $\mathcal{O}_{\mathcal{X}_{K}}$-modules
			\end{enumerate}
		\item[3.] Rel\`evement de Teichm¬\"uller, r\'eduction, sp\'ecialisation et espaces rigides analytiques
			\begin{enumerate}
			\item[3.1.] Rel\`evements de Teichm\" uller
			\item[3.2.] Morphismes de sp\'ecialisation et de r\'eduction
			\item[3.3.] $Spm\ K\{T\}$: la boule unit\'e de la droite affine rigide
			\item[3.4.] Exemples			
			\end{enumerate}
		\item[4.] Images directes de $F$-isocristaux convergents
			\begin{enumerate}
			\item[4.1.] Convergence des images directes
			\item[4.2.] Cas fini \'etale
			\end{enumerate}
		\end{enumerate}

\newpage

\noindent\textbf{R\'esum\'e}\\

Cet article est le deuxi\`eme  d'une s\'erie de trois articles consacr\'es aux images directes d'isocristaux: ici nous consid\'erons des isocristaux convergents avec structure de Frobenius.\\

Soit $\mathcal{V}$ un anneau de valuation discr\`ete complet, de corps r\'esiduel $k = \mathcal{V}/\mathfrak{m}$ de caract\'eristique $p > 0$ et de corps des fractions $K$ de caract\'eristique $0$. Dans un premier temps nous caract\'erisons les $F$-isocristaux convergents sur un sch\'ema affine et lisse sur $k$. Dans un deuxi\`eme temps, pour $k$ parfait et apr\`es avoir explicit\'e en d\'etail les rel\`evements de Teichm\" uller, notamment pour la droite affine rigide, nous en d\'eduisons l'existence d'isomorphismes de Frobenius sur les images directes de $F$-isocristaux convergents par un $k$-morphisme propre et lisse relevable.\\

\noindent\textbf{Abstract}\\

This article is the second one of a series of three articles devoted to direct images of isocrystals: here we consider convergent isocrystals with Frobenius structure.\\

Let $\mathcal{V}$ be a complete discrete valuation ring, with residue field $k = \mathcal{V}/\mathfrak{m}$ of characteristic $p > 0$ and fraction field $K$ of characteristic $0$. Firstly we characterize convergent $F$-isocrystals on a smooth affine $k$-scheme. Secondly, for perfect $k$ and after a detailed exposition of the Teichm\" uller liftings, especially for the affine rigid line, we derive the existence of Frobenius isomorphisms on the direct images of convergent $F$-isocrystals under a proper smooth and liftable $k$-morphism.\\


Mots cl\'es: $F$- isocristaux convergents, cohomologie cristalline, cohomologie convergente, cohomologie rigide, images directes.\\


\newpage 

\section*{0. Introduction}
Sauf mention contraire, dans tout cet article on d\'esigne par $\mathcal{V}$ un anneau de valuation discr\`ete complet, de corps r\'esiduel $k = \mathcal{V}/\mathfrak{m}$ de caract\'eristique $p > 0$, de corps des fractions $K$ de caract\'eristique $0$, d'indice de ramification $e$, et $\pi$ une uniformisante.\\

Cet article est le deuxi\`eme d'une s\'erie de trois articles consacr\'es aux images directes d'isocristaux. Ici on introduit une structure de Frobenius sur les isocristaux consid\'er\'es dans [Et 6], tout en restant dans le cadre \guillemotleft convergent\guillemotright. Les r\'esultats concernent trois aspects de la th\'eorie: d'une part l'allure des $F$-isocristaux convergents sur un sch\'ema lisse aux \S1 et \S2, d'autre part la convergence des images directes par un morphisme propre et lisse relevable $f: X\rightarrow S$ de $k$-sch\'emas au \S4, et dans l'\'etape interm\'ediaire du \S3 nous explicitons les rel\`evements de Teichm\" uller, leurs liens avec la sp\'ecialisation et la r\'eduction, en d\'eveloppant le cas de la boule unit\'e de la droite affine rigide. Ainsi les fonctions $L$ pour de telles images directes convergentes sont d\'efinies et on \'etudiera leurs propri\'et\'es ult\'erieurement.\\

Notre approche via les isocristaux n\'eglige la torsion dans la cohomologie cristalline; d'autres approches la prennent en compte, chez Ogus [O 4] via les $T$-cristaux, ou chez Shiho ([Shi 1], [Shi 2], [Shi 3]) via la cohomologie log-cristalline.\\

Au \S1 on obtient une description des $F$-isocristaux convergents sur un sch\'ema affine et lisse [cor (1.2.3)] qui est l'analogue de celle de type Monsky-Washnitzer obtenue par Berthelot dans le cas surconvergent [B 3]. Cette description s'\'etend au cas lisse relevable dans le \S2 [th\'eo (2.2.1)].\\

 Au \S3 on \'etudie en d\'etail les rel\`evements de Teichm\" uller qui joueront un r\^ole crucial au \S4 dans la preuve que le morphisme de Frobenius sur les images directes est un isomorphisme: une fois d\'efini le morphisme de Frobenius sur l'image directe (qui est un isocristal convergent), un r\'esultat fondamental de Bosch-G\" untzer-Remmert nous dit qu'il sera un isomorphisme s'il en est ainsi en passant aux fibres; le point cl\'e est que l'image inverse par ces rel\`evements de Teichm\" uller est la r\'ealisation du foncteur fibre qui commute aux Frobenius, ce qui permet de se ramener \`a la cohomologie rigide pour laquelle il est connu que le Frobenius est un isomorphisme. On est amen\'e \`a expliciter le morphisme de sp\'ecialisation utilis\'e par Berthelot et \`a le comparer avec celui utilis\'e par Tate [Ta] ou celui utilis\'e par Bosch-G\" untzer-Remmert [B-G-R]. On montre que le rel\`evement de Teichm\" uller est une section du morphisme de sp\'ecialisation et du morphisme de r\'eduction de [B-G-R]: une attention particuli\`ere est accord\'ee aux points \`a bonne r\'eduction. On termine ce \S3 par l'\'etude  de la boule unit\'e ferm\'ee de la droite affine rigide et par des exemples explicites de points dans les tubes. \\

 Si le morphisme propre et lisse $f: X\rightarrow S$ est relevable en un morphisme propre et lisse sur $\mathcal{V}$, les images directes de $F$-isocristaux convergents sont des $F$-isocristaux convergents [th\'eo (4.1.2)]. En fait, si le morphisme $f: X\rightarrow S$ est projectif et lisse et que, ou bien $f$ est relevable, ou bien que $X$ est une intersection compl\`ete relative dans des espaces projectifs sur $S$, ce r\'esultat demeure: la preuve repose sur le cas plongeable de [Et 6, 3.4]. De plus les images directes commutent au passage aux fibres. Dans le cas fini \'etale les restrictions pr\'ec\'edentes sont lev\'ees [th\'eo (4.2.1)], et dans le cas galoisien la cohomologie convergente commute aux points fixes sous le groupe du rev\^etement [th\'eo (4.2.1)].\\

\textbf {Notations:}
 Pour les notions sur les espaces rigides analytiques et la cohomologie rigide nous renvoyons le lecteur \`a [B 3], [B 4], [B-G-R], [C-T] et [LS].\\

On suppose donn\'e un entier $a \in \mathbb{N}^{\ast}$ et on d\'esigne par $C(k)$ un anneau de Cohen de $k$ de caract\'eristique $0$ [Bour, AC IX, $\S$ 2, $\no3$, prop 5] : $C(k)$ est un anneau de valuation discr\`ete complet d'id\'eal maximal $p\ C(k)$ [EGA $O_{IV}$, 19.8.5] et on note $K_{0}$ son corps des fractions, $K_{0} =$ Frac $(C(k))$. Il existe une injection fid\`element plate $C(k) \hookrightarrow \mathcal{V}$ qui fait de $\mathcal{V}$ un $C(k)$-module libre de rang $e$ [EGA $O_{IV}$, 19.8.6, 19.8.8] et [Bour, AC IX, $\S$ 2, $\no1$, prop 2]. On fixe un rel\`evement $\sigma : C(k) \rightarrow C(k)$ de la puissance $a^{i \grave{e}me}$ du Frobenius absolu de $k$ comme dans  [Et 3, I, 1.1]; on suppose que l'on peut \'etendre $\sigma$ en un endomorphisme de $\mathcal{V}$, encore not\'e $\sigma$, de telle sorte que $\sigma(\pi) = \pi$; on notera encore $\sigma$ l'extension naturelle de $\sigma$ \`a $K$: lorsque $k$ est parfait $C(k)$ est isomorphe \`a l'anneau $W(k)$ des vecteurs de Witt de $k$ et $\sigma$ est un automorphisme de $K$. Si $k \hookrightarrow k'$ est une extension de corps de caract\'eristique $p > 0$, $\mathcal{V'} : = \mathcal{V} \otimes_{C(k)} C(k'), K' = \operatorname{Frac} (\mathcal{V'})$, on peut relever la puissance $a^{i \grave{e}me}$ du Frobenius absolu de $k'$ en un morphisme $\sigma' : K' \rightarrow K'$ au-dessus de $\sigma : K \rightarrow K$ [Et 3, I, 1.1].\\

\section*{1. $F$-isocristaux convergents sur un sch\'ema affine et lisse}

\subsection*{1.1. Notations}

Soient $X = \mbox{Spec}\  A_{0}$ un $k$-sch\'ema affine et lisse, $\mathcal{A}$ une $C(k)$-alg\`ebre lisse relevant $A_{0}$ [E$\ell$, th\'eo 6] et $A = \mathcal{A}\ \otimes_{C(k)} \mathcal{V}$. Fixons une pr\'esentation

$$\mathcal{A} = C(k) [T_{1}, ..., T_{n}]\ /\ (f_{1}, ... , f_{s}) ; $$
soient $P$ le compl\'et\'e formel de la fermeture projective de $\mathcal{X} = \mbox{Spec}\  A$ dans $\mathbb{P}^n_{\mathcal{V}}$, $Y$ sa r\'eduction sur $k$ et $j : X\ \hookrightarrow\ Y. $ Alors $]Y[_{P}\  = P_{K}$ ; de plus $]X[_{P}$, qui est l'intersection de $\mathcal{X}^{an}_{K}$ avec la boule unit\'e $B(0,1^+) \subset \mathbb{A}^n_{K}$, est l'affino¬\"{\i}de 

$$]X[_{P}\  = \mbox{Spm}\ (K\{T_{1},...,T_{n}\}\ /\ (f_{1},...,f_{s})). $$ 
Notons $\hat{\mathcal{A}}$ (resp. $\hat{A}$) le s\'epar\'e compl\'et\'e $p$-adique de $\mathcal{A}$ (resp. $A$), $\mathcal{A}^{\dag} \subset \hat{\mathcal{A}}$ (resp. $A^{\dag} \subset \hat{A})$ le compl\'et\'e faible de $\mathcal{A}$ au-dessus de $(C(k), (p))$ (resp. de $A$ au-dessus de $(\mathcal{V},(\pi))$ [M-W, \S\ 1], et posons $A^{\dag}_{K} = A^{\dag}\ \otimes_{\mathcal{V}}\ K, \hat{A}_{K}\ =\ \hat{A}\ \otimes_{\mathcal{V}}\ K$. On a des isomorphismes\\

\qquad $\hat{\mathcal{A}}\ \simeq\ C(k) \{T_{1},...,T_{n} \} / (f_{1},...,f_{s}) $ , \\

\qquad $\hat{A}\ \simeq\ \mathcal{V} \{T_{1},...,T_{n} \} / (f_{1},...,f_{s})\ \simeq\  \hat{\mathcal{A}} \otimes_{C(k)} \mathcal{V}$ ,\\

\qquad $\mathcal{A}^{\dag}\ \simeq\ C(k) [T_{1},...,T_{n}]^{\dag} / (f_{1},...,f_{s})$ ,\\

\qquad $A^{\dag}\ \simeq\ \mathcal{V}[T_{1},...,T_{n}]^{\dag} / (f_{1},...,f_{s})\ \simeq\ \mathcal{A}^{\dag} \otimes_{C(k)} \mathcal{V}$ ,\\

\noindent et aussi [B 3, (2.1.2.4]\\

\qquad $\Gamma(P_{K}, j^{\dag} \mathcal{O}_{P_{K}})\ \simeq\ \Gamma(\mathcal{X}^{an}_{K}, j^{\dag} \mathcal{O}_{P_{K}})\ \simeq\ A^{\dag}_{K} $,\\

\qquad $\Gamma(] X [,  j^{\dag} \mathcal{O}_{P_{K}})\ \simeq\ \Gamma (]X[, \mathcal{O}_{P_{K}})\ =\  \Gamma(]X[, \mathcal{O}_{]X[})\ \simeq\ \hat{A}_{K}$.\\

On fixe un rel\`evement $F_{\mathcal{A}^{\dag}} : \mathcal{A}^{\dag}\ \rightarrow\ \mathcal{A}^{\dag}$ de l'\'el\'evation \`a la puissance $p^a$, $F_{A_{0}} : A_{0}\ \rightarrow\ A_{0}$, au-dessus de $\sigma$ [vdP, cor 2.4.3] : on peut choisir un tel rel\`evement $F_{\mathcal{A}^{\dag}}$ de mani\`ere compatible \`a une extension $k\ \hookrightarrow\ k'$ du corps de base [Et 3, I, 1.2]. Posons $F_{A^{\dag}} = F_{\mathcal{A}^{\dag}} \otimes_{\sigma\vert C(k)} \sigma $, et $F_{\hat{A}}$ le s\'epar\'e compl\'et\'e $p$-adique de $F_{A^{\dag}}$; d'o\`u des morphismes $F_{A^{\dag}_{K}} : A^{\dag}_{K}\ \rightarrow\  A^{\dag}_{K}$,
$F_{\hat{A}_{K}} : \hat{A}_{K}\ \rightarrow\  \hat{A}_{K}$ au-dessus de $\sigma : K \rightarrow\ K$. De m\^eme pour $\mathcal{V}'$ comme ci-dessus [cf notations en fin de l'introduction] et $A' := A\ \otimes_{\mathcal{V}}\ \mathcal{V}'$ il existe d'apr\`es [vdP, cor 2.4.3] un carr\'e commutatif

$$
\xymatrix{
A'^{\dag} \ar[r]^{F_{A'^{\dag}}} & A'^{\dag} \\
A^{\dag} \ar@{^{(}->}[u] \ar[r]_{F_{A^{\dag}}} & A^{\dag} \ar@{^{(}->}[u] 
}
$$

\noindent au-dessus du carr\'e commutatif

$$
\xymatrix{
\mathcal{V} \ar[r]^{\sigma'}& \mathcal{V}' \\
\mathcal{V} \ar@{^{(}->}[u] \ar[r]^{\sigma} & \mathcal{V}.  \ar@{^{(}->}[u] 
}
$$

On d\'esigne par $\mathbf{F^a \textrm{\bf-Mod}(A^{\dag}_{K})}$ (resp. $ \mathbf{F^a \textrm{\bf-Modloc}(A^{\dag}_{K})} $; 
resp. $\mathbf{F^a\textrm{\bf-Modlib}(A^{\dag}_{K})}$) la cat\'egorie des $A^{\dag}_{K}$-modules de type fini (resp. $A^{\dag}_{K}$-modules projectifs de type fini ; resp. $A^{\dag}_{K}$-modules libres de type fini) $M$ munis d'un morphisme

$$\phi_{M} : F^{\ast}_{A^{\dag}_{K}}(M) =: M^{\sigma} \longrightarrow\ M$$

\noindent appel\'e morphisme de Frobenius. Par analogie avec Wan [W 2, def 2.8] [W 3, def 2.1] nous dirons que $M$ est un $\mathbf{F^a\textrm{\bf-module surconvergent}}$ (resp. et $\textbf{projectif}$ ; resp. et $\textbf{libre}$) sur $A^{\dag}_{K}$. En consid\'erant les m\^emes d\'efinitions sur $\hat{A}_{K}$ au lieu de $A^{\dag}_{K}$ on dira que $\mathcal{M} \in \mathbf{F^a \textrm{\bf-Mod}(\hat{A}_{K})}$ (resp. 
$\mathcal{M} \in \mathbf{F^a\textrm{\bf-Modloc}(\hat{A}_{K})}$ ; resp. $\mathcal{M} \in \mathbf{F^a \textrm{\bf-Modlib}(\hat{A}_{K}))}$ muni de

$$\phi_{\mathcal{M}} : F^{\ast}_{\hat{A}_{K}}(\mathcal{M}) =: \mathcal{M}^{\sigma}\ \longrightarrow\ \mathcal{M}$$

\noindent est un $\mathbf{F^a\textrm{\bf-module convergent}}$ (resp. et \textbf{projectif} ; resp. et \textbf{libre}) sur $\hat{A}_{K}$.\\

Lorsque le Frobenius est un isomorphisme on dira que l'on a une \textbf{structure de Frobenius forte}. On a des notions analogues sur $A^{\dag}$ et $\hat{A}$, sans tensoriser par $K$ : on utilisera alors les notations $\mathbf{F^a\textrm{\bf-Mod}(A^{\dag})},... \mathbf{F^a\textrm{\bf-Modlib}(\hat{A})}$.\\

Soit $\Omega^1_{A^{\dag}}$ le module des $\mathcal{V}$-diff\'erentielles de $A^{\dag}$ au sens de Monsky-Washnitzer [M-W, theo 4.2]

$$\Omega^1_{A^{\dag}} :=  \Omega^1_{A^{\dag}/\mathcal{V}}\ /\ \displaystyle \mathop{\bigcap_{n}}\ \mathfrak{m}^n\ \Omega^1_{A^{\dag}/\mathcal{V}}\  ,$$

$$\Omega^1_{A^{\dag}_{K}} :=  \Omega^1_{A^{\dag}} \otimes_{\mathcal{V}} K\ ,\  \Omega^1_{\hat{A}} : = \widehat{\Omega^1_{\hat{A}/\mathcal{V}}}\ ,\ \Omega^1_{\hat{A}_{K}} := \Omega^1_{\hat{A}} \otimes_{\mathcal{V}} K.$$

\noindent Notons $\mathbf{F^a\textrm{\bf-Conn}^{\dag}(A^{\dag}_{K})}$ (resp. $\mathbf{F^a\textrm{\bf-Conn}(\hat{A}_{K}))}$ la cat\'egorie des $A^{\dag}_{K}$-modules (resp. $\hat{A}_{K}$-modules) projectifs de type fini $M$ (resp. $\mathcal{M})$ \`a connexion int\'egrable

$$\nabla : M \longrightarrow M \otimes_{A^{\dag}_{K}}  \Omega^1_{A^{\dag}_{K}}$$

$$( \mbox{resp.} \hat{\nabla} : \mathcal{M} \longrightarrow \mathcal{M} \otimes_{\hat{A}_{K}}  \Omega^1_{\hat{A}_{K}})$$

\noindent et munis d'un isomorphisme horizontal

$$\phi^{\dag} : (F^{\ast}_{A^{\dag}_{K}}(M), \ F^{\ast}_{A^{\dag}_{K}}(\nabla)) =: (M^{\sigma}, \nabla^{\sigma})\  \tilde{\longrightarrow}\ (M, \nabla) $$

$$(\mbox{resp.}\  \hat{\phi} : (F^{\ast}_{\hat{A}_{K}}(\mathcal{M}),\  F^{\ast}_{\hat{A}_{K}}(\hat{\nabla}) =: (\mathcal{M}^{\sigma}, \hat{\nabla}^{\sigma})\ \tilde{\longrightarrow}\ (\mathcal{M}, \hat{\nabla})    ).$$

\noindent On note $\mathbf{\textrm{\bf Conn}^{\dag}({A^{\dag}_{K}})}$ (resp. $\mathbf{\textrm{\bf Conn}\hat{ }(\hat{A}_{K})}$) la cat\'egorie des $A^{\dag}_{K}$-modules (resp. $\hat{A}_{K}$-modules) projectifs de type fini $M$ (resp. $\mathcal{M})$ \`a connexion int\'egrable dont la s\'erie de Taylor converge sur un voisinage strict du tube de la diagonale dans $\mathcal{X}^{\mbox{an}}_{K} \times \mathcal{X}^{\mbox{an}}_{K}$ (resp. dont la s\'erie de Taylor converge sur le tube de la diagonale dans $\mathcal{X}^{\mbox{an}}_{K} \times \mathcal{X}^{\mbox{an}}_{K}).$ On note $\mathbf{F^{a}\mbox{-}Isoc^{\dag} (X/K)}$ (resp $\mathbf{F^{a}\mbox{-}Isoc (X/K)}$) la cat\'egorie de $F$-isocristaux surconvergents (resp des  $F$-isocristaux convergents ) sur $X$, d\'efinie par Berthelot [B 3]: une d\'efinition rapide est donn\'ee par Crew [Cr, \S1]; cf aussi [LS] pour plus de d\'etails.\\

\noindent On utilisera un exposant $(\ )^{\circ}$ pour sp\'ecifier les sous-cat\'egories des objets unit\'es (i.e. tels qu'en tout point g\'eom\'etrique les pentes du Frobenius sont nulles), $F^{a}$-Mod$(A^{\dag}_{K})^{\circ}$, $F^{a}$-Mod$(\hat{A}_{K})^{\circ}$, $F^{a}$-Conn$^{\dag}(A^{\dag}_{K})^{\circ}$, $F^{a}$-Conn$(\hat{A}_{K})^{\circ}$, ... ou la restriction de foncteurs aux objets unit\'es.\\

 On dispose de foncteurs naturels rendant commutatif le diagramme [Et 3, I, \S\ 5]

$$
\xymatrix{
F^{a}\mbox{-Isoc}^{\dag} (X/K) \ar[r]^{\Gamma^{\dag}} \ar[d]^{\mathcal{F}} & F^{a}\mbox{-Conn}^{\dag}(A^{\dag}_{K}) \ar[r] \ar[d]^{\mathcal{G}} & F^{a}\mbox{-Modloc}(A^{\dag}_{K}) \ar[d]^{\mathcal{H}}\\
F^{a}\mbox{-Isoc} (X/K) \ar[r]^{\hat{\Gamma}}& F^{a}\mbox{-Conn}{(\hat{A}}_{K}) \ar[r] & F^{a}\mbox{-Modloc} {(\hat{A}}_{K}), 
}
$$

\noindent o\`u $\Gamma^{\dag} := \Gamma(\mathcal{X}^{\mbox{an}}_{K}, -)$ est une \'equivalence de cat\'egories [B 3, cor (2.5.8)], $\hat{\Gamma} := \Gamma(]X[, -)$ est pleinement fid\`ele [O\ 2, 2.15, 2.23)] et [Et 3, I, (5.2.2)] et le foncteur $\mathcal{G}$ (resp. $\mathcal{H}$) envoie un ${A^{\dag}_{K}}$-module projectif de type fini $M$ sur son s\'epar\'e compl\'et\'e $p$-adique $\mathcal{M} = M \otimes_{A^{\dag}_{K}} \hat{A}_{K} : \mathcal{G}$ et $\mathcal{H}$ sont fid\`eles [Bour, A II, \S\ 5, $\mbox{n}^\circ\  3$ , prop 7] et [Bour, AC\ I, \S\ 3, $\mbox{n}^\circ\  5$, prop 9\ c)].\\

La restriction $\mathcal{F}^\circ$ de $\mathcal{F}$ \`a $F^a$-Isoc$^{\dag}(X / K)^\circ$ est un foncteur pleinement fid\`ele [Et 3, th\'eo 5]

\centerline{ $\mathcal{F}^{\circ} : F^a$-Isoc$^{\dag}  (X / K)^{\circ}\ \longrightarrow\ F^a$-Isoc$ (X / K)^{\circ}$ ;}
\noindent et m\^eme $\mathcal{F}$ est pleinement fid\`ele [Ked, theo 1.1].\\

 Nous allons montrer en 1.2 ci-apr\`es que $\hat{\Gamma}$ est en fait une \'equivalence de cat\'egories: ainsi $\mathcal{G}$ sera pleinement fid\`ele.

\subsection*{1.2. Des \'equivalences de cat\'egories}

Avec les notations de 1.1 nous allons montrer que la donn\'ee d'un $F^a$-isocristal convergent sur $X$ \'equivaut \`a celle d'un $\hat{A}_{K}$-module projectif de type fini $\mathcal{M}$, muni d'une connexion int\'egrable

$$\nabla : \mathcal{M}\  \longrightarrow\  \mathcal{M}\ \otimes_{\hat{A}_{K}}\ \Omega^1_{{\hat{A}_{K}}}$$

\noindent et d'un isomorphisme horizontal

$$\phi : (\mathcal{M}^{\sigma}, \nabla^{\sigma})\ \tilde{\longrightarrow}\ (\mathcal{M}, \nabla), $$

\noindent o\`u $(\mathcal{M}^{\sigma}, \nabla^{\sigma})$ provient de $(\mathcal{M}, \nabla)$ en \'etendant les scalaires par $F_{\hat{A}_{K}}$.

\vskip 3mm
\noindent \textbf{Proposition (1.2.1)}. \textit{Avec les notations ci-dessus, le foncteur $\Gamma(]X[, -)$ induit une \'equivalence entre  }

\begin{enumerate}
\item[(i)] \textit{La cat\'egorie des $\mathcal{O}_{]X[}$-modules coh\'erents (resp. et localement libres), et celles des $\hat{A}_{K}$-modules de type fini (resp. et projectifs) ; }
\item[(ii)] \textit{La cat\'egorie des  $\mathcal{O}_{]X[}$-modules coh\'erents \`a connexion int\'egrable (resp. des isocristaux convergents sur $X$), et celle des $\hat{A}_{K}$-modules projectifs de type fini munis d'une connexion int\'egrable (resp. et dont la s\'erie de Taylor converge sur $]X[_{\mathcal{X}^2}$).}

\end{enumerate}

\noindent \textbf{Remarque}. Berthelot a fourni une description analogue pour les $j^{\dag} \mathcal{O}_{]X[-}$-modules coh\'erents [B 3, (2.5.2)].\\

\noindent \textit{D\'emonstration}. La d\'emonstration est semblable \`a celle de loc. cit. Remarquons simplement que $]X[$ \'etant affino¬\"{\i}de, la donn\'ee d'un $\mathcal{O}_{]X[}$-module coh\'erent $\mathcal{E}$ \'equivaut \`a celle du $\hat{A}_{K}$-module de type fini $\mathcal{M} = \Gamma(]X[, \mathcal{E})$ [B-G-R, 9.4.2]. De m\^eme l'assertion ``projectif" en (ii) r\'esulte du fait que $K$ est de caract\'eristique 0 [B 3, (2.2.3) (ii)] [P, 10.3.1]. $\square$

\vskip 3mm
\noindent \textbf{Th\'eor\`eme (1.2.2)}. \textit{Avec les notations de 1.1, soient $\mathcal{M}$ un $\hat{A}_{K}$-module de type fini, muni d'une connexion int\'egrable $\nabla$, et $(\mathcal{M}^{\sigma}, \nabla^{\sigma})$ le module \`a connexion int\'egrable d\'eduit de $(\mathcal{M}, \nabla)$ par l'extension des scalaires $F_{\hat{A}_{K}}$. On suppose qu'il existe un isomorphisme horizontal }

$$\phi : (\mathcal{M}^{\sigma}, \nabla^{\sigma})\ \tilde{\longrightarrow}\  (\mathcal{M}, \nabla).$$
\textit{Alors, si $(\mathcal{E, \nabla})$ est le $\mathcal{O}_{]X[}$-module correspondant \`a $(\mathcal{M}, \nabla)$ par l'\'equivalence de (1.2.1), la connexion $\nabla$ de $\mathcal{E}$ est convergente.}\\

\noindent \textit{D\'emonstration}. La preuve suit celle de [B 3, (2.5.7)]. Comme l'assertion est locale sur $X$ [B 3, (2.2.11)] on peut supposer $\Omega^1_{X}$ libre de base $d\overline{z}_{1},...,d\overline{z}_{m}$; de m\^eme, comme $\mathcal{M}$ est n\'ecessairement projectif d'apr\`es [P, lemme 10.3.1], on peut supposer $\mathcal{M}$ libre de base $e_{1},...,e_{r}$. Si $z_{1},...,z_{m}  \in A$ rel\`event $\overline{z}_{1},...\overline{z}_{m}$, $\Omega^1_{\hat{A}}$ est un $\hat{A}$-module libre de base $dz_{1},...,dz_{m}$; soient $\partial_{1},...,\partial_{m}$ les d\'erivations correspondantes. Pour tout $i$ et tout $k$, posons $\partial_{i}\ e_{k} = \displaystyle \mathop{\Sigma}_{j}\  b_{ijk}\ e_{j}$, et soit $B_{i} \in M_{r}(\hat{A}_{K})$ la matrice des $b_{i,j,k}$ ; la  connexion $\nabla$ est d\'etermin\'ee par les $B_{i}$. En utilisant l'isomorphisme horizontal $\phi$, on peut supposer, comme dans la preuve de [B 3, (2.5.7)] que les matrices $B_{i}$ sont \`a coefficients dans $\hat{A}$.\\

Soit $\eta_{0} \in\ ]0,p^{-(1/p-1)} [$. Pour tout $e \in \mathcal{M} = \Gamma(]X[, \mathcal{E}),\  \underline{\partial}\ ^{\underline{k}}\  e$ est \`a coeffficients dans $\hat{A}$, donc $\parallel \underline{\partial}\ ^{\underline{k}}\  e \parallel\  \leqslant 1$, et comme $\eta_{0} < p^{-(1/p-1)}$ on en d\'eduit\\

\noindent (1.2.2.1) \qquad  $\parallel \frac{1}{\underline{k}!}\ \underline{\partial}\ ^{\underline{k}}\  e \parallel\ \eta_{0}^{\mid \underline{k} \mid}\ \longrightarrow 0$ quand $\mid \underline{k} \mid\ \longrightarrow + \infty .$\\

\noindent Posons $\zeta_{i} = 1 \otimes z_{i} - z_{i} \otimes 1\ ,\ \tau_{j} = 1 \otimes T_{j} - T_{j} \otimes 1.$\\

\noindent Pour $\eta < 1$, soient\\

$U = \{x \in \mathcal{X}^{an}_{K} \times \mathcal{X}^{an}_{K} /\  \forall j,~ \mid(1 \otimes T_{j})(x) \mid\  \leqslant 1, ~Ê\mid (T_{j} \otimes 1)(x) \mid\  \leqslant 1\} ,$\\

$W_{\eta} = \{x \in \mathcal{X}^{an}_{K} \times \mathcal{X}^{an}_{K} /\ \forall j,  \mid \tau_{j}(x) \mid\  \leqslant \eta \},$\\

$W'_{\eta} = \{x \in \mathcal{X}^{an}_{K} \times \mathcal{X}^{an}_{K} /\ \forall i,  \mid \zeta_{i}(x) \mid\  \leqslant \eta \},$\\

\noindent et $V_{\eta} = W_{\eta} \cap U$. Pour toute suite croissante $\underline{\eta}$ de limite 1, l'ouvert $V_{\underline{\eta}} :=  \displaystyle\mathop{\cup}_{{n}}\ V_{{\eta}_{n}}$  est \'egal \`a $] X [_{\mathcal{X}^2}$, d'apr\`es l'exemple de [B 3, (1.3.10)]. Nous allons construire une suite $\underline{\eta}$ telle qu'il existe sur l'ouvert $V_{\underline{\eta}} = ] X [_{\mathcal{X}^2}$ un isomorphisme $\epsilon : p^{\ast}_{2}\ \mathcal{E} \tilde{\longrightarrow} p^{\ast}_{1}\  \mathcal{E}$ induisant sur les voisinages infinit\'esimaux de la diagonale les isomorphismes $\epsilon_{n}$ d\'efinis par $\nabla$.\\

Puisque les $T_{j}$ engendrent $A$, il existe des relations $\zeta_{i} = \displaystyle \mathop{\Sigma}_{j}\ \beta_{ij}\  \tau_{j}$, avec $\beta_{ij} \in\ A \otimes_{\mathcal{V}} A$ ; d'o\`u $\parallel \zeta_{i} \parallel\  \leqslant \mbox{Sup}_{j} \parallel \tau_{j} \parallel$ , la norme \'etant la norme spectrale sur $U$. Par suite on a $V_{\eta_{0}} \subset W'_{\eta_{0}} \cap U$. On proc\`ede alors comme dans la d\'emonstration de [B 3, (2.2.13)] pour d\'efinir sur $V_{\eta_{0}}$ un isomorphisme $\varepsilon : p^{\ast}_{2}\ \mathcal{E}\  \tilde{\longrightarrow}\ 
p^{\ast}_{1}\ \mathcal{E}$ induisant sur les voisinages infinit\'esimaux de la diagonale les isomorphismes $\varepsilon_{n}$ d\'efinis par $\nabla$, en posant

$$\varepsilon(p^{\ast}_{2}(e)) = \displaystyle \mathop{\Sigma}_{\underline{k}}\ \frac{1}{\underline{k}!}\ \underline{\partial}^{\underline{k}}\ e \otimes \underline{\zeta}^{\underline{k}} ,$$

\noindent la s\'erie convergeant dans $\Gamma(V_{\eta_{0}}, p^{\ast}_{1} (\mathcal{E}))$ gr\^ace \`a (1.2.2.1).  On va utiliser ensuite l'action de Frobenius pour prolonger l'isomorphisme $\varepsilon$ de $V_{\eta_{0}}$ \`a $V_{\underline{\eta}}$, pour une suite $\underline{\eta}$ convenable : tout d'abord, on peut supposer d'apr\`es (1.2.1) qu'il existe un morphisme $F_{K} : \ ]X[ \ = \mbox{Spm} (\hat{A}_{K}) \longrightarrow\  ]X [$ tel que l'homomorphisme $\Gamma(]X[, \mathcal{O}_{]X[}) \longrightarrow\ \Gamma(]X[, \mathcal{O}_{]X[})$ induit par $F_{K}$ soit \'egal \`a $F_{\hat{A}_{K}} \otimes\  \sigma$. Le Frobenius $\phi$ de $\mathcal{E}$ est alors d\'efini par un isomorphisme $\phi : F^{\ast}_{K}\ \mathcal{E}\  \tilde{\longrightarrow}\   \mathcal{E}$ sur $]X[$, ce qui fournit des isomorphismes

$$\phi_{i} = p^{\ast}_{i}(\phi) : (F_{K} \times F_{K})^{\ast}\  (p^{\ast}_{i}\  \mathcal{E})\ \tilde{\longrightarrow}\ p^{\ast}_{i}\ \mathcal{E}\  \mbox{sur}\  V_{\eta_{0}} .$$

\noindent On d\'efinit une suite croissante $\underline{\eta}$ de limite 1 en posant

$$\eta_{n+1} = \mbox{min}( \vert\pi\vert^{-1}\eta_{n}    ,\eta^{1/p^{a}}_{n})\ .$$

\noindent Montrons que $(F_{K} \times F_{K}) (V_{\eta_{n+1}}) \subset V_{\eta_{n}}$. Posons

$$F_{\hat{A}}(T_{j}) = T^{p^{a}}_{j} + \pi\ a_{j}, $$

\noindent avec $a_{j} \in \hat{A}$. Pour $x \in U$ on a
$$
\begin{array}{cccc}
\mid (T_{j} \otimes 1) ((F_{K} \times F_{K})(x)) \mid &=&\  \mid (F_{\hat{A}}(T_{j}) \otimes 1) (x) \mid& \\
&=&\ \mid ((T^{p^ a}_{j} + \pi\ a_{j}) \otimes 1)(x) \mid & \leqslant 1,
\end{array}
 $$

\noindent de m\^eme pour $1 \otimes T_{j}$, de sorte que  $(F_{K} \times F_{K})(x) \in U$. D'autre part, si $J = \mbox{Ker} (A \otimes_{\mathcal{V}} A \rightarrow A)$, on a dans $(A \otimes_{\mathcal{V}} A)^{\wedge}Ê$  la relation\\

\qquad  \qquad   ${(F_{\hat{A}} \times F_{\hat{A}}) (\tau_{j}) = 1 \otimes (T^{p^ a}_{j} +  \pi\  a_{j})\ -\ (T^{p^ a}_{j} +  \pi\  a_{j})\ \otimes\ 1} $

$\qquad \qquad \qquad \qquad \quad \quad ~ = \tau^{p^ a}_{j} + \pi\ \alpha_{j} ,$\\

\noindent avec $\alpha_{j} \in J(A \otimes_{\mathcal{V}} A)^{\wedge} $. Dans $(A \otimes_{\mathcal{V}} A)^{\wedge}$ on peut \'ecrire $\alpha_{j}$ sous la forme $\alpha_{j} = \displaystyle \mathop{\Sigma}_{i}\  \gamma_{ij}\  \tau_{i}$, avec $\| \gamma_{ij} \| \leqslant 1$. Alors, pour $x \in V_{\eta_{n+1}}$, on obtient

$$| \pi\ \alpha_{j}(x) | \leqslant | \pi |\  \eta_{n+1} \leqslant \eta_{n} $$

\noindent et

$$| \tau^{p^ a}_{j} (x) | \leqslant \eta^{p^a}_{n+1} \leqslant \eta_{n}.$$

\noindent Par cons\'equent on a bien

$$(F_{K} \times F_{K}) (V_{\eta_{n+1}}) \subset V_{\eta_{n}}. $$

\vskip 2mm
Supposons construit sur $V_{n} := \displaystyle \mathop{\cup}_{i \leqslant n} V_{\eta_{i}}$ un isomorphisme $\varepsilon^{(n)} : p^{\ast}_{2}\  \mathcal{E} \tilde{\longrightarrow}\  p^{\ast}_{1}\ \mathcal{E}$ ; en utilisant l'action du Frobenius, on conclut alors, comme dans la preuve de [B 3, (2.5.7)], \`a l'existence d'un isomorphisme sur $V = \displaystyle \mathop{\cup}_{i } V_{\eta_{i}}  =\  ]X[_{\mathcal{X}^2}$, d'o\`u la convergence de $\nabla .~ \square$

\vskip 3mm
\noindent \textbf{Corollaire (1.2.3)}. \textit{Soient $X = \mbox{Spec}\  A_{0}$ un $k$-sch\'ema affine et lisse, $A$ une $\mathcal{V}$-alg\`ebre lisse relevant $A_{0}$, $\hat{A}$ le s\'epar\'e compl\'et\'e $p$-adique de $A$, $F_{\hat{A}} : \hat{A} \rightarrow \hat{A}$ un rel\`evement de la puissance $a^{\mbox{i\`eme}}$ de l'endomorphisme de Frobenius de $A_{0}$ au-dessus de $\sigma$. Alors la cat\'egorie des $F^{a}$-isocristaux convergents sur $X$ est \'equivalente \`a la cat\'egorie des $\hat{A}_{K}$-modules (n\'ecessairement projectifs) de type fini $\mathcal{M}$, munis d'une connexion int\'egrable $\nabla$ et d'un isomorphisme horizontal $\phi : (\mathcal{M}^{\sigma}, \nabla^{\sigma})\ \tilde{\longrightarrow}\ (\mathcal{M}, \nabla)$. }\\

\noindent \textit{D\'emonstration}. D'apr\`es la proposition (1.2.1) la cat\'egorie des isocristaux convergents sur $X$ est \'equivalente \`a celle des $\hat{A}_{K}$-modules projectifs de type fini, munis d'une connexion convergente $\nabla$. De plus cette \'equivalence est fonctorielle par rapport \`a $\hat{A}$ par des arguments analogues \`a [B 3, (2.5.6)]. On conclut par le th\'eor\`eme (1.2.2).\ $\square$\\

\newpage
\section*{2. $F$-isocristaux convergents sur un sch\'ema lisse formellement relevable}

Nous allons g\'en\'eraliser au cas relevable l'\'equivalence de cat\'egories du corollaire (1.2.3) pr\'ec\'edent.
\subsection*{2.1. Espaces rigides associ\'es aux sch\'emas formels}

Dans ce \S\ 2 on se donne $f : X \rightarrow \mbox{Spec}\  k$ un $k$-sch\'ema lisse tel qu'il existe un $\mathcal{V}$-sch\'ema formel lisse $h : \mathcal{X} \rightarrow Spf\  \mathcal{V}$ relevant $f$. Par une construction de Raynaud on sait associer \`a $\mathcal{X}$ (resp. \`a $h$) un espace rigide analytique not\'e $\mathcal{X}_{K}$ [Bo-L¬\"u 1 et 2] [B 3] (resp. un morphisme $h_{K} : \mathcal{X}_{K}
\rightarrow \mbox{Spm}\ K)$. Si $h$ est propre alors $h_{K}$ est propre [L¬\"u]. L'espace $\mathcal{X}_{K}$ est muni d'une topologie de Grothendieck [B 3, (0.1.2), (0.2)] et d'un faisceau d'anneaux $\mathcal{O}_{\mathcal{X}_{K}}$ : nous dirons que $(\mathcal{X}_{K}, \mathcal{O}_{\mathcal{X}_{K}})$ est un $G$-espace annel\'e [B-G-R, 9.3.1].

\vskip 3mm
\noindent \textbf{Proposition (2.1.1)}. \textit{Sous les hypoth\`eses 2.1 il existe un $\mathcal{V}$-sch\'ema formel lisse $h' : \mathcal{X}' \rightarrow \mbox{Spf}\  \mathcal{V}$, un $\mathcal{V}$-isomorphisme $\mathcal{X'} \tilde{\longrightarrow} \mathcal{X}$ et des recouvrements par des ouverts lisses $\mathcal{X}' = \displaystyle \mathop{\cup}_{\alpha}\ \mbox{Spf}\ \hat{A}_{\alpha}$, $X = \displaystyle \mathop{\cup}_{\alpha}\ Spec\ A_{\alpha,0}$, o\`u les $A_{\alpha}$ sont des $\mathcal{V}$-alg\`ebres lisses et $A_{\alpha,0} := A_{\alpha} / \pi A_{\alpha}$.}\\

\noindent \textit{D\'emonstration}. Dire que $\mathcal{X}$ est un $\mathcal{V}$-sch\'ema formel lisse signifie qu'il existe un recouvrement $\mathcal{X} = \displaystyle \mathop{\cup}_{\alpha}\  \mbox{Spf}\  \mathcal{B}_{\alpha}$, o\`u les $\mathcal{B}_{\alpha}$ sont des $\mathcal{V}$-alg\`ebres plates s\'epar\'ees et compl\`etes pour la topologie $\pi$-adique et formellement lisses pour les topologies discr\`etes sur $\mathcal{V}$ et $\mathcal{B}_{\alpha} $ : les  $\mathcal{B}_{\alpha}$ sont donc des $\mathcal{V}$-alg\`ebres formellement lisses pour les topologies $\pi$-adiques sur $\mathcal{V}$ et $\mathcal{B}_{\alpha}$ [EGA O$_{IV}$, (19.3.1)]. Pour tout $\alpha$ notons $A_{\alpha,0} : = \mathcal{B}_{\alpha}\ /\ \pi\ \mathcal{B}_{\alpha}$, et $A_{\alpha}$ une $\mathcal{V}$-alg\`ebre lisse relevant $A_{\alpha,0}$ [E$\ell$, th\'eo 6] : d'apr\`es [Et 2, cor 1 du th\'eo 4] il existe, pour tout $\alpha$, un $\mathcal{V}$-isomorphisme

$$\hat{A}_{\alpha}\  \tilde{\longleftarrow}\ \mathcal{B}_{\alpha}, \mbox{o\`u}\  \hat{A}_{\alpha} := 
\displaystyle \mathop{\mbox{lim}}_{{\leftarrow} \atop{n}} A_{\alpha}\ /\  \pi^n\ A_{\alpha}.$$

\noindent Posons
$$P_{\alpha} = \mbox{Spf}\ \mathcal{B}_{\alpha}, \ P'_{\alpha} = \mbox{Spf}\ \hat{A}_{\alpha},\  P'_{\alpha \beta}\  := P'_{\alpha}\ \times_{P_{\alpha}} (P_{\alpha}\  \cap P_{\beta});$$
\noindent d\'esignons respectivement par $\psi_{\alpha}$ le $\mathcal{V}$-isomorphisme $\psi_{\alpha} : P'_{\alpha}\ \tilde{\rightarrow}\ P_{\alpha}$, et par $\psi_{{\alpha}_{\beta}}$ le $\mathcal{V}$-isomorphisme 
$$\psi_{{\alpha}_{\beta}}:P'_{\alpha \beta}\ \tilde{\rightarrow}\ P_{\alpha}\ \cap\ P_{\beta}$$
 \noindent d\'eduit de $\psi_{\alpha}$ par le changement de base $P_{\alpha}\ \cap\ P_{\beta}\ \hookrightarrow\ P_{\alpha}.$ Le $\mathcal{V}$-isomorphisme

$$\varphi_{{\alpha}_{\beta}} :  = \psi^{-1}_{\beta \alpha}\ \circ\  \psi_{{\alpha}_{\beta}} : P'_{\alpha \beta}\ \tilde{\rightarrow}\ P'_{\beta \alpha} $$

\noindent induit des $\mathcal{V}$-isomorphismes

$$\varphi_{{\alpha \beta \gamma}} : P'_{\alpha \beta}\ \cap\ P'_{{\alpha \gamma}}\ \tilde{\rightarrow}\  P'_{\beta \alpha}\ \cap\ P'_{\beta \gamma}$$

\noindent et on a les identit\'es

$$\varphi_{\alpha \beta}\ \circ\ \varphi_{\beta \alpha} = Id\ ,\ \varphi_{\alpha \alpha} = Id, $$

$$\varphi_{\alpha \beta \gamma} =  \varphi_{\gamma \beta \alpha} \circ\ \varphi_{\alpha \gamma \beta}\ ,$$

 \noindent gr\^ace au fait que les $P_{\alpha}$ se recollent pour former $\mathcal{X}$. Par cons\'equent on peut recoller les $P'_{\alpha}$ le long des $P'_{\alpha \beta}$ : le sch\'ema formel ainsi obtenu est le sch\'ema formel $\mathcal{X}'$ cherch\'e. $\square$

\vskip 3mm
\noindent \textbf{Corollaire (2.1.2)}. \textit{Sous les hypoth\`eses 2.1 le morphisme d'espaces rigides $h_{K} : \mathcal{X}_{K}\ \rightarrow\ \mbox{Spm}\ K$ est lisse. }\\

\noindent \textit{D\'emonstration}. Compte-tenu de l'isomorphisme $\mathcal{X}'\ \tilde{\rightarrow}\ \mathcal{X}$ de (2.1.1) il suffit d'appliquer le crit\`ere jacobien [B 3, (0.1.11)]. $\square$

\vskip 3mm
\noindent \textbf{Corollaire (2.1.3)}. \textit{Soient $S$ un $k$-sch\'ema lisse et $f : X\ \rightarrow\ S$ un $k$-morphisme lisse et supposons donn\'es un $\mathcal{V}$-sch\'ema formel lisse $\mathcal{S}$ et un $\mathcal{V}$-morphisme lisse $h : \mathcal{X}\ \rightarrow\ \mathcal{S}$ de sch\'emas formels relevant $f$. Alors }

\begin{enumerate}
\item[(i)] \textit{Il existe un $\mathcal{V}$-sch\'ema formel lisse $\mathcal{S}' = \displaystyle \mathop{\cup}_{\alpha}\ \mathcal{S}'_{\alpha}$\ , $\mathcal{S}'_{\alpha} = \mbox{Spf}\ \hat{A}_{\alpha}$ o\`u les $A_{\alpha}$ sont des $\mathcal{V}$-alg\`ebres lisses et un $\mathcal{V}$-isomorphisme $\mathcal{S}'\ \tilde{\rightarrow}\ \mathcal{S}$.}
\item[(ii)] \textit{Si l'on pose $\mathcal{X}' : = \mathcal{X} \times_{\mathcal{S}}\ \mathcal{S}'$, $\mathcal{X}'_{\alpha} : = \mathcal{X}' \times_{\mathcal{S}} \mathcal{S}'_{\alpha} $, il existe de plus un $\mathcal{S}'$-sch\'ema formel lisse $\mathcal{X}''$ et un $\mathcal{S}'$-isomorphisme $\mathcal{X}''\ \tilde{\rightarrow}\ \mathcal{X}'$ tel que}

$$\mathcal{X}'' = \displaystyle \mathop{\cup}_{\alpha} \mathcal{X}''_{\alpha}\ , \mbox{\textit{{o\`u}}}\  \mathcal{X}''_{\alpha} = \displaystyle \mathop{\cup}_{\beta}\ \mathcal{X}''_{\alpha, \beta} = \displaystyle \mathop{\cup}_{\beta}\  \mbox{Spf}\  \hat{B}_{\alpha, \beta}\ ,$$ 

\textit{les $B_{\alpha, \beta}$  \'etant des $\hat{A}_{\alpha}$-alg\`ebres lisses (resp. des $A^{\dag}_{\alpha}$-alg\`ebres lisses)}

\item[(iii)] \textit{Avec les notations du (ii) il existe aussi un $\mathcal{V}$-sch\'ema formel lisse $\mathcal{X}'''$ et un $\mathcal{V}$-isomorphisme $\mathcal{X}'''\ \tilde{\rightarrow}\ \mathcal{X}''$ tel que  }

$$\mathcal{X}'''\ = \displaystyle \mathop{\cup}_{\alpha}\ \mathcal{X}'''_{\alpha}\ ,\ \mbox{o\`u}\  \mathcal{X}'''_{\alpha} = \displaystyle \mathop{\cup}_{\beta}\ \mbox{Spf}\ \hat{C}_{\alpha, \beta}$$

\noindent \textit{et $ C_{\alpha,\beta}$ est une $\mathcal{V}$-alg\`ebre lisse munie d'un $\mathcal{V}$-isomorphisme}

$$\hat{B}_{\alpha, \beta} \simeq\ \hat{C}_{\alpha, \beta}.$$

\end{enumerate}

\noindent \textit{D\'emonstration}.
 La proposition (2.1.1.) fournit le (i).\\
 Pour (ii) et (iii) on utilise encore [Et 2, th\'eo 4 et son cor 1]: on peut recoller les $\mathcal{X}''_{\alpha,\beta}$ (resp. les $\mathcal{X}''_{\alpha}$, resp. les $\mathcal{X}'''_{\alpha})$ gr\^ace \`a l'existence globale de $\mathcal{X}'_{\alpha}$ (resp. de $\mathcal{X}'$, resp. de $\mathcal{X}'').\  \square$

\vskip 3mm
\noindent \textbf{Corollaire (2.1.4)}. \textit{Sous les hypoth\`eses (2.1.3) le morphisme d'espaces rigides $h_{K} : \mathcal{X}_{K} \rightarrow\mathcal{S}_{K}$ est lisse. Si de plus $h$ est propre alors $h_{K}$ est propre. }\\

\noindent \textit{D\'emonstration}. La lissit\'e de $h_{K}$ r\'esulte de (2.1.3) (ii) et du crit\`ere jacobien [B 3, (0.1.11]. \\
La d\'efinition d'un morphisme propre d'espaces rigides est donn\'ee dans [L¬\"u, 2.4] : la propret\'e de $h$ entra\^{\i}ne celle de $h_{K}$ [L¬\"u, theo 3.1].\ $\square$\\

\subsection*{2.2. F-isocristaux convergents et $\mathcal{O}_{\mathcal{X}_{K}}$-modules }
\vskip 2mm
Rappelons que les hypoth\`eses (2.1) sont satisfaites.\\

La donn\'ee de $\mathcal{E} \in \mbox{Isoc}(X/K)$ \'equivaut  \`a celle d'un $\mathcal{O}_{\mathcal{X}_{K}}$-module localement libre de type fini $\mathcal{E}_{\mathcal{X}} $ [B 3, (2.3.2), (2.2.3) (ii)] muni d'une connexion $\nabla$ relativement \`a $K$, int\'egrable et convergente. D'apr\`es (2.1.1) on a un $\mathcal{V}$-isomorphisme $\mathcal{X} \simeq \displaystyle \mathop{\cup}_{\alpha}\ \mbox{Spf}\ \hat{A}_{\alpha} = : \displaystyle \mathop{\cup}_{\alpha}\ \mathcal{X}_{\alpha}$ o\`u les $A
_{\alpha}$ sont des $\mathcal{V}$-alg\`ebres lisses ; si $F_{\alpha_{1}}, F_{\alpha_{2}}: \mathcal{X}_{\alpha} \rightarrow \mathcal{X}_{\alpha}$ sont deux rel\`evements de la puissance $a^{\mbox{i\`eme}}$ du Frobenius absolu de $X_{\alpha} = \mathcal{X}_{\alpha}\ \mbox{mod}\ \pi$, alors on a un isomorphisme canonique [B 3, (2.2.17)]

$$F^{\ast}_{\alpha_{1}}(\mathcal{E}_{\mathcal{X}_{\alpha}})\  \tilde{\longrightarrow}\ F^{\ast}_{\alpha_{2} }(\mathcal{E}_{\mathcal{X}_{\alpha}})\ , 
\leqno{(2.2.0)}
$$

\noindent o\`u $(\mathcal{E}_{\mathcal{X}_{\alpha}}, \nabla_{\alpha})$ est la restriction de $(\mathcal{E}_{\mathcal{X}}, \nabla)$ \`a $\mathcal{X}_{\alpha K} = \mbox{Spm}\ (\hat{A}_{\alpha K})$. \\

\noindent Supposons fix\'e pour chaque $\alpha$ un rel\`evement $F_{\alpha} : \mathcal{X}_{\alpha} \rightarrow \mathcal{X}_{\alpha}$, de la puissance $a^{\mbox{i\`eme}}$ du Frobenius absolu de $X_{\alpha} = \mathcal{X}_{\alpha}\  \mbox{mod}\  \pi$, et soit $\mathcal{E} \in F^a{\mbox{-}}\textrm{Isoc}(X/K)$. La structure de Frobenius sur $\mathcal{E}$ fournit pour tout $\alpha$ un isomorphisme [cor (1.2.3)]

$$
\phi_{\alpha} : (F^{\ast}_{\alpha}(\mathcal{E}_{\mathcal{X}_{\alpha}}), F^{\ast}_{\alpha}(\nabla_{\alpha}))  \displaystyle \mathop{\longrightarrow}^{\sim} (\mathcal{E}_{\mathcal{X}_{\alpha}}, \nabla_{\alpha}),
$$

\noindent avec compatibilit\'es \'evidentes quand $\alpha$ varie : d'apr\`es (1.2.3) la connaissance de $\mathcal{E}_{\mathcal{X}_{\alpha}}$ \'equivaut \`a celle de $\Gamma(\mathcal{X}_{\alpha K}, \mathcal{E}_{\mathcal{X}_{\alpha}})$ avec m\^emes donn\'ees.\\

Notons $\mathbf{F^a\textrm{\bf -Conn}(\mathcal{X}_{K})}$ la cat\'egorie des $\mathcal{O}_{\mathcal{X}_{K}}$-modules localement libres de type fini $\mathcal{M}$ munis d'une connexion $\nabla$ relativement \`a $K$, int\'egrable et d'une famille d'isomorphismes

$$
\phi_{\alpha} : (F^{\ast}_{\alpha}(\mathcal{M}_{\alpha}), F^{\ast}_{\alpha}(\nabla_{\alpha}))  \displaystyle \mathop{\longrightarrow}^{\sim} (\mathcal{M}_{\alpha}, \nabla_{\alpha}) := (\mathcal{M}, \nabla)_{|\mathcal{X}_{\alpha K}}\ ,
$$
\noindent compatibles aux isomorphismes de (2.2.0).
\vskip 3mm
\noindent \textbf{Th\'eor\`eme (2.2.1)}. \textit{ Avec les notations de (2.2) la fl\`eche naturelle
$$
F^{a}{\mbox{-}}Isoc(X/K) \longrightarrow F^{a}{\mbox{-}}Conn(\mathcal{X}_{K})
$$
$$ \mathcal{E}  \longmapsto \mathcal{E}_{\mathcal{X}}$$
est une \'equivalence de cat\'egories.}\\

\noindent \textit{D\'emonstration}. Nous venons de montrer que si $\mathcal{E} \in F^{a}$-$Isoc(X/K)$ alors $\mathcal{E}_{\mathcal{X}} \in F^{a}\mbox{-}Conn(\mathcal{X}_{K})$.\\

R\'eciproquement soit $\mathcal{M} \in F^{a}\mbox{-}Conn(\mathcal{X}_{K})$ ; puisque $\mathcal{M}$ est localement libre de type fini et que $\mathcal{O}_{\mathcal{X}_{K}}$ est coh\'erent [B 3, (2.1.9)], alors $\mathcal{M}$ est un $\mathcal{O}_{\mathcal{X}_{K}}$-module coh\'erent.  D'apr\`es (1.2.3) l'existence des isomorphismes horizontaux

$$
\phi_{\alpha} : (F^{\ast}_{\alpha}(\mathcal{M}_{\alpha}), F^{\ast}_{\alpha}(\nabla_{\alpha})) \displaystyle \mathop {\longrightarrow}^{\sim} (\mathcal{M}_{\alpha}, \nabla_{\alpha})
$$
prouve que la connexion est convergente, car la donn\'ee de $(\mathcal{M}_{\alpha}, \nabla_{\alpha}, \phi_{\alpha})$ \'equivaut \`a celle d'un \'el\'ement de $F^{a}\mbox{-}Isoc(X_{\alpha}/K)$ : puisque ces donn\'ees se recollent la connexion $\nabla$ est convergente [B 3, (2.2.11)], d'o\`u un \'el\'ement $\mathcal{E}$ de Isoc$(X/K)$ [B 3, (2.3.2)] tel que $\mathcal{E}_{\mathcal{X}} = \mathcal{M}$. Enfin les $\phi_{\alpha}$ d\'efinissent un isomorphisme de Isoc$(X/K)$

$$
\phi_{\mathcal{E}} : F^{\ast}_{\sigma}\  \mathcal{E} \displaystyle \mathop {\longrightarrow}^{\sim} \mathcal{E}
$$
qui est le Frobenius de $\mathcal{E}$ [B 3, (2.3.7)] et $(\phi_{\mathcal{E}})_{\mathcal{X}} = \phi_{\mathcal{M}}. \ \square $\\

Suite aux observations en fin de \S1.1 on en d\'eduit:
\vskip 3mm
\noindent \textbf{Corollaire (2.2.2)}. \textit{ Avec les notations de (2.2) et 1.1 le foncteur
$$
\mathcal{G}:F^{a}\mbox{-Conn}^{\dag}(A^{\dag}_{K})\longrightarrow F^{a}\mbox{-Conn}{(\hat{A}}_{K})
$$
est pleinement fid\`ele.}

\newpage
\section*{3. Rel\`evement de Teichm¬\"uller, r\'eduction, sp\'ecialisation et espaces rigides analytiques}

\subsection*{3.0.}
Avec les notations de l'introduction, on suppose dor\'enavant jusqu'\`a la fin de cet article que $k$ est un corps parfait de caract\'eristique $p>0,\ q = p^a$, et on fixe une cl\^oture alg\'ebrique $\overline{k}$ de $k$.\\

On suppose relev\'ee l'\'el\'evation \`a  la puissance $q$ sur $k$ (qui est un automorphisme puisque $k$ est parfait) en un automorphisme $\sigma$ de $\mathcal{V}$, tel que $\sigma(\pi) = \pi$: on note encore $\sigma$ son extension \`a $K$ et ce $\sigma$ est un automorphisme de $K$.

\subsection*{3.1. Rel\`evement de Teichm¬\"uller}

Soit $X_{0} = \textrm{Spec}\ A_{0}$ un $k$-sch\'ema lisse. Soient $k'\subset \overline{k}$ une extension finie de $k$ de degr\'e $d$ et $x_{0}\in X_{0}(k')$: le $k$-homomorphisme correspondant \`a $x_{0}$,  $s:A_{0}\rightarrow k(x_{0})=k'$, se factorise en $s(x_{0}):A_{0}(x_{0})=A_{0}\otimes_{k}k(x_{0})\rightarrow k(x_{0})$. Notons $W = W(k)$ (resp. $W(x_{0}) = W(k(x_{0}))$ l'anneau des vecteurs de Witt \`a coefficients dans $k$ (resp. $k(x_{0})$),
$$
\mathcal{V}(x_{0}) = W(x_{0}) \otimes_{W} \mathcal{V} \simeq W(x_{0}) [\pi]\ ,
$$
\noindent $K_{0} = \textrm{Frac}\ W$, $K_{0}(x_{0}) = \textrm{Frac}(W(x_{0})), K(x_{0}) = \textrm{Frac}(\mathcal{V}(x_{0}))$, $\sigma_{x_{0}}$ la puissance $p^a$ sur $k(x_{0})$, $\sigma_{W(x_{0})} = W(\sigma_{x_{0}})$ le rel\`evement canonique de $\sigma_{x_{0}}$ \`a $W(x_{0})$, $\sigma_{\mathcal{V}(x_{0})} = \sigma_{W(x_{0})} \otimes_{\sigma\vert W} \sigma$ et $\sigma_{K(x_{0})}$ (resp. $\sigma_{K_{0}(x_{0})})$ son extension naturelle \`a $K(x_{0})$ (resp. $K_{0}(x_{0}))$ d\'efinie par $\sigma_{K(x_{0})}(u/v) = \sigma_{\mathcal{V}(x_{0})}(u)/\sigma_{\mathcal{V}(x_{0})}(v)$ (resp. $\sigma_{K_{0}(x_{0})}(u/v) = \sigma_{W(x_{0})}(u)/\sigma_{W(x_{0})}v))$. Le morphisme $\sigma_{K(x_{0})}$ co¬\"{\i}ncide, d'apr\`es [Et 3, I.1.1] et [B-M 2, (1.2.7) (ii)], avec le morphisme $\sigma' : K' \rightarrow K'$ (au-dessus de $\sigma : K \rightarrow K)$ de [Et 3, I.1.1] pour $k' = k(x_{0})$.\\

Notons qu'il existe un unique polyn\^ome unitaire irr\'eductible $\overline{g}\in k[T]$ tel que l'on ait un isomorphisme de $k$-alg\`ebres $k(x_{0})\overset{\sim}{\rightarrow} k[T]/(\overline{g})$, et qu'\`a chaque rel\`evement (n\'ecessairement irr\'eductible) $g\in \mathcal{V}[T]$ de $\overline{g}$ est associ\'e un unique isomorphisme de $\mathcal{V}$-alg\`ebres relevant le pr\'ec\'edent [EGA IV, (18.3.2)]
$$
\mathcal{V}(x_{0})\overset{\sim}{\rightarrow} \mathcal{V}[T]/(g),
$$

\noindent en particulier pour $g=Teich(\overline{g})$, o\`u $Teich$ est le rel\`evement de Teichm\" uller d\'efini comme \'etant le morphisme compos\'e
$$
Teich: k[T]\rightarrow W[T]\hookrightarrow\mathcal{V}[T]
$$
$$
\overline{P}=\sum_{i=0}^{i=r}a_{i}T^{i}\mapsto Teich(\overline{P})=\sum_{i=0}^{i=r}b_{i}T^{i}
$$
 avec $b_{i}\in W(k)\overset{\iota}{\hookrightarrow} \mathcal{V}$ [Se, II, \S5, th\'eo 4] le rel\`evement de Teichm\" uller de $a_{i}$: d'apr\`es [Se, loc. cit.] et [I$\ell$, (0.1.3.16) et (0.1.3.18)], $Teich$ est un morphisme d'anneaux.\\

Toutes les consid\'erations pr\'ec\'edentes s'appliqueront aux points ferm\'es $x_{0} \in  \vert X_{0} \vert  :=  \{ \textrm{points ferm\'es de}\  X_{0} \}$; dans ce cas, le $k$-homomorphisme surjectif correspondant \`a $x_{0}$,  $s:A_{0}\twoheadrightarrow k(x_{0})$, se factorise en $s(x_{0}):A_{0}(x_{0})=A_{0}\otimes_{k}k(x_{0})\twoheadrightarrow k(x_{0})$ et $k(x_{0}) = A_{0}/\mathfrak{m}_{x_{0}}$ est une extension finie \'etale de $k$ de degr\'e $deg\ x_{0} = [k(x_{0}) : k]$.

Soit $A$ une $\mathcal{V}$-alg\`ebre lisse relevant $A_{0}$ et fixons une pr\'esentation $A = \mathcal{V}[T_{1},...,T_{n}]/(f_{1},...,f_{m})$. On d\'esigne par $\hat{A}$ (resp. $A^{\dag}$) le s\'epar\'e compl\'et\'e (resp. le compl\'et\'e faible) $\mathfrak{m}$-adique de $A$.  D'apr\`es [Et 3, (1.2.1)] il existe un carr\'e commutatif

$$
\xymatrix{
A^{\dag} \otimes_{\mathcal{V}} \mathcal{V}(x_{0}) =: A^{\dag}(x_{0})  \ar[rr]^{\qquad \quad F_{A^{\dag}(x_{0})}} & & A^{\dag}(x_{0}) \\
A^{\dag} \ar[u] \ar[rr]^{F_{A^{\dag}}}   & & A^{\dag}  \ar[u]
}
$$

\noindent au-dessus du carr\'e commutatif

$$
\xymatrix{
\mathcal{V}(x_{0})  \ar[r]^{\sigma_{\mathcal{V}(x_{0})}}  & \mathcal{V}(x_{0}) &\\
\mathcal{V} \ar[u] \ar[r]^{\sigma}  & \mathcal{V} \ar[u] &,
}
$$

\noindent o\`u $F_{A^{\dag}}$ est un rel\`evement \`a $A^{\dag}$ du Frobenius (puissance $q$) de $A_{0}$ ; d'o\`u un diagramme commutatif

$$
\xymatrix{
\mathcal{V}(x_{0})  \ar[r] \ar[d]_{\sigma_{\mathcal{V}(x_{0})}} & A^{\dag}(x_{0}) \ar[r] \ar[d]_{F_{A^{\dag}(x_{0})}}  & \hat{A}(x_{0}) := \hat{A} \otimes_{\mathcal{V}} \mathcal{V}(x_{0}) \ar[d]^{F_{\hat{A}(x_{0})}} &\\
\mathcal{V}(x_{0}) \ar[r]   & A^{\dag}(x_{0}) \ar[r] & \hat{A}(x_{0}) & ,
}
$$
\noindent avec pour fl\`eches horizontales les morphismes canoniques.
Nous allons montrer que le $k(x_{0})$-homomorphisme (surjectif si $x_{0}$ est un point ferm\'e)
$$s(x_{0}):A_{0}(x_{0})=A_{0}\otimes_{k}k(x_{0})\rightarrow k(x_{0})$$ 
se rel\`eve de mani\`ere canonique  en un $\mathcal{V}(x_{0})$-homomorphisme (surjectif si $x_{0}$ est un point ferm\'e)

$$
\tau_{\hat{A}(x_{0})} : \hat{A}(x_{0}) \rightarrow \mathcal{V}(x_{0})
$$

\noindent tel que le diagramme

$$
\xymatrix{
\hat{A}(x_{0})  \ar[r]^{\tau_{\hat{A}(x_{0})}} \ar[d]_{F_{\hat{A}(x_{0})}} & \mathcal{V}(x_{0}) \ar[d]^{\sigma_{\mathcal{V}(x_{0})}}\\
\hat{A}(x_{0})  \ar[r]^{\tau_{\hat{A}(x_{0})}}  & \mathcal{V}(x_{0})
}
$$

\noindent commute. Dans un premier temps nous traiterons le cas particulier d'un corps fini $k=\mathbb{F}_{q}$; dans un deuxi\`eme temps nous traiterons le cas g\'en\'eral d'un corps parfait de caract\'eristique $p>0$, en explicitant le morphisme $\tau_{\hat{A}(x_{0})}$.\\

\textbf{1\ier cas: cas particulier o\`u $k=\mathbb{F}_{q}$ est un corps fini.}\\

Par le th\'eor\`eme [Mo 2, theo 3.3] de Monsky et [Mo 2, def 5.2] (cf aussi [K1] dans le cas $\mathcal{V}=W$) le $k(x_{0})$-homomorphisme (surjectif si $x_{0}$ est un point ferm\'e)

$$s(x_{0}):A_{0}(x_{0})\rightarrow k(x_{0})$$

\noindent se rel\`eve de mani\`ere unique en un $\mathcal{V}(x_{0})$-homomorphisme (surjectif si $x_{0}$ est un point ferm\'e)
$$
\tau_{\hat{A}(x_{0})} : \hat{A}(x_{0}) \rightarrow \mathcal{V}(x_{0})
$$
tel que 
$$
\tau_{\hat{A}(x_{0})}\circ(F_{\hat{A}(x_{0})})^{deg\ x_{0}}=\tau_{\hat{A}(x_{0})}.
\leqno{(3.1.1)}
$$

\noindent Le morphisme $\tau_{\hat{A}(x_{0})}$ est appel\'e le \textbf{rel\`evement de Teichm\"{u}ller de $s(x_{0})$ (ou de $x_{0}$) relativement \`a $F_{\hat{A}(x_{0})}$}.\\
Or en notant $F_{A_{0}(x_{0})}$ (resp $F_{k(x_{0})}=\sigma_{k(x_{0})}=\sigma_{x_{0}}$) le Frobenius (\'el\'evation \`a la puissance $q$) de $A_{0}(x_{0})$ (resp de $k(x_{0})$), le diagramme suivant commute\\
$$
\xymatrix{
A_{0}(x_{0})  \ar[r]^{s(x_{0})} \ar[d]_{F_{A_{0}(x_{0})}} & k(x_{0}) \ar[d]^{\sigma_{x_{0}}}\\
A_{0}(x_{0})  \ar[r]^{s(x_{0})}  & k(x_{0}),
}
\leqno{(3.1.2)}
$$
et se rel\`eve en le diagramme

$$
\xymatrix{
\hat{A}(x_{0})  \ar[r]^{\tau_{\hat{A}(x_{0})}} \ar[d]_{F_{\hat{A}(x_{0})}} & \mathcal{V}(x_{0}) \ar[d]^{\sigma_{\mathcal{V}(x_{0})}}\\
\hat{A}(x_{0})  \ar[r]^{\tau_{\hat{A}(x_{0})}}  & \mathcal{V}(x_{0})\  ,
}
\leqno{(3.1.3)}
$$
dans lequel on v\'erifie que les deux morphismes $\tau_{\hat{A}(x_{0})}\circ F_{\hat{A}(x_{0})}$ et $\sigma_{\mathcal{V}(x_{0})}\circ \tau_{\hat{A}(x_{0})}$ sont des rel\`evements de Teichm\" uller de
$$
s(x_{0})\circ F_{A_{0}(x_{0})}=\sigma_{x_{0}}\circ s(x_{0})
$$
relativement \`a $F_{\hat{A}(x_{0})}$, puisque l'on a
$$
\tau_{\hat{A}(x_{0})}\circ F_{\hat{A}(x_{0})}\circ (F_{\hat{A}(x_{0})})^{deg\ x_{0}}=\tau_{\hat{A}(x_{0})}\circ F_{\hat{A}(x_{0})}
$$
\noindent et
$$
\sigma_{\mathcal{V}(x_{0})}\circ \tau_{\hat{A}(x_{0})}\circ (F_{\hat{A}(x_{0})})^{deg\ x_{0}}=\sigma_{\mathcal{V}(x_{0})}\circ \tau_{\hat{A}(x_{0})}\ .
$$
\noindent Par unicit\'e du rel\`evement de Teichm\" uller on en d\'eduit que (3.1.3) commute, i.e.
$$
\tau_{\hat{A}(x_{0})}\circ F_{\hat{A}(x_{0})}=\sigma_{\mathcal{V}(x_{0})}\circ \tau_{\hat{A}(x_{0})}\ .
$$

\textbf{2\ieme cas: cas g\'en\'eral d'un corps parfait $k$ de caract\'eristique $p> 0$.}

 Puisque $\hat{A}(x_{0})$ est sans $p$-torsion il existe d'apr\`es Illusie [I$\ell$, (0.1.3.16)] un unique homomorphisme d'anneaux
$$
s_{F_{\hat{A}(x_{0})}}: \hat{A}(x_{0})\rightarrow W(\hat{A}(x_{0}))
$$
\noindent qui soit section de la projection canonique $W(\hat{A}(x_{0}))\rightarrow\hat{A}(x_{0})$ et tel que
$$s_{F_{\hat{A}(x_{0})}}\circ F_{\hat{A}(x_{0})}= W(F_{\hat{A}(x_{0})})\circ s_{F_{\hat{A}(x_{0})}} \ , $$
\noindent o\`u $W(F_{\hat{A}(x_{0})}):W(\hat{A}(x_{0}))\rightarrow W(\hat{A}(x_{0}))$
est induit par fonctorialit\'e par $F_{\hat{A}(x_{0})}$ entre anneaux de vecteurs de Witt [I$\ell$, (0.1.3.19)].  De plus, pour $y\in \hat{A}(x_{0})$, $s_{F_{\hat{A}(x_{0})}}(y)$ est l'unique solution d'\'equations explicites [I$\ell$, (0.1.3.17)].\\
De m\^eme, puisque $A_{0}(x_{0})$et $k(x_{0})$ sont sans $p$-torsion, on d\'efinit $s_{F_{A_{0}(x_{0})}}$ et $s_{F_{k(x_{0})}}=s_{\sigma_{x_{0}}}$. Il est \`a remarquer que $s_{\sigma_{x_{0}}}$ est le rel\`evement de Teichm\" uller usuel pour un corps parfait de caract\'eristique $p>0$ [I$\ell$, (0.1.3.18)] et qu'il s'ins\`ere dans un diagramme commutatif
$$
\xymatrix@C=2cm @R=1cm
{
k(x_{0})\ar^{\sim}[r]\ar_{s_{\sigma_{x_{0}}}}[d]&k[T]/(\overline{g})\ar[d]^{Teich}\\
W(k(x_{0}))\ar^{\sim}[r]\ar_{\iota(x_{0})}@{^{(}->}[d]&W[T]/(Teich(\overline{g}))\ar@{^{(}->}[d]^{\iota}\\
\mathcal{V}(x_{0})\ar^{\sim}[r]&\mathcal{V}[T]/(Teich(\overline{g}))&,
}
\leqno{(3.1.4)}
$$
o\`u $\iota(x_{0}): W(k(x_{0}))\hookrightarrow \mathcal{V}(x_{0})$ est l'injection canonique [Se, II,\S5, th\'eo 4] et $\iota$ est induit par l'injection canonique $W\hookrightarrow \mathcal{V}$ [loc. cit.].\\

 En notant $W(u):W(R_{1})\rightarrow W(R_{2})$ le morphisme d'anneaux de vecteurs de Witt induit par un morphisme d'anneaux $R_{1}\rightarrow R_{2}$ on en d\'eduit un diagramme commutatif

$$
\xymatrix @C=2cm @R=1,5cm
{
&\hat{A}(x_{0})\ar[d]_{s_{F_{\hat{A}(x_{0})}}}\ar[r]^{proj} & A_{0}(x_{0})\ar[d]_{s_{F_{A_{0}(x_{0})}}}\ar[r]^{s(x_{0})}&k(x_{0})\ar[d]^{s_{\sigma_{x_{0}}}}\\
\hat{A}(x_{0})  \ar[r]^{s_{F_{\hat{A}(x_{0})}}} \ar[d]_{F_{\hat{A}(x_{0})}} &W(\hat{A}(x_{0})) \ar[r]^{W(proj)} \ar[d]_{W(F_{\hat{A}(x_{0})})}  & W(A_{0}(x_{0}))\ar[r]^{W(s(x_{0}))} \ar[d]_{W(F_{A_{0}(x_{0})})} &W(k(x_{0}))\ar[d]^{\sigma_{W(x_{0})}}\\
\hat{A}(x_{0})  \ar[r]_{s_{F_{\hat{A}(x_{0})}}}  &W(\hat{A}(x_{0})) \ar[r]_{W(proj)}   & W(A_{0}(x_{0}))\ar[r]_{W(s(x_{0}))}  &W(k(x_{0})) ,\\ 
}
\leqno{(3.1.5)}
$$
\noindent dans lequel $\sigma_{W(x_{0})}=W(\sigma_{x_{0}})$; et on a d\'ej\`a vu que 
$$
\sigma_{\mathcal{V}(x_{0})}\circ \iota(x_{0})=\iota(x_{0})\circ \sigma_{W(x_{0})}\ .
\leqno{(3.1.6)}
$$

\noindent On note $\tau_{\hat{A}(x_{0})}$ le morphisme compos\'e 
$$
\tau_{\hat{A}(x_{0})}:=\iota(x_{0})\circ W(s(x_{0}))\circ W(proj)\circ s_{F_{\hat{A}(x_{0})}}=\iota(x_{0})\circ s_{\sigma_{x_{0}}}\circ s(x_{0})\circ proj\ ;
\leqno{(3.1.7)}
$$ 
\noindent$\tau_{\hat{A}(x_{0})}$ est appel\'e le \textbf{rel\`evement de Teichm\"{u}ller de $s(x_{0})$ (ou de $x_{0}$) relativement \`a $F_{\hat{A}(x_{0})}$}.\\
\noindent Compte tenu de la commutativit\'e du diagramme (3.1.5) et de (3.1.6) on obtient l'\'egalit\'e
$$
\tau_{\hat{A}(x_{0})}\circ F_{\hat{A}(x_{0})}=\sigma_{\mathcal{V}(x_{0})}\circ \tau_{\hat{A}(x_{0})}\ .
\leqno{(3.1.8)}
$$
\noindent De plus l'\'egalit\'e (3.1.7) prouve par r\'eduction modulo $\pi$ que $\tau_{\hat{A}(x_{0})}$ est bien un rel\`evement de $s(x_{0})$.\\

Lorsque $k$ est un corps fini cette construction redonne bien celle du premier cas \'etudi\'e: en effet, gr\^ace \`a (3.1.8) on obtient l'\'egalit\'e (3.1.1), car $({\sigma_{W(x_{0})}})^{deg\ x_{0}}= id$ fournit
$$
({\sigma_{\mathcal{V}(x_{0})}})^{deg\ x_{0}}\circ \tau_{\hat{A}(x_{0})}=\iota(x_{0})\circ({\sigma_{W(x_{0})}})^{deg\ x_{0}}\circ s_{\sigma_{x_{0}}}\circ s(x_{0})\circ proj=\tau_{\hat{A}(x_{0})}\ ,
$$

\noindent d'o\`u la co\" incidence des deux constructions par l'unicit\'e dans (3.1.1).\\

 Les morphismes compos\'es

$$
\tau_{\hat{A}}(x_{0}) : \hat{A} \hookrightarrow \hat{A}(x_{0}) \xrightarrow[]{\tau_{ \hat{A}(x_{0})}} \mathcal{V}(x_{0})\ ,
\leqno{(3.1.9)}
$$

$$
\tau_{\hat{A}_{K}}(x_{0}):= \tau_{\hat{A}}(x_{0})\otimes_{\mathcal{V}}K: \hat{A}_{K}\rightarrow K(x_{0})=Frac(\mathcal{V}(x_{0}))\ ,
\leqno{(3.1.9)'}
$$

$$
\tau_{A^{\dag}}(x_{0}) : A^{\dag} \hookrightarrow \hat{A}\xrightarrow[]{\tau_{\hat{A}}(x_{0})} \mathcal{V}(x_{0})\ ,
\leqno{(3.1.10)}
$$

$$
\tau_{A^{\dag}_{K}}(x_{0}):= \tau_{A^{\dag}}(x_{0})\otimes_{\mathcal{V}}K: A^{\dag}_{K}\rightarrow K(x_{0})=Frac(\mathcal{V}(x_{0}))\ ,
\leqno{(3.1.10)'}
$$

\noindent sont appel\'es respectivement \textbf{
rel\`evement de Teichm\"{u}ller de $s$ (ou de $x_{0}$) 
relativement \`a $F_{\hat{A}}$, $F_{\hat{A}_{K}}$, $F_{A^{\dag}}$, $F_{A^{\dag}_{K}}:\tau_{\hat{A}}(x_{0})\ (\mbox{resp.} \ \tau_{\hat{A}_{K}}(x_{0}),\tau_{A^{\dag}}(x_{0}), \tau_{A^{\dag}_{K}}(x_{0})$
 est un point de Teichm\"{u}ller de $\hat{A}$ (resp. de $\hat{A}_{K}$, $A^{\dag}$, $A^{\dag}_{K}$) relativement \`a $F_{\hat{A}}$ (resp. $F_{\hat{A}_{K}}$, $F_{A^{\dag}}$, $F_{A^{\dag}_{K}}$)
	}, tel que
	
	$$
\tau_{\hat{A}}(x_{0})\circ F_{\hat{A}}=\sigma_{\mathcal{V}(x_{0})}\circ \tau_{\hat{A}}(x_{0})\ ,
\leqno{(3.1.11)}
$$

$$
\tau_{\hat{A}_{K}}(x_{0})\circ F_{\hat{A}_{K}}=\sigma_{K(x_{0})}\circ \tau_{\hat{A}_{K}}(x_{0})\ ,
\leqno{(3.1.11)'}
$$

$$
\tau_{A^{\dag}}(x_{0})\circ F_{A^{\dag}}=\sigma_{\mathcal{V}(x_{0})}\circ \tau_{A^{\dag}}(x_{0})\ ,
\leqno{(3.1.12)}
$$

$$
\tau_{A^{\dag}_{K}}(x_{0})\circ F_{A^{\dag}_{K}}=\sigma_{K(x_{0})}\circ \tau_{A^{\dag}_{K}}(x_{0})\ .
\leqno{(3.1.12)'}
$$

La correspondance entre $\tau_{\hat{A}(x_{0})}$ et $\tau_{\hat{A}}(x_{0})$ est bijective gr\^ace au $\mathbb{Z}$-isomorphisme d'adjonction
$$
Hom_{\mathcal{V}(x_{0})}(\hat{A}\otimes_{\mathcal{V}}\mathcal{V}(x_{0}), \mathcal{V}(x_{0}))\simeq Hom_{\mathcal{V}}(\hat{A}, \mathcal{V})\ ;
$$
\noindent de m\^eme le passage de $\tau_{A^{\dag}}(x_{0})$ \`a $\tau_{\hat{A}}(x_{0})$ est bijectif comme passage au compl\'et\'e.\\

Au lieu de consid\'erer $x_{0}$ comme point de $X_{0}$ \`a valeur dans $k'$ on peut aussi le consid\'erer comme point de $\mathbb{A}_{k}^{n}= Spec\ k[\underline{T}]$ \`a valeur dans $k'$, o\`u $\underline{T}= (T_{1},...,T_{n})$: on va voir qu'alors  $\tau_{A^{\dag}}(x_{0})$ s'\'etend en un point de Teichm\"uller de $\mathcal{V}[\underline{T}]^{\dag}$. Notons $R=\mathcal{V}[\underline{T}],\ R^{\dag}$ (resp. $\hat{R}$) son compl\'et\'e faible (resp. son s\'epar\'e compl\'et\'e) et $I$ le noyau de la surjection canonique
$$
\mu: R^{\dag}\twoheadrightarrow A^{\dag}\ .
$$
\noindent On peut relever (de mani\`ere non unique) le Frobenius $F_{A^{\dag}}$ de $A^{\dag}$ en un endomorphisme $F_{R^{\dag}}$ de $R^{\dag}$ de la fa\c con suivante: il suffit de choisir des \'el\'ements $F_{R^{\dag}}(T_{i})$ de $R^{\dag}$ tels que

$$
F_{R^{\dag}}(T_{i})\equiv T_{i}^{q} (\mbox{mod}\ \pi),\ F_{R^{\dag}}(T_{i})\in \mu^{-1}(F_{A^{\dag}}(\mu(T_{i})))\ ,
$$
\noindent choix qu'il est possible de faire car $F_{A^{\dag}}(\mu(T_{i}))\equiv T^{q}_{i}$ ( mod ($\pi, I$)). On \'etend ce choix d'\'el\'ements $F_{A^{\dag}}(\mu(T_{i}))$ en un endomorphisme $F_{R^{\dag}}$ de la $\mathcal{V}$-alg\`ebre $R^{\dag}$ tel que le diagramme
$$
\xymatrix{
R^{\dag} \ar@{->>}[r]^{\mu} \ar[d]_{F_{R^{\dag}}} & A^{\dag} \ar[d]^{F_{A^{\dag}}}\\
R^{\dag} \ar@{->>}[r]^{\mu}  &A^{\dag}
}
$$
\noindent commute; en particulier on a $F_{R^{\dag}}(I)\subset I$.\\
Les morphismes $F_{R^{\dag}}$ et $F_{A^{\dag}}$ sont finis et fid\`elements plats puisque leur r\'eduction mod $ \pi$ le sont [Et 2, th\'eo 17]; par changement de base de $R^{\dag}$ \`a $\hat{R}$ appliqu\'e au diagramme pr\'ec\'edent on en d\'eduit le diagramme commutatif

$$
\xymatrix{
\hat{R} \ar@{->>}[r]^{\hat{\mu}} \ar[d]_{F_{\hat{R}}} &\hat{A} \ar[d]^{F_{\hat{A}}}\\
\hat{R} \ar@{->>}[r]^{\hat{\mu}}  &\hat{A}&. } 
$$

\noindent Par cons\'equent le $k$-morphisme compos\'e
$$
k[\underline{T}]\twoheadrightarrow A_{0}  \displaystyle \mathop{\longrightarrow}^{s} k(x_{0})
$$  
\noindent correspondant au point $x_{0}$ de $\mathbb{A}^{n}_{k}$ se rel\`eve en un $\mathcal{V}$-morphisme

$$
\tau_{\hat{R}}(x_{0}) : \hat{R} \hookrightarrow \hat{A} \displaystyle \mathop{\longrightarrow}^{\tau_{\hat{A}}(x_{0})} \mathcal{V}(x_{0})
$$
\noindent tel que le diagramme
$$
\xymatrix@C=1,5cm
{
\hat{R}  \ar[r]^{\tau_{\hat{R}}(x_{0})} \ar[d]_{F_{\hat{R}}} & \mathcal{V}(x_{0}) \ar[d]^{\sigma_{\mathcal{V}(x_{0})}}\\
\hat{R}  \ar[r]^{\tau_{\hat{R}}(x_{0})}  & \mathcal{V}(x_{0})
}
$$

\noindent commute : $\tau_{\hat{R}}(x_{0})$ est appel\'e le \textbf{rel\`evement de Teichm\"{u}ller de  $x_{0}$ relativement \`a $F_{\hat{R}}$}. Le morphisme compos\'e
$$
\tau_{R^{\dag}}(x_{0}) : R^{\dag} \hookrightarrow \hat{R} \displaystyle \mathop{\longrightarrow}^{\tau_{\hat{R}}(x_{0})} \mathcal{V}(x_{0})
$$
\noindent est appel\'e le \textbf{rel\`evement de Teichm\"{u}ller de  $x_{0}$ relativement \`a $F_{R^{\dag}}$}.\\

Lorsque $k$ est un corps fini nous avons vu qu'il y a unicit\'e des rel\`evements de Teichm\"uller: dans ce cas il y a  \textbf{bijection entre les points de Teichm\"uller de $A^{\dag}$ relativement \`a $F_{A^{\dag}}$ et les points de Teichm\"uller de $R^{\dag}$ relativement \`a $F_{R^{\dag}}$ qui se r\'eduisent mod $ \pi$ en des points de $X_{0}$}. Ainsi nous avons prouv\'e:

\vskip 3mm
\noindent \textbf{Proposition (3.1.13)}. \textit{
Sous les hypoth\`eses pr\'ec\'edentes, soit $k \hookrightarrow k'$ une extension finie de corps finis. Alors il y a une bijection entre les trois ensembles suivants:}
\begin{enumerate}
\item[(i)]\textit{l'ensemble des points $x_{0}\in X_{0}(k')$= \{points de $X_{0}$ \`a valeur dans $k'$\} ,}
\item[(ii)] \textit{l'ensemble des points de Teichm\"uller de $A^{\dag}$ ( resp de $\hat{A}$) relativement \`a $F_{A^{\dag}}$ (resp $F_{\hat{A}}$) \`a valeur dans $\mathcal{V}(k'):=\ W(k')\otimes_{W} \mathcal{V}$ ,}
\item[(iii)]\textit{l'ensemble des points de Teichm\"uller de $R^{\dag}$ (resp $\hat{R}$) relativement \`a $F_{R^{\dag}}$ (resp $F_{\hat{R}}$) \`a valeur dans $\mathcal{V}(k'):=\ W(k')\otimes_{W} \mathcal{V}$ qui se r\'eduisent mod $ \pi$ en des points de $X_{0}$ \`a valeur dans $k'$.}
\end{enumerate}

\vskip 5mm
 \textbf{Supposons dor\'enavant jusqu'\`a la fin de ce \S3.1 que $x_{0}\in\vert X_{0}\vert$ est un point ferm\'e de $X_{0}$}. Alors les morphismes $\tau_{\hat{A}}(x_{0})$  et $\tau_{A^{\dag}}(x_{0})$ (resp. $\tau_{\hat{R}}(x_{0})$  et $\tau_{R^{\dag}}(x_{0})$) sont surjectifs car la r\'eduction de $\tau_{A^{\dag}}(x_{0}) \ \textrm{mod}\  \pi$ (resp.$\tau_{R^{\dag}}(x_{0}) \ \textrm{mod}\  \pi$) est le morphisme surjectif $s : A_{0} \twoheadrightarrow k(x_{0})=A_{0}/\mathfrak{m}_{x_{0}}$ de d\'epart (resp. le morphisme surjectif $k[\underline{T}]\twoheadrightarrow k(x_{0})$) [M-W, theo 3.2]. Donc $\mathcal{V}(x_{0})$ est un quotient de $\hat{A}$ et $\mathcal{V}(x_{0}) \simeq W(x_{0}) [\pi]$, qui est un anneau de valuation discr\`ete, est une extension finie \'etale de $\mathcal{V}$ de rang $deg\ x_{0}$. Le noyau du morphisme surjectif
$$
\tau_{\hat{A}_{K}}(x_{0}):= \tau_{\hat{A}}(x_{0})\otimes_{\mathcal{V}}K: \hat{A}_{K}\rightarrow K(x_{0})=Frac(\mathcal{V}(x_{0}))\ 
$$
\noindent est ainsi un id\'eal maximal $\mathfrak{q}_{x}$ de $\hat{A}_{K}$. Remarquons que le morphisme
$$
\tau_{A^{\dag}_{K}}(x_{0}):= \tau_{A^{\dag}}(x_{0})\otimes_{\mathcal{V}}K: A^{\dag}_{K} \displaystyle \mathop{\longrightarrow}^{\varphi}  \hat{A}_{K} \displaystyle -\hspace{-10pt}-\hspace{-10pt}-\hspace{-10pt}-\hspace{-10pt}-\hspace{-10pt} \mathop{\twoheadrightarrow}^{\hspace{-10pt}\tau_{\hat{A}_{K}}(x_{0})}  K(x_{0})
$$
\noindent est aussi surjectif car $\tau_{A^{\dag}}(x_{0})$ l'est; cette surjectivit\'e provient aussi du fait 
que $\varphi$ induit un isomorphisme [G-K 2, prop 1.5]
$$
A^{\dag}_{K}/Ker\ \tau_{A^{\dag}_{K}}(x_{0})\displaystyle \mathop{\longrightarrow}^{\sim} \hat{A}_{K}/Ker\ \tau_{\hat{A}_{K}}(x_{0})\ .
$$

 En d\'efinissant l'application (encore appel\'ee \textbf{rel\`evement de Teichm¬\"uller})
 $$
 Teich_{\hat{A}_{K}} : Max\ A_{0} \rightarrow Max\ \hat{A}_{K}
 \leqno{(3.1.14)}
 $$

\noindent entre l'ensemble des id\'eaux maximaux de $A_{0}$ et ceux de $\hat{A}_{K}$ par $Teich_{\hat{A}_{K}}(x_{0}) = \{ \textrm{Ker}\ \tau_{\hat{A}_{K}}(x_{0}) \} = \{ \mathfrak{q}_{x} \}$, nous montrerons au \S3.2 que $Teich_{\hat{A}_{K}}$ est \textbf{une section du  morphisme de sp\'ecialisation} [B 3,(0.2.2.1)], consid\'er\'e comme une application

$$
\overline{sp} : Max\ \hat{A}_{K} \ \rightarrow Max\ A_{0} 
$$
\noindent entre l'espace sous-jacent \`a $Spm\ \hat{A}_{K}$ et l'ensemble des points ferm\'es de $X_{0}=\ Spec\ A_{0}$\ .\\

\subsection*{3.2. Morphismes de sp\'ecialisation et de r\'eduction}

Avec les notations de 3.1 soient $\mathcal{X}=Spf\ \hat{A}$ et $\mathcal{X}_{K}= Spm\ \hat{A}_{K}$. Nous allons examiner de plus pr\`es la description donn\'ee par Berthelot [B3, (0.2.2)] de $\mathcal{X}_{K}$ et du morphisme de sp\'ecialisation
$$
\overline{sp}:\ Spm\ \hat{A}_{K}\rightarrow Spec\ A_{0}\ ,
$$
\noindent vu ici comme une application entre les ensembles sous-jacents de points ferm\'es
$$
\overline{sp}:\ \vert Spec\ \hat{A}_{K}\vert= Max\ \hat{A}_{K} \rightarrow \vert Spec\ A_{0}\vert = Max\ A_{0} ,\
$$
\noindent pour pouvoir montrer que $T_{\hat{A}_{K}}$ fournit une section de $\overline{sp}$.\\

Soient $x\in Spm\ \hat{A}_{K}$ correspondant \`a l'id\'eal maximal $\mathfrak{m}_{x}\subset \hat{A}_{K}$, $\overset{\circ}{\mathfrak{m}}_{x}=j^{-1}(\mathfrak{m}_{x})=\mathfrak{m}_{x}\cap\hat{A}$ l'image inverse de $\mathfrak{m}_{x}$ par l'injection canonique $j: \ \hat{A}\hookrightarrow\hat{A}_{K}$ et $R_{x}=\hat{A}/\overset{\circ}{\mathfrak{m}}_{x}:\ R_{x}$ s'identifie \`a l'image canonique de $\hat{A}$ dans $K_{x}=\hat{A}_{K}/\mathfrak{m}_{x}$ et le corps des fractions de $R_{x}$ est $R_{x}\otimes_{\mathcal{V}}K\simeq K_{x}$; on note $j_{x}: R_{x}\hookrightarrow K_{x}$ l'injection obtenue par changement de base de $\hat{A}$ \`a $R_{x}$. Berthelot montre que la correspondance
$$
x \mapsto R_{x}=\hat{A}/\overset{\circ}{\mathfrak{m}}_{x}
$$
\noindent est une bijection [B3, (0.2.2)] entre les points de $\mathcal{X}_{K}$ et les quotients $R_{x}$ de $\hat{A}$ qui sont int\`egres finis et plats sur $\mathcal{V}$.\\

La valuation discr\`ete de $K$ s'\'etend de mani\`ere unique en une valuation discr\`ete de $K_{x}$, dont l'anneau de valuation $\mathcal{V}_{x}$ est un $\mathcal{V}$-module libre de rang $n= [K_{x}:K]$ et $\mathcal{V}_{x}$ est la fermeture int\'egrale de $\mathcal{V}$ dans $K_{x}$ [Se, chap II, \S2, prop 3]. 
On dispose ainsi d'une injection canonique
$$
i_{x}:\ R_{x}\hookrightarrow\mathcal{V}_{x}
$$
\noindent qui est finie; par suite le morphisme
$$
i^{\#}_{x}:=Spec\ (i_{x}): \ Spec\ (\mathcal{V}_{x})\rightarrow Spec\ (R_{x})
$$
est fini et surjectif, donc $R_{x}$ n'a que deux id\'eaux premiers $(0)$ et $\mathfrak{m}_{R_{x}}$ qui sont distincts car l'id\'eal maximal $\mathfrak{m}_{R_{x}}$ est le seul id\'eal premier de $R_{x}$ au-dessus de 
$\mathfrak{m}_{\mathcal{V}}=\pi\mathcal{V}$ [R, chap I, $\S$1, prop 1]: en particulier on a $\pi R_{x}\subset \mathfrak{m}_{R_{x}}$ et $\mathcal{V}$ s'injecte dans l'anneau hens\'elien $R_{x}$ [$EGA\ IV$, (18.5.10)] de dimension 1 [$EGA\ II$, (7.1.5)], $Spec\ R_{x}=\{ (0),\mathfrak{m}_{R_{x}}\}$. Puisque $K_{x}$ est aussi le corps des fractions de $R_{x}$, $R_{x}$ qui est fini et plat sur $\mathcal{V}$ est un sous-$\mathcal{V}$-module libre de $\mathcal{V}_{x}$ de rang $n= [K_{x}:K]$ et l'anneau de valuation $\mathcal{V}_{x}$ domine l'anneau local $R_{x}$, i.e. $\mathfrak{m}_{R_{x}}=\mathfrak{m}_{\mathcal{V}_{x}}\cap R_{x}$\ .\\

Consid\'erons le diagramme commutatif suivant dans lequel $\mathfrak{m}_{\mathcal{V}_{x}}$ est l'id\'eal maximal de $\mathcal{V}_{x}$ et les fl\`eches sont les fl\`eches canoniques

$$
\xymatrix@C=1,2cm  
{
A_{0} \ar@{->>}[r]^{\widetilde{\psi}_{x}\ \ \ \  \ \ \ \ \ \ \ \ \ \ \ \ \ \ \ \ \ \ \   }\ar@{=}[d]&k(R_{x})=R_{x}/\mathfrak{m}_{R_{x}}=A_{0}/\mathfrak{m}_{\overline{sp}(x)}=\hat{A}/\mathfrak{m}_{sp(x)}\ar@{^{(}->}[r]^{ \ \ \ \  \  \ \ \ \ \ \ \ \ \ \ \  \  \widetilde{i}_{x}}&k(\mathcal{V}_{x})=\mathcal{V}_{x}/ \mathfrak{m}_{\mathcal{V}_{x}}\\
 A_{0}\ar@{->>}[r]^{\overline{\psi}_{x} \  \  \ \ \ \ \ \  \ \ \ \ \ \ \  \ \  \  \  \  }&R_{x}/\pi R_{x}=A_{0}/\widetilde{\mathfrak{m}}_{x}=\hat{A}/<\pi, \overset{\circ}{\mathfrak{m}}_{x}>\ar[r]^{ \  \  \  \ \ \ \ \ \ \  \   \  \  \overline{i}_{x}}\ar@{->>}[u]&\mathcal{V}_{x}/\pi \mathcal{V}_{x}\ar@{->>}[u]\\
 \hat{A} \ar@{->>}[r]^{{\psi}_{x}}\ar@{->>}[u]^{\rho}&R_{x}=\hat{A}/\overset{\circ}{\mathfrak{m}}_{x} \ar@{^{(}->}[r]^{i_{x}}\ar@{->>}[u]^{\rho_{x}}&\mathcal{V}_{x}\ar@{->>}[u]\\
 \hat{A}_{K}\ar@{->>}[r]_{\theta_{x}}\ar@{<-^{)}}[u]^{j}&K_{x}=\hat{A}_{K}/\mathfrak{m}_{x}\ar@{=}_{id}[r]\ar@{<-^{)}}[u]^{j_{x}}&K_{x}\ar@{<-^{)}}[u]^{j'_{x}}
}
\leqno{(3.2.0)}
$$
\noindent et o\`u $\overline{\psi}_{x}$, $\overline{i}_{x}$ sont obtenus \`a partir de ${\psi}_{x}$, $i_{x}$ par r\'eduction modulo $\pi$.\\
Les fl\`eches horizontales compos\'ees du diagramme sont not\'ees
$$
\theta_{x} : = id\circ\theta_{x} \ ,\ \overset{\circ}{\theta}_{x}: =i_{x}\circ\psi_{x}\ ,\ \overline{\overset{\circ}{\theta}_{x}}:=\overline{i}_{x}\circ\overline{\psi}_{x} \ ,\ \widetilde{\theta}_{x} :=\widetilde{i}_{x}\circ\widetilde{\psi}_{x} \ ,
$$

\noindent et les noyaux successifs
$$
Ker\ \theta_{x}= \mathfrak{m}_{x}\ ,
$$
$$
Ker\  \big(\ \overset{\circ}{\theta}_{x}\big)= Ker\ \psi_{x}= \overset{\circ}{\mathfrak{m}}_{x}=\overset{\circ}{\left(Ker\ \theta_{x} \right)} \hookrightarrow \big( \overset{\circ}{\theta}_{x}\big)^{-1}(\mathfrak{m}_{\mathcal{V}_{x}})=\psi^{-1}_{x}(\mathfrak{m}_{R_{x}})=:\mathfrak{m}_{sp(x)}\ ,
$$
$$
Ker\ \overline{ \psi}_{x}=\widetilde{\mathfrak{m}}_{x}=\widetilde{Ker\ \theta}_{x}= \overset{\circ}{\mathfrak{m}}_{x}/ \overset{\circ}{\mathfrak{m}}_{x}\cap\pi \hat{A}=\rho( \overset{\circ}{\mathfrak{m}}_{x})= \overset{\circ}{\mathfrak{m}}_{x}/\pi \overset{\circ}{\mathfrak{m}}_{x}\ \  [ \overset{\circ}{\mathfrak{m}}_{x}\  \mbox{est}\ 1^{ier} \mbox{et}\  \pi\notin \overset{\circ}{\mathfrak{m}}_{x}]
$$
$$
\widetilde{\mathfrak{m}}_{x}=\widetilde{Ker\ \theta}_{x}\hookrightarrow Ker\ \widetilde{\theta}_{x}= Ker\ \widetilde{\psi}_{x}=\mathfrak{m}_{\overline{sp}(x)}=Rad\ (\widetilde{\mathfrak{m}}_{x})\ .
$$

Par le morphisme de sp\'ecialisation
$$
sp:\ Max\ {\hat{A}_{K}}\rightarrow Max\ {\hat{A}}
$$
\noindent l'image de $\mathfrak{m}_{x}$ est l'id\'eal maximal
$$
\mathfrak{m}_{sp(x)}:=\psi^{-1}_{x}(\mathfrak{m}_{R_{x}})=\big( \overset{\circ}{\theta}_{x}\big)^{-1}(\mathfrak{m}_{\mathcal{V}_{x}})\ .
$$
\noindent Comme la surjection canonique
$$
\rho:\ \hat{A}\twoheadrightarrow A_{0}=\hat{A}/\pi\ \hat{A}
$$
\noindent induit une bijection
$$
\rho^{\#}:\ Max\ A_{0}\overset{\sim}{\rightarrow}Max\ \hat{A}\ ,
$$
\noindent on notera $\overline{sp}$ le compos\'e
$$
\overline{sp}: \ Max\ {\hat{A}_{K}}\overset{sp}{\longrightarrow} Max\ {\hat{A}}\overset{(\rho^{\#})^{-1}}{\longrightarrow}Max\ A_{0}\ ,
$$
\noindent et l'image de $\mathfrak{m}_{x}$ par $\overline{sp}$ est l'id\'eal maximal
$$
\mathfrak{m}_{\overline{sp}(x)}= \rho(\mathfrak{m}_{sp(x)})\ .
$$
\subsubsection*{3.2.1. Interpr\'etation g\'eom\'etrique du morphisme de sp\'ecialisation}
Avec les notations de 3.2 consid\'erons le diagramme commutatif suivant
$$
\xymatrix{
Spec\ R_{x}/ \pi R_{x}\ \ar@{^{(}->}[r]^{\ \ \ \ \ \ \rho_{x}^{\#}}& Spec\ R_{x}\ar@{<-^{)}}[r]^{j_{x}^{\#}}&Spec\ K_{x}\\
Spec\ A_{0}\ \ar@{^{(}->}[r]^{\ \ \ \rho^{\#}}\ar[d]^{pr_{0}}\ar@{<-^{)}}[u]^{\overline{\psi}_{x}^{\#}}& Spec\ \hat{A}\ar@{<-^{)}}[r]^{j^{\#}}\ar[d]^{pr}\ar@{<-^{)}}[u]^{\psi_{x}^{\#}}&Spec\ \hat{A}_{K}\ar[d]^{pr_{K}}\ar@{<-^{)}}[u]^{\theta_{x}^{\#}}\\
Spec\ k\ \ar@{^{(}->}[r]&Spec\ \mathcal{V}\ar@{<-^{)}}[r]& Spec\ K
}
\leqno{(3.2.1.1)}
$$
\noindent dans lequel $\overline{\psi}_{x}^{\#}, \psi_{x}^{\#},\theta_{x}^{\#}$ sont des immersions ferm\'ees et les fl\`eches $pr$ les projections canoniques. \\
L'image sch\'ematique de $x=\{\mathfrak{m}_{x}\}$ par $j^{\#}$ n'est autre que l'image sch\'ematique de $Spec\ K_{x}$ par le morphisme dominant $j_{x}^{\#}$, \`a savoir
$$
Spec\ (\hat{A}/\overset{\circ}{\mathfrak{m}}_{x})= Spec\ R_{x}=\{ (0),\mathfrak{m}_{R_{x}}\}\overset{\ \ \psi^{\#}_{x}}{\hooklongrightarrow}Spec\ \hat{A}\ ,
$$
\noindent et l'intersection de cette adh\'erence sch\'ematique de $\{\mathfrak{m}_{x}\}$ dans $Spec\ \hat{A}$ avec la fibre sp\'eciale $Spec\ A_{0}$ consiste en un unique point ferm\'e $\{\mathfrak{m}_{R_{x}}\}$ (o\`u $\mathfrak{m}_{R_{x}}$ est identifi\'e \`a $\mathfrak{m}_{sp(x)}:=\psi_{x}^{-1}(\mathfrak{m}_{R_{x}})$ ou \`a $\mathfrak{m}_{\overline{sp}(x)}= \rho(\mathfrak{m}_{sp(x)})$ ) qui est la sp\'ecialisation de $\{\mathfrak{m}_{x}\}$. Cette intersection s'\'ecrivant aussi $Spec (R_{x}/\pi R_{x})$, qui est r\'eduite \`a un seul point, l'anneau $\overline{R}_{x}:=R_{x}/\pi R_{x}$ est donc un anneau artinien local d'id\'eal maximal $\mathfrak{m}_{\overline{R}_{x}}$, image de $\mathfrak{m}_{R_{x}}$ par la surjection canonique 
$$R_{x}\twoheadrightarrow\overline{R}_{x}=R_{x}/\pi R_{x}\ .$$
Au passage, comme $\mathfrak{m}_{\overline{R}_{x}}$ est nilpotent [Bour, A, chap 8, $\S$6, $\no$4, cor du th\'eo 3], remarquons que la topologie $\pi$-adique et la topologie $\mathfrak{m}_{R_{x}}$-adique co\"incident sur l'anneau $R_{x}$, i.e. il existe un entier $r\in \mathbb{N}$ tel que 
$$
\mathfrak{m}^{r}_{R_{x}}\subset \pi R_{x}\subset\mathfrak{m}_{R_{x}} \ .
$$
Comme l'adh\'erence de Zariski de $\{\overset{\circ}{\mathfrak{m}}_{x}\}$ consiste en les id\'eaux premiers $\mathfrak{A}$ de $\hat{A}$ tels que $\mathfrak{A}\supset\overset{\circ}{\mathfrak{m}}_{x}$ on retrouve le fait que $\{\mathfrak{m}_{sp(x)}\}$ est \guillemotleft une\guillemotright\ sp\'ecialisation, au sens classique du terme [$EGA\ 0_{I},\ (2.1.1)$], de $\{\overset{\circ}{\mathfrak{m}}_{x}\}$.

\subsubsection*{3.2.2. Lien entre Teichm\" uller et sp\'ecialisation}
Toujours avec les notations de 3.2, soit $x_{0}$ un point ferm\'e de $Spec\ A_{0}$, correspondant \`a un id\'eal maximal $\mathfrak{m}_{x_{0}}$ de $A_{0}$ et \`a un id\'eal maximal $\mathfrak{q}_{x_{0}}$ de $\hat{A}$. On note $x=Teich_{\hat{A}_{K}}(x_{0})$ l'image de $x_{0}$ par le rel\`evement de Teichm\" uller $Teich_{\hat{A}_{K}} : Max\ A_{0} \rightarrow Max\ \hat{A}_{K}$; soient $\mathfrak{m}_{x}\subset\hat{A}_{K}$ l'id\'eal correspondant \`a $x$ et $\overset{\circ}{\mathfrak{m}}_{x}=j^{-1}(\mathfrak{m}_{x})=\mathfrak{m}_{x}\cap\hat{A}$.\\

Cette fois-ci $K_{x}=K(x_{0})=\hat{A}_{K}/\mathfrak{m}_{x}$ est une extension finie non ramifi\'ee de $K$ et on a $\theta_{x}=\tau_{\hat{A}_{K}}(x_{0}), \psi_{x}=\overset{\circ}{\theta}_{x}=\tau_{\hat{A}}(x_{0}): \hat{A}\twoheadrightarrow\mathcal{V}_{x}=\mathcal{V}(x_{0}), i_{x}$ est un isomorphisme et $R_{x}=\mathcal{V}_{x}$, qui est non ramifi\'e sur $\mathcal{V}$, admet $\pi$ pour uniformisante: l'id\'eal $<\pi,\overset{\circ}{\mathfrak{m}}_{x}>$ de $\hat{A}$ engendr\'e par $\pi$ et $\overset{\circ}{\mathfrak{m}}_{x}$ est \'egal \`a $\mathfrak{m}_{sp(x)}$ et l'anneau $\overline{R}_{x}=R_{x}/\pi R_{x}$ est non seulement artinien local, c'est m\^eme le corps $k(sp(x))=\hat{A}/\mathfrak{m}_{sp(x)}=A_{0}/\mathfrak{m}_{\overline{sp}(x)}$ qui est aussi le corps r\'esiduel $k(\mathcal{V}_{x})$ de $\mathcal{V}_{x}=\mathcal{V}(x_{0})$. D'o\`u $\mathfrak{m}_{sp(x)}=\mathfrak{q}_{x_{0}}$ et $\overline{sp}(x)=x_{0}$; ainsi
$$
Teich_{\hat{A}_{K}} : \vert Spec\ A_{0}\vert :=Max\ A_{0} \rightarrow\vert Spec\ \hat{A}_{K}\vert:= Max\ \hat{A}_{K}
$$
\noindent est une section du morphisme de sp\'ecialisation vu comme une application entre ensembles de points ferm\'es
$$
\overline{sp}:\ \vert Spec\ \hat{A}_{K}\vert= Max\ \hat{A}_{K} \rightarrow \vert Spec\ A_{0}\vert = Max\ A_{0} \ .
$$
\noindent Dans notre situation d'un sch\'ema affine et lisse sur $\mathcal{V}$ \textbf{ceci red\'emontre la surjectivit\'e du morphisme de sp\'ecialisation} \'etablie par Berthelot [B3, (1.1.5)].

\subsubsection*{3.2.3. Morphisme \guillemotleft$\alpha$\guillemotright\  de Tate et sp\'ecialisation}
L'anneau local $R_{x}$ \'etant domin\'e par l'anneau de valuation $\mathcal{V}_{x}$ de $K_{x}$ on a 
$i_{x}^{-1}(\mathfrak{m}_{\mathcal{V}_{x}})=\mathfrak{m}_{R_{x}}$. Ainsi \textbf{le morphisme \guillemotleft$\alpha$\guillemotright\  de Tate} [Ta, theo 6.4]  \textbf{est le morphisme}
$$
\begin{array}{ccccc}
Max\ {\hat{A}_{K}}&\longrightarrow& Max\ \hat{A}\\
\mathfrak{m}_{x}& \longmapsto &\big( \overset{\circ}{\theta}_{x} \big)^{-1}(\mathfrak{m}_{\mathcal{V}_{x}})=&\psi_{x}^{-1}(\mathfrak{m}_{R_{x}})=&\mathfrak{m}_{sp(x)}\ ;
\end{array}
$$
c'est donc le morphisme \textbf{de sp\'ecialisation}
$$
sp:\ Max\ {\hat{A}_{K}}\rightarrow Max\ {\hat{A}}\ , x\mapsto \{\mathfrak{m}_{sp(x)} \}
$$
de Berthelot [B3, (0.2.2)], dont on retrouve la surjectivit\'e dans [Ta, theo 6.4].

\subsubsection*{3.2.4. Morphisme \guillemotleft$\pi$\guillemotright\ de Bosch-G\" untzer-Remmert, morphisme de r\'eduction et sp\'ecialisation}
Pour un annneau $B$, on note $Nilrad(B)$ (resp $Rad(B)$) le nilradical de $B$ (resp le radical de Jacobson de $B$), intersection des id\'eaux premiers (resp maximaux) de $B$. Dans le cas de notre $\mathcal{V}$-alg\`ebre lisse $A$ on note $\overset{\circ}{\left(\hat{A}_{K}\right)}$ l'ensemble des \'el\'ements \`a puissances born\'ees de $\hat{A}_{K}$ et $\overset{\lor}{\left(\hat{A}_{K}\right)}$ l'ensemble des \'el\'ements topologiquement nilpotents de $\hat{A}_{K}$. D'apr\`es Tate [Ta, Theo 5.1, 5.2] on a 
$$
\overset{\circ}{\left(\hat{A}_{K}\right)}=\{ f\in\hat{A}_{K}/\  \forall \ x\in Spm\ \hat{A}_{K},\ \vert f(x)\vert\leqslant 1   \}\ ,
$$
$$
\overset{\lor}{\left(\hat{A}_{K}\right)}=\{ f\in\hat{A}_{K}/\  \forall \ x\in Spm\ \hat{A}_{K},\ \vert f(x)\vert< 1   \}\ .
$$
\vskip 3mm
\noindent \textbf{Proposition (3.2.4.1)}. \textit{
Sous les hypoth\`eses pr\'ec\'edentes on a:}
\begin{enumerate}
\item[(i)] $\overset{\circ}{\left(\hat{A}_{K}\right)}=\hat{A}\ .$
\item[(ii)]$\overset{\lor}{\left(\hat{A}_{K}\right)}= Rad\ \hat{A}=\pi\ \hat{A}\ .$
\end{enumerate}
\vskip 3mm
\noindent\textit{D\'emonstration}.\\
Comme $A$ est normal noeth\'erien on peut d\'ecomposer $S=Spec\ A$ en la somme finie de ses composantes connexes $\displaystyle \mathop{\coprod_{i\in\llbracket{1,r}\rrbracket}} Spec\ A_{i}$ o\`u $A_{i}$ est int\`egre [EGA I, (4.5.5)] et m\^eme int\'egralement clos [Et 2, prop11(1)(ii)]; par suite $\hat{A}=\prod\limits_{i\in\llbracket{1,r}\rrbracket}\hat{A}_{i}$ o\`u $\hat{A}_{i}$ est int\'egralement clos [loc. cit.] et $\hat{A}_{K}=\prod\limits_{i\in\llbracket{1,r}\rrbracket}\hat{A}_{iK}$ o\`u $\hat{A}_{iK}=\hat{A}_{i}\otimes_{\mathcal{V}}K$ est int\'egralement clos [Bour, AC V, \S1, \no5, cor 1 de prop 16]. Puisqu'un id\'eal maximal $\mathfrak{m}_{x}\subset \hat{A}_{K}$ est de la forme
$$
\mathfrak{m}_{x}=\hat{A}_{1,K}\times\ ...\times\hat{A}_{i-1,K}\times\mathfrak{m}_{x_{i}}\times\hat{A}_{i+1,K}\times\ ...\times\hat{A}_{r,K}
$$
\noindent o\`u $\mathfrak{m}_{x_{i}}\subset \hat{A}_{i,K}$ est maximal, on a, pour $f=(f_{1},...,f_{r})\in\hat{A}_{K}=\prod\limits_{i\in\llbracket{1,r}\rrbracket}\hat{A}_{i,K}$\ ,\\

$f(x):=f \ \mbox{mod}\ \mathfrak{m}_{x}\  \mbox{dans}\  K_{x}:=\hat{A}_{K}/\mathfrak{m}_{x}\simeq\hat{A}_{i,K}/\mathfrak{m}_{x_{i}}=K_{x_{i}}$\\

\qquad $ \ \ =f_{i}\  \mbox{mod}\  \mathfrak{m}_{x_{i}}=:f_{i}(x_{i})\  \mbox{dans}\ K_{x_{i}}\ .$\\

\noindent On en d\'eduit donc que
$$
\overset{\circ}{\left(\hat{A}_{K}\right)}=\prod\limits_{i\in\llbracket{1,r}\rrbracket}\overset{\circ}{\left(\hat{A}_{i,K}\right)}
$$
\noindent et 
$$
\overset{\lor}{\left(\hat{A}_{K}\right)}=\prod\limits_{i\in\llbracket{1,r}\rrbracket}\overset{\lor}{\left(\hat{A}_{i,K}\right)}\ .
$$
\noindent Pour d\'emontrer la proposition on est donc ramen\'e au cas $\hat{A}$ int\'egralement clos. En se donnant une pr\'esentation de la $\mathcal{V}$-alg\`ebre $A$
$$
\hat{A}\simeq \mathcal{V}\{ T_{1}, ...,T_{n}\}/ \mathfrak{A}\ ,
$$
on note $t_{i}$ l'image de $T_{i}$ dans $\hat{A}$; la $K$-alg\`ebre de Tate $\hat{A}_{K}$ s'\'ecrit
$$
\hat{A}_{K}\simeq K\{ T_{1}, ...,T_{n}\}/ \mathfrak{A}_{K}\ ,
$$
 de sorte que
$$
\hat{A}_{K}\simeq K\{ t_{1}, ...,t_{n}\}\ .
$$
D'apr\`es Tate [Ta, remarques entre les lemmes 6.2 et 6.3 affirmant la validit\'e du corollaire 2, \S5, dans le cas $\mathcal{V}$ noeth\'erien] le sous-ensemble $\overset{\circ}{\left(\hat{A}_{K}\right)}$ de $\hat{A}_{K}$ est la fermeture int\'egrale de $\mathcal{V}\{t_{1},...,t_{n}\}$ dans $\hat{A}_{K}$: or les injections \'evidentes
$$
\mathcal{V}\{t_{1},...,t_{n}\}\hookrightarrow \hat{A}\hookrightarrow\overset{\circ}{\left(\hat{A}_{K}\right)}
$$
et le fait que $\hat{A}$ est int\'egralement clos prouvent que 
$$
\overset{\circ}{\left(\hat{A}_{K}\right)} =\hat{A}\ .
$$

Puisque $\overset{\lor}{\left(\hat{A}_{K}\right)}$ est un sous-ensemble de $\overset{\circ}{\left(\hat{A}_{K}\right)}$ on peut \'ecrire 
\begin{eqnarray*}
\overset{\lor}{\left(\hat{A}_{K}\right)}&=&\{ f\in\hat{A}_{K}/\  \forall \ x\in Spm\ \hat{A}_{K},\ \vert f(x)\vert< 1   \} \\
&= &\{ f\in\hat{A}_{K}/\  \forall \ x\in Spm\ \hat{A}_{K},\ f\in\mathfrak{m}_{sp(x)}\} \\
&=&Rad\ \hat{A}\ ,
\end{eqnarray*}
la derni\`ere \'egalit\'e provenant de la surjectivit\'e du morphisme de sp\'ecialisation
$$
sp:\ Max\ {\hat{A}_{K}}\rightarrow Max\ {\hat{A}}\ , x\mapsto \{\mathfrak{m}_{sp(x)} \}\ .
$$
Or $A_{0}$ \'etant un anneau de Jacobson [Bour, AC V, \S3, \no4, d\'ef 1 et th\'eo 3] on a 
$$
Nilrad(A_{0})= Rad\ A_{0}\ ,
$$
d'o\`u $Rad\ A_{0}= \{0\}$ puisque $A_{0}$ est r\'eduit, ce qui \'equivaut \`a $Rad\ \hat{A}=\pi\ \hat{A}$; par suite
$$
\overset{\lor}{\left(\hat{A}_{K}\right)}=Rad\ \hat{A}=\pi\ \hat{A} \ . \qquad\qquad\square
$$\\

Si $B$ est un anneau, $\mathcal{J}(B)$ d\'esigne l'ensemble des id\'eaux de $B$ et, pour $I\in \mathcal{J}(B)$, $Rad(I):=Nilrad(B/I)$ est le radical de $I$
$$
Rad(I)=\{b\in B/ \exists\ \alpha\in{\mathbb{N}}, b^{\alpha}\in I \}=\bigcap_{\substack{\mathfrak{P}\supset I\\
\mathfrak{P}\in Spec B  }}\mathfrak{P}\ .
$$
Pour $I\in \mathcal{J}(\hat{A}_{K})$ on note avec [B-G-R, 1.2.5]
$$
\begin{array}{ccc}
\overset{\circ}{I}=&I\cap \overset{\circ}{\big(\hat{A}_{K}\big)}=& I\cap \hat{A}\ ,\\
\overset{\lor}{I}=&I \cap \overset{\lor}{\big(\hat{A}_{K}\big)}=&I\cap(\pi\hat{A})\ ,\\
\widetilde{I}=&\overset{\circ}{I}/\overset{\lor}{I}&\ .
\end{array}
$$
Cette notation \'etend celle du d\'ebut du \S3.2 utilis\'ee pour $\mathfrak{m}_{x}\in Max\ \hat{A}_{K}$:
$$
\begin{array}{ccccc}
\overset{\lor}{\mathfrak{m}}_{x}=&\mathfrak{m}_{x}\cap (\pi\hat{A})&\hookrightarrow &\overset{\circ}{\mathfrak{m}}_{x}=&\mathfrak{m}_{x}\cap\hat{A}\\
\widetilde{\mathfrak{m}}_{x}=&\overset{\circ}{\mathfrak{m}}_{x}/\overset{\lor}{\mathfrak{m}}_{x}\ .
\end{array}
$$
Par analogie avec [B-G-R, 2.1.10, 7.1.5] nous d\'efinissons l'application de r\'eduction $r\acute{e}d$ par
$$
\begin{array}{ccccccc}
r\acute{e}d:& Spm\ \hat{A}_{K}&\longrightarrow&\mathcal{J}(A_{0})&&&\\
         &x=\{\mathfrak{m}_{x}\} &\longmapsto &r\acute{e}d(x):=&\widetilde{\mathfrak{m}}_{x}=&\rho(\overset{\circ}{\mathfrak{m}}_{x})=&\rho(<\pi, \overset{\circ}{\mathfrak{m}}_{x}>)\ .
\end{array}
$$\\

D'apr\`es nos descriptions pr\'ec\'edentes des \'el\'ements \`a puissances born\'ees (resp des \'el\'ements topologiquement nilpotents) de $\hat{A}_{K}$, les r\'eductions $\widetilde{\hat{A}_{K}}$ et $\widetilde{K_{x}}$ au sens de [B-G-R, 1.2.5, 2.1.10] de $\hat{A}_{K}$ et $K_{x}=\hat{A}_{K}/\mathfrak{m}_{x}$ respectivement (pour $x\in Spm\ \hat{A}_{K}$) s'\'ecrivent
$$
\begin{array}{cccc}
\widetilde{\hat{A}_{K}}:=& \overset{\circ}{\left(\hat{A}_{K}\right)}/\overset{\lor}{\left(\hat{A}_{K}\right)}&=\hat{A}/\pi\ \hat{A}&=A_{0} \ , \\
\widetilde{K_{x}}:=&\mathcal{V}_{x}/\mathfrak{m}_{\mathcal{V}_{x}}&=:k(\mathcal{V}_{x})\  ,&
\end{array}
$$
et  la r\'eduction $\widetilde{\theta}_{x}$ au sens de [B-G-R, \S 6.3] du morphisme canonique
$$
\theta_{x}: \hat{A}_{K}\twoheadrightarrow K_{x}=\hat{A}_{K}/\mathfrak{m}_{x}
$$
s'identifie au morphisme compos\'e
$$
\widetilde{\theta}_{x}:A_{0}=\hat{A}/\pi \hat{A}\overset{\widetilde{\psi}_{x}}{\twoheadrightarrow}k(R_{x})\overset{\widetilde{i}_{x}}{\hookrightarrow}k(\mathcal{V}_{x})=\mathcal{V}_{x}/\mathfrak{m}_{\mathcal{V}_{x}}\ .
$$
Le morphisme \guillemotleft$\pi$\guillemotright\  de [B-G-R, \S 7.1.5]  est le morphisme
$$
\begin{array}{ccccc}
\mbox{\guillemotleft}\pi\mbox{\guillemotright} :&Spm\ {\hat{A}_{K}}&\longrightarrow& Max\ A_{0}\\
&x=\{\mathfrak{m}_{x}\}& \longmapsto &Ker\ \widetilde{\theta}_{x}=&Ker\ \widetilde{\psi}_{x}\ ;
\end{array}
$$
or
$$
Ker\ \widetilde{\theta}_{x}= Ker\ \widetilde{\psi}_{x}=\rho(\mathfrak{m}_{sp(x)})=:\mathfrak{m}_{\overline{sp}(x)}\ ,
$$
 donc \textbf{ce morphisme \guillemotleft$\pi$\guillemotright\  de Bosch-G\" untzer-Remmert [B-G-R]} (qu'il ne faut pas confondre avec notre uniformisante $\pi$) \textbf{s'identifie au morphisme de sp\'ecialisation}
$$
\overline{sp}: Max\ {\hat{A}_{K}}\rightarrow Max\ A_{0}
$$
dont on retrouve la surjectivit\'e dans [B-G-R, \S 7.1.5, theo 4].\\

Nous allons comparer dans la proposition (3.2.4.2) ci-dessous les morphismes $r\acute{e}d$ et  \guillemotleft$\pi$\guillemotright\  en donnant une CNS pour que $\widetilde{\mathfrak{m}}_{x}$ soit maximal et de ce fait co\" incide avec $\mathfrak{m}_{\overline{sp}(x)}$, i.e. on donnera une CNS pour que l'injection canonique
$$
\widetilde{\mathfrak{m}}_{x}=\widetilde{Ker\ \theta}_{x}\hookrightarrow Ker\ \widetilde{\theta}_{x}= Ker\ \widetilde{\psi}_{x}=\mathfrak{m}_{\overline{sp}(x)}
$$
soit un isomorphisme; la proposition 1 de [B-G-R, \S7.1.5] qui donne la condition suffisante $n=[K_{x}:K]=1$ en est un cas particulier.\\
Dans le \S3.2.2 sur le lien entre Teichm\" uller et la sp\'ecialisation on a constat\'e qu'avec $x=Teich_{\hat{A}_{K}}(x_{0})$ on a 
$$
\begin{array}{ccc}
(i)&\ <\pi,\overset{\circ}{\mathfrak{m}}_{x}>=\mathfrak{m}_{sp(x)}&\\
\\
(ii)&		\left\{
		\begin{array}{c}
		K_{x}\ \mbox{est une extension alg\'ebrique finie non ramifi\'ee de}\ K\ , \\
		\mbox{le corps r\'esiduel}\ \widetilde{K_{x}}\ \mbox{de}\ K_{x}\ \mbox{est un quotient de}\ A_{0}\ .
		\end{array}
		\right.
\end{array}
$$
En fait ces deux propri\'et\'es $(i)$ et $(ii)$ sont \'equivalentes comme nous allons le voir aussi dans la proposition (3.2.4.2) ci-dessous.

\vskip 3mm
\noindent \textbf{Proposition (3.2.4.2)}. \textit{Soient $A$ une $\mathcal{V}$-alg\`ebre lisse, $x\in\mathcal{X}_{K}=Spm\ \hat{A}_{K}$, \\
correspondant \`a l'id\'eal maximal $\mathfrak{m}_{x}\subset\hat{A}_{K},\ K_{x}=\hat{A}_{K}/\mathfrak{m}_{x}, \ \overset{\circ}{\mathfrak{m}}_{x}= \mathfrak{m}_{x}\cap\hat{A},\\
 R_{x}=\hat{A}/\overset{\circ}{\mathfrak{m}}_{x},\ \mathfrak{m}_{R_{x}}$ l'id\'eal maximal de $R_{x}$, $\mathcal{V}_{x}$ l'anneau de valuation de $K_{x}$, $\widetilde{K_{x}}:=\mathcal{V}_{x}/\mathfrak{m}_{\mathcal{V}_{x}}=:k(\mathcal{V}_{x})$ le corps r\'esiduel de $\mathcal{V}_{x},\ \mathfrak{m}_{\overline{sp}(x)}:=\rho(\mathfrak{m}_{sp(x)}), \\
  \widetilde{\mathfrak{m}}_{x}=\rho(\overset{\circ}{\mathfrak{m}}_{x})=\rho(<\pi, \overset{\circ}{\mathfrak{m}}_{x}>)$.\ Alors on a:}
\begin{enumerate}
\item[1)]$\mathfrak{m}_{\overline{sp}(x)}=Rad\ (\widetilde{\mathfrak{m}}_{x}).$
\item[2)]\textit{Les propri\'et\'es suivantes sont \'equivalentes:}
	\begin{enumerate}
	\item[$(i)$]$\widetilde{\mathfrak{m}}_{x}=\mathfrak{m}_{\overline{sp}(x)},\  i.e. \ r \acute{e}d(x)=\overline	{sp}(x),\ i.e.\ \widetilde{Ker\ \theta}_{x}\simeq Ker\ \widetilde{\theta}_{x} \ .$
	\item[$(i)'$]$\widetilde{\mathfrak{m}}_{x}=Rad\ (\widetilde{\mathfrak{m}}_{x}).$
	\item[$(ii)$]$<\pi,\overset{\circ}{\mathfrak{m}}_{x}>=\mathfrak{m}_{sp(x)},\ i.e.\ <\pi,\overset{\circ}{\mathfrak{m}}_{x}>\ est\ \acute{e}gal\ \grave{a}\ son \ radical.$
	\item[$(iii)$]$\widetilde{\mathfrak{m}}_{x}\ est \ un\ id\acute{e}al\ maximal.$
	\item[$(iii)'$]$<\pi,\overset{\circ}{\mathfrak{m}}_{x}>\ est \ un\ id\acute{e}al\ maximal.$
	\item[$(iv)$]$\widetilde{\mathfrak{m}}_{x}\ est \ un\ id\acute{e}al\ premier.$
	\item[$(iv)'$]$<\pi,\overset{\circ}{\mathfrak{m}}_{x}>\ est \ un\ id\acute{e}al\ premier.$
	\item[$(v)$]$R_{x}/\pi R_{x}=A_{0}/\widetilde{\mathfrak{m}}_{x}\ est \ un\ corps,\ i.e.\ \pi\ engendre\ l'id\acute{e}al\ maximal\ de\ R_{x} \ . $
	\item[$(v)'$]$R_{x}/\pi R_{x}=A_{0}/\widetilde{\mathfrak{m}}_{x}\ est \ un\ anneau\ int\grave{e}gre.$
	\item[$(v)''$]$R_{x}/\pi R_{x}=A_{0}/\widetilde{\mathfrak{m}}_{x}\ est \ un\ anneau\ r\acute{e}duit.$
	\item[$(vi)$]$L'adh\acute{e}rence\ sch\acute{e}matique\ de\ x\ dans\ Spec \hat{A}\ intersecte\ la \ fibre\ sp\acute{e}ciale\\
	 Spec\ A_{0}\ selon\ un\ sch\acute{e}ma\ int\grave{e}gre.$
	\item[$(vi)'$]$L'adh\acute{e}rence\ sch\acute{e}matique\ de\ x\ dans\ Spec \hat{A}\ intersecte\ la \ fibre\ sp\acute{e}ciale\\
	 Spec\ A_{0}\ selon\ un\ sch\acute{e}ma\ r\acute{e}duit.$
	\item[$(vii)$]$K_{x}\ est\ une\ extension\ alg\acute{e}brique\ finie\ non\ ramifi\acute{e}e\ de\ \ K\  et\ le\ corps\\
	 r\acute{e}siduel\ \widetilde{K_{x}}\ de\ K_{x}\ est\ un\ quotient\ de\ A_{0} \ .$
	 \item[$(vii)'$]$K_{x}\ est\ une\ extension\ alg\acute{e}brique\ finie\ non\ ramifi\acute{e}e\ de\ \ K\  et\ la\\
	  r\acute{e}duction\ \widetilde{K_{x}}\ de\ K_{x}\ est\ le\ corps\	 r\acute{e}siduel\ d'un\ point\ ferm\acute{e}\ de\ la\\
	   r\acute{e}duction  \ X_{0}=Spec A_{0}\ de\ \mathcal{X}_{K}=Spm \hat{A}_{K}.$
	\item[$(viii)$]$\mathcal{V}_{x}\ est \ une \ extension\ non\ ramifi\acute{e}e\ de\ \mathcal{V}\ et\ c'est\ un \ quotient\ de \  \hat{A}.$
	\item[$(viii)'$]$\mathcal{V}_{x}\ est \ une \ extension\ non\ ramifi\acute{e}e\ de\ \mathcal{V} \ et\ k(R_{x})=k(\mathcal{V}_{x}).$
	\item[$(ix)$]$R_{x}\ est \ un \ anneau \ de\ valuation\ discr\grave{e}te\ avec\ \pi\ pour\ uniformisante.$
	\item[$(ix)'$]$R_{x}\ est \ une \ extension\ non\ ramifi\acute{e}e\ de\  \mathcal{V}.$
	\item[$(ix)''$]$R_{x}\ est \ une \ extension\ finie\ \acute{e}tale \ de\  \mathcal{V}.$
	\end{enumerate}
\end{enumerate}

\vskip 3mm
\noindent\textit{D\'emonstration}.\\

1) Puisque $\pi R_{x}\subset\mathfrak{m}_{R_{x}}$ on a $<\pi,\overset{\circ}{\mathfrak{m}}_{x}>\ \subset\mathfrak{m}_{sp(x)}$, d'o\`u 
$$ 
\widetilde{\mathfrak{m}}_{x}=\rho(\overset{\circ}{\mathfrak{m}}_{x})=\rho(<\pi, \overset{\circ}{\mathfrak{m}}_{x}>)\subset \mathfrak{m}_{\overline{sp}(x)}:=\rho(\mathfrak{m}_{sp(x)})\ ;
$$
or il r\'esulte de l'\'egalit\'e $Spec\ R_{x}=\{(0),\mathfrak{m}_{R_{x}}\}$ que $\mathfrak{m}_{sp(x)}$ est le seul id\'eal premier de $\hat{A}$ contenant $<\pi,\overset{\circ}{\mathfrak{m}}_{x}>$, d'o\`u l'\'egalit\'e
$$
\mathfrak{m}_{sp(x)}=Rad\ (<\pi,\overset{\circ}{\mathfrak{m}}_{x}>)
$$
gr\^ace \`a [B-G-R, 7.1.3, prop 1], et par suite
$$
\mathfrak{m}_{\overline{sp}(x)}=Rad\ (\widetilde{\mathfrak{m}}_{x})\ .
$$

2) L'\'equivalence de $(i)$ et $(i)'$ r\'esulte du 1) et celle de $(i)$ et $(ii)$ est tautologique.\\

Dire que $  \widetilde{\mathfrak{m}}_{x}=\rho(\overset{\circ}{\mathfrak{m}}_{x})=\rho(<\pi, \overset{\circ}{\mathfrak{m}}_{x}>)$ est premier (resp maximal) \'equivaut \`a dire que $<\pi, \overset{\circ}{\mathfrak{m}}_{x}>$ est premier (resp maximal), i.e. l'anneau artinien local $R_{x}/\pi R_{x}=A_{0}/\widetilde{\mathfrak{m}}_{x}=\hat{A}/<\pi, \overset{\circ}{\mathfrak{m}}_{x}>$ est int\`egre (resp un corps); d'o\`u l'\'equivalence de $(iii),\ (iii)', \ (iv), \ (iv)',\ (v),\ (v)'.$\\

L'id\'eal maximal $\mathfrak{m}_{\overline{R}_{x}}$ de l'anneau artinien local $\overline{R}_{x}=R_{x}/\pi R_{x}$ est l'ensemble des \'el\'ements nilpotents de $\overline{R}_{x}$ [Bour, A, chap 8, \S6, \no4, cor du th\'eo 3]; d'o\`u l'\'equivalence de $ (v), (v)''.$\\

L'assertion $(v)$ qui \'equivaut \`a dire que  $\pi$ engendre l'id\'eal maximal $\mathfrak{m}_{R_{x}}$ de l'anneau local $R_{x}$, ou encore, par d\'efinition de $\mathfrak{m}_{sp(x)}:=\psi_{x}^{-1}(\mathfrak{m}_{R_{x}})$, que $<\pi,\overset{\circ}{\mathfrak{m}}_{x}>=\mathfrak{m}_{sp(x)}$; d'o\`u l'\'equivalence de $(ii)$ et $(v)$.\\

Gr\^ace \`a l'interpr\'etation g\'eom\'etrique de la sp\'ecialisation vue ci-dessus l'intersection de l'adh\'erence sch\'ematique de $x$ dans $Spec\hat{A}$ avec la fibre sp\'eciale $Spec A_{0}$ est le sch\'ema artinien local $SpecR_{x}/\pi R_{x}$; d'o\`u l'\'equivalence de $(v)''$ et $(vi)'$ et celle de $(v)'$ et $(vi)$.\\ 

Montrons $(ii)\Longrightarrow(vii)$. D'apr\`es $(ii)$ l'id\'eal maximal de $R_{x}$ est engendr\'e par $\pi$, donc l'anneau local noeth\'erien int\`egre $R_{x}$ est un anneau de valuation discr\`ete [Se, chap 1, \S2, prop 2]. Comme $K_{x}$ est le corps des fractions de $R_{x}$, il en r\'esulte que $R_{x}= \mathcal{V}_{x}$ [Bour,, AC, VI, \S1, \no2, th\'eo 1]; et du fait que le corps r\'esiduel $k(sp(x))=R_{x}/\mathfrak{m}_{R_{x}}$ de $R_{x}$ est un quotient de $A_{0}$, celui de $\mathcal{V}_{x}$ en est un aussi. De plus on sait [B3, (0.2.2) ] ou [B-G-R, \S6.1.2, cor 3 du theo 1] que $K_{x}$ est une extension alg\'ebrique finie de $K$, et on vient de voir qu'elle est non ramifi\'ee puisque $\pi$ est une uniformisante de $\mathcal{V}_{x}$. D'o\`u $(vii).$\\

R\'eciproquement montrons que $(vii)\Longrightarrow(ii)$. D'apr\`es $(vii)$ le corps r\'esiduel $k(\mathcal{V}_{x})=\mathcal{V}_{x}/\pi\mathcal{V}_{x}$ de $\mathcal{V}_{x}$ est un quotient de $A_{0}$; par cons\'equent le morphisme
$$
\overline{i}_{x}: R_{x}/\pi R_{x}\rightarrow \mathcal{V}_{x}/\pi\mathcal{V}_{x}
$$
est une surjection entre deux $k$-espaces vectoriels de m\^eme dimension $n=[K_{x}:K]$, donc c'est un isomorphisme. Par suite [$EGA\ 0_{I}, (6.6.21)$]
$$
i_{x}: R_{x}\hookrightarrow \mathcal{V}_{x}
$$
est un isomorphisme; donc $\pi$ engendre $\mathfrak{m}_{R_{x}}$, d'o\`u $(ii).$\\

L'\'equivalence de $(vii)$ et $(vii)'$ est tautologique.\\

Clairement $(viii)\Longrightarrow(vii)$. R\'eciproquement supposons $(vii)$; alors $\mathcal{V}_{x}$ est non ramifi\'e sur $\mathcal{V}$. Or on a vu ci-dessus dans la d\'emonstration de  l'\'equivalence de $(ii)$ et $(vii)$ que la propri\'et\'e $(vii)$ implique que 
$$
i_{x}: R_{x}\hookrightarrow \mathcal{V}_{x}
$$
est un isomorphisme, donc $\mathcal{V}_{x}$ est un quotient de $\hat{A}$ au m\^eme titre que $R_{x}=\hat{A}/\overset{\circ}{\mathfrak{m}}_{x}$; d'o\`u $(viii)$.\\  

Montrons $(viii)\Longrightarrow(viii)'$. Comme $\mathcal{V}_{x}$ est un quotient de $\hat{A}$, l'injection
$$
\widetilde{i}_{x}: k(R_{x})\hookrightarrow k(\mathcal{V}_{x})
$$
est une surjection, donc c'est un isomorphisme.\\

R\'eciproquement montrons que $(viii)'\Longrightarrow(viii)$. Puisque $\mathcal{V}_{x}$ est  une  extension non ramifi\'ee de $\mathcal{V}$, l'injection $ i_{x}:R_{x}\hookrightarrow\mathcal{V}_{x}$ est non ramifi\'ee [EGA IV, (17.3.3)(v)], et comme $\widetilde{i}_{x}$ est bijective par hypoth\`ese, l'injection
$$
i_{x}: R_{x}\hookrightarrow \mathcal{V}_{x}
$$
est surjective [SGA1, I, cor 7.5], d'o\`u $(viii)$.  \\ 

 L'\'equivalence de $(viii)$ et $(ix)$ provient du fait que $\mathcal{V}_{x}$ est un anneau de valuation qui domine l'anneau local int\`egre $R_{x}$ en ayant m\^eme corps des fractions.\\
 
 L'\'equivalence de  $(ix)$ et  $(ix)'$ provient de [EGA IV, (17.4.1)$d''$)] et [Se, chap I, \S2, prop 2] et celle de $(ix)'$ et  $(ix)''$ est claire.  $\square$\\

Comme nous allons le voir ci-dessous les points $x$ qui v\'erifient les propri\'et\'es \'equivalentes de (3.2.4.2)2) sont les points \`a bonne r\'eduction au sens suivant:
 \vskip 3mm
\noindent \textbf{D\'efinition (3.2.4.3)}. \textit{Soit $A$ une $\mathcal{V}$-alg\`ebre lisse, $\mathcal{X}_{K}:=Spm\ \hat{A}_{K}$ et $x\in \mathcal{X}_{K}$ correspondant \`a un morphisme $\hat{A}_{K}\twoheadrightarrow K_{x}:=\hat{A}_{K}/\mathfrak{m}_{x}$. On dit que $x$ a bonne r\'eduction si et seulement si il existe une $\mathcal{V}$-alg\`ebre lisse $R$ qui est un quotient de $\hat{A}$ et un $K$-isomorphisme $R_{K}:=R\otimes_{\mathcal{V}}K\simeq K_{x}$}.

\vskip 3mm
\noindent \textbf{Proposition (3.2.4.4)}. \textit{Avec les notations pr\'ec\'edentes on a:}
\begin{enumerate}
	\item[$(i)$]\textit{Si $x\in \mathcal{X}_{K}$ a bonne r\'eduction, alors:}
		\begin{enumerate}
		\item[]\textit{Il existe un $\mathcal{V}$-isomorphisme $R\simeq R_{x}:=\hat{A}/\overset{\circ}{\mathfrak{m}}_{x}$,
		\item[]$R_{x}$ est une $\mathcal{V}$-alg\`ebre int\`egre finie \'etale,
		\item[]$\overline{sp}(x)=r\acute{e}d(x)$.}
		\end{enumerate}
	\item[$(ii)$]\textit{R\'eciproquement, si $\overline{sp}(x)=r\acute{e}d(x)$ alors $x$ a bonne r\'eduction.}
\end{enumerate}

\vskip 3mm
\noindent\textit{D\'emonstration}.\\

$(i)$ Supposons que $x\in \mathcal{X}_{K}$ a bonne r\'eduction et choisissons une pr\'esentation $R\simeq \mathcal{V}[X_{1},...,X_{m}]/\mathfrak{J}$ de la $\mathcal{V}$-alg\`ebre lisse $R$. Comme $R$ est plate sur $\mathcal{V}$, elle s'injecte dans $R_{K}\simeq K_{x}$, donc $R$ est int\`egre; de plus $R$ \'etant normal [EGA IV, (17.5.7)] il en r\'esulte que $R$ est int\'egralement clos de corps des fractions la $K$-alg\`ebre finie $K_{x}$. Notons $x_{i}\in R\hookrightarrow K_{x}$ l'image de $X_{i}$ dans $R$: $x_{i}$ est entier sur $K$, donc il existe un entier $r_{i}\in \mathbb{N}$ et des \'el\'ements $\lambda_{i,j}\in K$ pour $j\in \llbracket 0,r_{i}\rrbracket$ tels que $\lambda_{i,r_{i}}=1$ et 
$$
\sum_{j=0}^{j=r_{i}} \lambda_{i,j}x_{i}^{j}=0.
\leqno{ \bigstar}
$$	
En multipliant par une puissance ad hoc de $\pi$ il existe pour tout $i$ un entier $d_{i}\in \mathbb{N}$ tel que les $\lambda_{i,j}\pi^{d_{i}},\ j\in \llbracket 0,r_{i}\rrbracket$, soient des \'el\'ements de $\mathcal{V}$	premiers entre eux. Si $\lambda_{i,0}\pi^{d_{i}}=:\mu_{i,0}$ \'etait le seul des $\lambda_{i,j}\pi^{d_{i}},\ j\in \llbracket 0,r_{i}\rrbracket$, \`a ne pas \^etre multiple de $\pi$, en r\'eduisant modulo $\pi$ on aurait $\overline{\mu_{i,0}}=\overline{0}$, ce qui est absurde. Donc au moins un des $\lambda_{i,j}\pi^{d_{i}},\ j\in \llbracket 1,r_{i}\rrbracket$, est inversible: en r\'eduisant modulo $\pi$ l'\'equation $\bigstar$ on en d\'eduit que $\overline{x}_{i}:=x_{i}\ mod\ \pi$ est entier sur $k$. En faisant cela pour chaque $i$ il en r\'esulte que la $k$-alg\`ebre de type fini $R/\pi R$ est enti\`ere, donc elle est finie. Or $R_{K}\simeq K_{x}$ est fini sur $K$, donc $R$ est quasi-fini sur $\mathcal{V}$ [R, chap IV, prop 3]. En notant $B$ la fermeture int\'egrale de $\mathcal{V}$ dans $R$ on sait  alors par Raynaud [R, chap IV, cor 2 du th\'eo 1] que $Spec\ R\hookrightarrow Spec\ B$ est une immersion ouverte. Puisque $R$ est int\'egralement clos de corps des fractions $K_{x}$, $B$ est aussi la fermeture int\'egrale de $\mathcal{V}$ dans $K_{x}$, donc $B$ est un anneau de valuation discr\`ete [Se, chap II, \S2, prop 3], c'est l'anneau de valuation $B=\mathcal{V}_{x}$ de $K_{x}$ et c'est aussi une $\mathcal{V}$-alg\`ebre finie [loc cit].\\
Le morphisme compos\'e
$$
Spec\ R/\pi R\hookrightarrow Spec\ B/\pi B\rightarrow Spec\ k
$$
\'etant fini \'etale, l'immersion ouverte
$$
Spec\ R/\pi R\hookrightarrow Spec\ B/\pi B
$$
est finie, c'est donc une immersion ferm\'ee; par suite $R/\pi R$ est un quotient de l'alg\`ebre locale $B/\pi B$ [Bour, A I, \S8, prop 5 d)] et c'est aussi un anneau local [loc cit]. De plus le morphisme fini $Spec\ B \rightarrow Spec\ \mathcal{V}$ envoie le point ferm\'e $\{\mathfrak{m}_{B}\}=\{\pi_{B}B\}$ sur le point ferm\'e $\{\mathfrak{m}_{\mathcal{V}}\}=\{\pi \mathcal{V}\}$, i.e. $\mathfrak{m}_{B}\cap\mathcal{V}=\mathfrak{m}_{\mathcal{V}}$, d'o\`u $\pi\in\mathfrak{m}_{\mathcal{V}}$, ce qui fournit une \'egalit\'e entre id\'eaux de 
$B$, $<\pi>=<\pi_{B}^{e}>$ pour un certain entier $e\in\mathbb{N}$. Ainsi $B/\pi B=B/\pi_{B}^{e}B$ est artinien local, donc $R/\pi R$ est \'egalement artinien local. En particulier l'immersion ouverte
$$
Spec\ R/\pi R=\{\mathfrak{P}\}\hookrightarrow Spec\ B/\pi B=\{\mathfrak{P}'\}
$$
est surjective, c'est donc un isomorphisme, d'o\`u un isomorphisme:
$$
B/\pi B\overset{\sim}{\rightarrow} R/\pi R.
\leqno{\bigstar\bigstar}
$$
Or $R$, muni de la topologie $\pi$-adique, est un sous-espace topologique de l'espace topologique $K_{x}$ qui est s\'epar\'e (et complet) pour la topologie $\pi$-adique, donc $R$ est s\'epar\'e pour la topologie $\pi$-adique. De plus, comme $R$ est plat sur $\mathcal{V}$ et que $B$ est s\'epar\'e et complet pour la topologie $\pi$-adique en tant que $\mathcal{V}$-alg\`ebre finie, l'isomorphisme $\bigstar\bigstar$ nous permet de conclure que l'injection
$$
B\hookrightarrow R
$$  
est un isomorphisme [$EGA\ 0_{I},\ (6.6.21)$]. En particulier $R$ est un anneau de valuation discr\`ete et une $\mathcal{V}$-alg\`ebre finie \'etale: ainsi $\pi$ est une uniformisante de $R$. Comme la fl\`eche compos\'ee
$$
\hat{A}\twoheadrightarrow R_{x}\hookrightarrow B=R=\mathcal{V}_{x}
$$
est surjective, on obtient l'isomorphisme 
$$
R_{x}\overset{\sim}{\rightarrow}\mathcal{V}_{x}\ ;
$$
d'o\`u $\overline{sp}(x)=r\acute{e}d(x)$ d'apr\`es (3.2.4.2)$(ix)'$, ce qui d\'emontre (i).\\

$(ii)$ R\'eciproquement si $\overline{sp}(x)=r\acute{e}d(x)$, alors $R_{x}$ est un anneau de valuation discr\`ete fini \'etale sur $\mathcal{V}$ d'apr\`es (3.2.4.2)$(ix)''$, donc $x$ a bonne r\'eduction.\ ${\square}$

\vskip 3mm
\noindent \textbf{Proposition (3.2.4.5)}. \textit{Sous les hypoth\`eses et notations de} (3.2.4.2) \textit{les propri\'et\'es suivantes sont \'equivalentes:}
\begin{enumerate}
	\item[$(i)$]$R_{x}=\mathcal{V}_{x}.$
	\item[$(i)'$]$\mathcal{V}_{x}\ est \ un\ quotient\ de\ \hat{A}.$
	\item[$(ii)$]$R_{x}\ est \ un \ anneau \ de \ valuation\ discr\grave{e}te.$
	\item[$(iii)$]$L'injection \ i_{x}:R_{x}\hookrightarrow\mathcal{V}_{x}\ est\ non\ ramifi\acute{e}e\ et \ k(R_{x})=k(\mathcal{V}_{x}).$ 
	\item[$(iv)$]$Le\ morphisme \ i_{x}^{\#}:Spec\ \mathcal{V}_{x} \rightarrow Spec\ R_{x}\ est\ non\ ramifi\acute{e}\ et \ radiciel.$ 
\end{enumerate}

\vskip 3mm
\noindent\textit{D\'emonstration}. L'\'equivalence de $(i)$ et $(i)'$ est claire et celle de $(i)$ et $(iii)$ r\'esulte de [SGA1, I, cor 7.5]. L'\'equivalence de $(i)$ et $(ii)$ provient du fait que $\mathcal{V}_{x}$ est un anneau de valuation qui domine l'anneau local int\`egre $R_{x}$ et a m\^eme corps des fractions que lui. L'\'equivalence de $(iii)$ et $(iv)$ r\'esulte de [$EGA\ 0_{IV}$, (23.2.2)] puisque $k$ est un corps parfait.\  
$\square$\\

\vskip 3mm
\noindent \textbf{Corollaire (3.2.4.6)}. \textit{Sous les hypoth\`eses et notations de} (3.2.4.2) $soit\\
 x_{0}\in\vert Spec\ A_{0}\vert \ un\ point\ ferm\acute{e}\ de\ Spec\ A_{0}$. \textit{Alors les propri\'et\'es suivantes sont \'equivalentes:}
\begin{enumerate}
	\item[$(i)$]$x\in r\acute{e}d^{-1}(x_{0}):=\{x\in\mathcal{X}_{K}/r\acute{e}d(x)=x_{0}\}.$
	\item[$(ii)$]$R_{x}\ est \ un \ anneau \ de \ valuation\ discr\grave{e}te\ d'uniformisante\ \pi\ et\\
	 k(R_{x})=k(x_{0}).$
	 \item[$(iii)$]$R_{x}\ est\ non \ ramifi\acute{e}\ sur\ \mathcal{V}\ et\  k(R_{x})= k(x_{0}).$
	\item[$(iv)$]$\mathcal{V}_{x} \ est\ non \ ramifi\acute{e}\ sur\ \mathcal{V} \ et\ k(\mathcal{V}_{x})=k(x_{0}).$
	\item[$(v)$]$K_{x}\  est\ non \ ramifi\acute{e}\ sur \ K\  et\  \widetilde{K_{x}}= k(x_{0}).$
\end{enumerate}
\vskip 3mm
\noindent\textit{D\'emonstration}. R\'esulte de (3.2.4.2) et (3.2.4.5). \ $\square$\\

\textbf{Description des tubes via des morphismes}
\vskip 3mm
\noindent \textbf{Proposition (3.2.4.7)}. \textit{Sous les hypoth\`eses et notations de (3.2.4.2) soit
 $x_{0}\in\vert Spec\ A_{0}\vert$  un point ferm\'e de $Spec\ A_{0}$. Alors} 
 \begin{enumerate}
\item[$(i)$]$r\acute{e}d^{-1}(x_{0})$\textit{ s'identifie \`a l'ensemble (non vide) des rel\`evements (n\'eces
sairement surjectifs) $\hat{A}\twoheadrightarrow \mathcal{V}(x_{0})$ du  morphisme $A_{0}\twoheadrightarrow k(x_{0})$ donn\'e avec $x_{0}$};
  \item[$(ii)$]\textit{En particulier   $\overline{sp}:\  Max\ \hat{A}_{K} \rightarrow  Max\ A_{0}$ est  surjectif. 
 \item[$(iii)$]Parmi ces  rel\`evements du (i) il en existe un et  un seul qui  commute aux Frobenius,  c'est le  rel\`evement  de Teichm\" uller:} 
\end{enumerate}
$$
\begin{array}{c}
Teich_{\hat{A}_{K}}(x_{0})\in r\acute{e}d^{-1}(x_{0})\subset\ ]x_{0}[:=\overline{sp}^{-1}(x_{0})\\
\\
r\acute{e}d\ (Teich_{\hat{A}_{K}}(x_{0}))=\overline{sp}\ (Teich_{\hat{A}_{K}}(x_{0}))=x_{0}\ .
\end{array}
$$

\vskip 3mm
\noindent\textit{D\'emonstration}. Si $x\in r\acute{e}d^{-1}(x_{0})$, alors $\hat{A}\twoheadrightarrow R_{x}$ est un rel\`evement de $A_{0}\twoheadrightarrow k(x_{0})$ et il existe un unique $\mathcal{V}$-isomorphisme $R_{x}\simeq \mathcal{V}(x_{0})$ relevant l'identit\'e de $k(x_{0})$ [EGA IV, (18.3.2)]. R\'eciproquement la lissit\'e formelle de $\hat{A}$ sur $\mathcal{V}$ nous assure de l'existence de tels rel\`evements $\hat{A}\twoheadrightarrow \mathcal{V}(x_{0})$ du  morphisme $A_{0}\twoheadrightarrow k(x_{0})$; alors l'image $R_{x}$ de $\hat{A}$ dans $\mathcal{V}(x_{0})$ est un anneau de valuation discr\`ete d'uniformisante $\pi$ et $k(R_{x})=k(x_{0})$, donc $x\in r\acute{e}d^{-1}(x_{0})$ d'apr\`es (3.2.4.6).\\
Le reste de la proposition est clair car $x=Teich_{\hat{A}_{K}}(x_{0})$ v\'erifie les propri\'et\'es de (3.2.4.6) (cf aussi 3.2.2).  $\square$\\

Pour $ x_{0}\in\vert Spec\ A_{0}\vert$  un point ferm\'e de $Spec\ A_{0}$, nous allons donner une description du tube $]x_{0}[:=\overline{sp}^{-1}(x_{0})$ qui fournira encore la surjectivit\'e de $\overline{sp}$ en exhibant des ant\'ec\'edents de $\overline{sp}$ en dehors de $r\acute{e}d^{-1}(x_{0})$.\\

Soit donc $ x_{0}\in\vert Spec\ A_{0}\vert$. Si $x\in ]x_{0}[$, on a vu que $R_{x}=\hat{A}/\overset{\circ}{\mathfrak{m}}_{x}$ est une $\mathcal{V}$-alg\`ebre int\`egre finie et plate (ce qui implique que $R_{x}$ est local et $Spec\ R_{x}=\{(0),\mathfrak{m}_{R_{x}}\}$) \`a corps r\'esiduel $R_{x}/\mathfrak{m}_{R_{x}}=k(x_{0}).$\\
Inversement, soit $\mathcal{V}'$ une $\mathcal{V}$-alg\`ebre int\`egre finie et plate, donc locale, de corps r\'esiduel $k(x_{0})$. Notons que sur $\mathcal{V}'$ les topologies $\pi$-adiques et $\mathfrak{m}_{\mathcal{V}'}$ co\"incident puisque l'id\'eal maximal de l'anneau local artinien $\mathcal{V}'/\pi \mathcal{V}'$ est nilpotent. Par lissit\'e formelle de $\hat{A}$ sur $\mathcal{V}$, il existe (au moins) un morphisme $\hat{A}\rightarrow \mathcal{V}'$ rendant commutatif le diagramme
$$
\xymatrix
{
A_{0} \ar@{->>}[r]&k(x_{0})\\
\hat{A}\ar[r]\ar@{->>}[u]&\mathcal{V}' \ar@{->>}[u]&.
}
$$
L'image $R'$ de $\hat{A}$ dans $\mathcal{V}'$ est elle aussi une $\mathcal{V}$-alg\`ebre int\`egre finie et plate, donc locale: de plus le corps r\'esiduel $k(R')=R'/\mathfrak{m}_{R'}$ est isomorphe \`a $k(x_{0})=\mathcal{V}'/\mathfrak{m}_{\mathcal{V}'}$ car dans le diagramme commutatif
$$
\xymatrix
{
A_{0} \ar@{->>}[r]&R'/\mathfrak{m}_{R'}\ \ar@{^{(}->}[r]&k(x_{0})=\mathcal{V}'/\mathfrak{m}_{\mathcal{V}'}\\
A_{0}\ar@{->>}[r]\ar@{=}[u]&R'/\pi R'\ar[r]\ar@{->>}[u]&\mathcal{V}'/\pi \mathcal{V}'\ar@{->>}[u]\\
\hat{A}\ar@{->>}[r]\ar@{->>}[u]&R'  \ \ar@{^{(}->}[r]\ar@{->>}[u]&\mathcal{V}'  \ar@{^{(}->} [u]
}
$$
la fl\`eche compos\'ee horizontale du haut est surjective par hypoth\`ese. Ainsi on a prouv\'e:

\vskip 3mm
\noindent \textbf{Proposition (3.2.4.8)}. \textit{Sous les hypoth\`eses et notations de (3.2.4.2) soit
 $x_{0}\in\vert Spec\ A_{0}\vert $ un point ferm\'e de $Spec\ A_{0}$. Alors:} 
 \begin{enumerate}
\item[]$]x_{0}[=sp^{-1}(x_{0})$ \textit{s'identifie \`a l'ensemble (non vide) des morphismes\\
 $\hat{A}\rightarrow\mathcal{V}'$ qui factorisent la surjection $\hat{A}\twoheadrightarrow k(x_{0}$) donn\'ee avec $x_{0}$,
  \item[] o\`u $\mathcal{V}'$ parcourt l'ensemble des  $\mathcal{V}$-alg\`ebres int\`egres finies et plates \`a corps r\'esiduel $k(x_{0})$.}
\end{enumerate}
L'\'etude de la boule unit\'e ferm\'ee de la droite affine rigide que nous allons entreprendre au \S3.3 nous fournira ensuite au \S3.4 des exemples explicites de points
$$
x\in ]x_{0}[\ \setminus\  r\acute{e}d^{-1}(x_{0})
$$ 
qui sont dans le tube de $x_{0}$ sans \^etre dans $r\acute{e}d^{-1}(x_{0})$.

\vskip 3mm
\textbf{Description des tubes via des \'equations}\\

Soit $ x_{0}\in\vert Spec\ A_{0}\vert$  un point ferm\'e de $Spec\ A_{0}$, $x=Teich_{\hat{A}_{K}}(x_{0})$ et $f_{1},...,f_{r}$ des g\'en\'erateurs de l'id\'eal (de type fini) $\overset{\circ}{\mathfrak{m}}_{x}$ de $\hat{A}$. Puisque 
$$
\mathfrak{m}_{x_{0}}=\mathfrak{m}_{\overline{sp}(x)}=\rho(\overset{\circ}{\mathfrak{m}}_{x})
$$
on a la description suivante du tube $]x_{0}[:=\overline{sp}^{-1}(x_{0})$ de $x_{0}$\ [B3, (1.1.1)]:
$$
]x_{0}[=\{z\in\mathcal{X}_{K}/\forall i=1,...,,r,\ \vert f_{i}(z)\vert<1\} \ .
$$
Pour $z=x=Teich_{\hat{A}_{K}}(x_{0})$ on a $f_{i}(x):=f_{i} \ mod\ \overset{\circ}{\mathfrak{m}}_{x}=0$; donc non seulement on a prouv\'e que
$$
x=Teich_{\hat{A}_{K}}(x_{0})\in ]x_{0}[=\{z\in\mathcal{X}_{K}/\forall i=1,...,,r,\ \vert f_{i}(z)\vert<1\} \ ,
$$
ce qui est vrai pour toute section de $\overline{sp}$, mais m\^eme que $x=Teich_{\hat{A}_{K}}(x_{0})$ est \guillemotleft au centre\guillemotright \ du tube de $]x_{0}[$, i.e. que 
$$
\forall i=1,...,r \ f_{i}(x)=0\ .
$$
Plus g\'en\'eralement, cette propri\'et\'e s'\'etend \`a tous les $y\in r\acute{e}d^{-1}(x_{0})$ \`a condition de changer de g\'en\'erateurs pour $\overset{\circ}{\mathfrak{m}}_{y}$: en effet, pour $g_{1},...,g_{s}$ des g\'en\'erateurs de $\overset{\circ}{\mathfrak{m}}_{y}$, on a:
$$
]x_{0}[=\{z\in\mathcal{X}_{K}/\forall i=1,...,,r,\ \vert f_{i}(z)\vert<1\} =\{z\in\mathcal{X}_{K}/\forall j=1,...,,s,\ \vert g_{j}(z)\vert<1\} 
$$
et
 $$
 \forall j=1,...,s, \ g_{j}(y):=g_{j} \ mod\  \overset{\circ}{\mathfrak{m}}_{y}=0\ .
$$
On peut formaliser cette constatation dans la proposition suivante:
\vskip 3mm
\noindent \textbf{Proposition (3.2.4.9)}. \textit{Sous les hypoth\`eses et notations de} (3.2.4.2) $soit\ x_{0}\in\vert Spec\ A_{0}\vert \ un\ point\ ferm\acute{e}\ de\ Spec\ A_{0}$. \textit{Alors}
$$
\begin{array}{ccc}
r\acute{e}d^{-1}(x_{0})=& \{x\in \mathcal{X}_{K}/ \exists J_{x}\subset\mathbb{N}\ fini,\ \forall i \in J_{x}\ \exists f_{i}\in \hat{A},\ tel\ que\\
\\
 &x_{0}=V(\overline{f}_{i},\ i\in J_{x})\ o\grave{u}\  \overline{f}_{i}:=f_{i}\ mod \ \pi\hat{A}\\
 \\
  & et\ x\in V(f_{i},\ i\in J_{x})\}\ .
\end{array}
$$
\vskip 3mm
\noindent\textit{D\'emonstration}. En notant $\mathcal{E}$ l'ensemble second membre dans l'\'egalit\'e ci-dessus, on a d\'ej\`a prouv\'e que $r\acute{e}d^{-1}(x_{0})\subset \mathcal{E}$. R\'eciproquement, soit $x\in\mathcal{E}$. Les $\overline{f}_{i},\ i\in J_{x}$, forment par hypoth\`ese un syst\`eme de g\'en\'erateurs de l'id\'eal $\mathfrak{m}_{x_{0}}$ de $A_{0}$ tels que $f_{i}(x):=f_{i}\ mod\ \overset{\circ}{\mathfrak{m}}_{x}=0$ donc, pour tout $i\in J_{x}$, $f_{i}\in \overset{\circ}{\mathfrak{m}}_{x}$ et $x\in ]x_{0}[=\{z\in \mathcal{X}_{K}/ \forall i\in J_{x},\ \vert f_{i}(x)\vert <1\}$. En particulier $\overline{sp}(x)=x_{0}$ et 
$$
\mathfrak{m}_{x_{0}}\subset \rho (\overset{\circ}{\mathfrak{m}}_{x})=:\widetilde{\mathfrak{m}}_{x}\ ;
$$
 comme on a toujours l'inclusion
 $$
 \widetilde{\mathfrak{m}}_{x}=\rho (\overset{\circ}{\mathfrak{m}}_{x})\subset \mathfrak{m}_{\overline{sp}(x)}=\mathfrak{m}_{x_{0}}
 $$
 on a l'\'egalit\'e $\widetilde{\mathfrak{m}}_{x}=\mathfrak{m}_{x_{0}}$, i.e. $r\acute{e}d(x)=x_{0}$\ . 
 $\square$ \\
 
De mani\`ere imag\'ee on peut dire que les $x\in r\acute{e}d^{-1}(x_{0})$ sont \guillemotleft au centre\guillemotright\ du tube $]x_{0}[$, au sens o\`u, pour un choix judicieux de g\'en\'erateurs $f_{i}$ de $\mathfrak{m}_{x}$, on a $f_{i}(x)=0$, alors que 
$$
]x_{0}[=\{z\in\mathcal{X}_{K}/\forall i=1,...,r,\ \vert f_{i}(z)\vert<1\} \ .
$$
 
 \subsection*{3.3. $Spm\ K\{T\}$: la boule unit\'e de la droite affine rigide}
 \subsubsection*{3.3.1. Les points de $Spm\ K\{T\}$}
 Rappelons pour m\'emoire le lemme \'evident suivant:
 \vskip 3mm
\noindent \textbf{Lemme (3.3.1.1)}. \textit{Les \'el\'ements inversibles de $\mathcal{V}\{T\}$ sont les \'el\'ements de 
$$
\mathcal{V}\{T\}^{\ast}=U+\pi\mathcal{V}\{T\}
$$
 o\`u $U=\mathcal{V}\setminus\pi\mathcal{V}$ est l'ensemble des \'el\'ements inversibles de $\mathcal{V}$.}\\

 Nous avons vu [B3, (0.2.2)] que les points de $Spm\ K\{T\}$ sont en bijection avec les quotients $\mathcal{V}\{T\}/I$ qui sont des $\mathcal{V}$-alg\`ebres int\`egres finies et plates.\\
 Par exemple, soit $P(T)\in \mathcal{V}[T]$ un polyn\^ome unitaire irr\'eductible dans $\mathcal{V}[T]$ et posons $I=P(T)\mathcal{V}\{T\},\ \widetilde{I}=P(T)\mathcal{V} [T]$; par factorialit\'e de $\mathcal{V}[T]$, $\widetilde{I}$ est premier et  la fl\`eche naturelle
$$
\widetilde{R}:=\mathcal{V}[T]/ \widetilde{I}\overset{\sim}{\longrightarrow}\mathcal{V}\{T\}/I:=R
$$
est un isomorphisme de $\mathcal{V}$-alg\`ebres finies int\`egres; et puisque $P(T)$ est unitaire, $\widetilde{R}$ est plate sur $\mathcal{V}$ [Mi, chap I, Rk 2.6 a)], donc $R$ correspond \`a un point de $Spm\ \mathcal{V}\{T\}$. Nous allons montrer ci-dessous qu'en fait tous les points de $Spm\ K\{T\}$ sont de cette forme.\\

 Consid\'erons en effet un id\'eal $I\subset \mathcal{V}\{T\}$ tel que $\mathcal{V}\{T\}/I$ soit une $\mathcal{V}$-alg\`ebre int\`egre finie et plate.; notons $\overline{I}$ l'image canonique de $I$ dans $k[T]$. Par platitude de $\mathcal{V}\{T\}/I$ sur $\mathcal{V}$, on a $Tor_{1}^{\mathcal{V}}(\mathcal{V}\{T\}/I,\ \mathcal{V}/\pi\mathcal{V})=0$, donc on a un isomorphisme
 $$
 I\otimes_{\mathcal{V}}k=I/\pi I\overset{\sim}{\rightarrow} \overline{I}\subset k[T],
 $$
  ce qui permet d'identifier $I/\pi I$ \`a son image $\overline{I}$ dans $k[T]$ qui est un id\'eal principal: on choisit pour g\'en\'erateur de l'id\'eal $\overline{I}$ le g\'en\'erateur unitaire $\overline{f}=\sum_{i=0}^{i=d}\overline{a}_{i}T^{i},\  \overline{a}_{i}\in k, \overline{a}_{d}=1$. Via Nakayama, puisque $Rad(\mathcal{V}\{T\})=\pi\ \mathcal{V}\{T\} $ [(3.2.4.1)(ii)], que $I$ est de type fini (par noeth\'erianit\'e de $\mathcal{V}\{T\}$ [Sal] ) et que $\overline{I}$ est principal engendr\'e par  $\overline{f}$, $I$ est lui aussi principal, engendr\'e par un $g\in I\subset \mathcal{V}\{T\}$ relevant $\overline{f}$. On \'ecrit $g(T)=\pi^{r}h(T)$, o\`u $r\in \mathbb{N}$ et la norme de Gau\ss \ $\vert\vert h\vert\vert$ de $h$ est \'egale \`a 1 ([B-G-R, \S1.4.1] si $h(T)=\sum_{i=0}^{+\infty}\alpha_{i}T^{i}\in \mathcal{V}\{T\}, \lim_{i \to \infty}\alpha_{i} =0,\  \vert\vert h\vert\vert:=max \vert \alpha_{i}\vert$), et puisque $\overline{f}\not= 0$ on a $r=0$. Alors [B-G-R, 5.2.1, def 1] $g(T)=h(T)$ est $T$-distingu\'e d'un certain degr\'e not\'e $s$; par le th\'eor\`eme de pr\'eparation de Weierstra\ss  \ [B-G-R, 5.2.2, theo 1] il existe un unique \'el\'ement inversible $e(T)\in K\{T\}$ et un unique polyn\^ome unitaire $\omega(T)\in \mathcal{V}[T]$ de degr\'e $s$ tels que $g(T)=e(T)\omega (T)$. Puisque $\vert\vert \omega \vert\vert=1$ et $1=\vert\vert g \vert\vert=\vert\vert e \vert\vert\times\vert\vert \omega \vert\vert$ [B-G-R, 5.1.2,  prop 1] on a $\vert\vert e \vert\vert=1$ et en fait $e(T)\in \mathcal{V}\{T\}$: si $e'(T)$ est un inverse de $e(T)$ dans $K\{T\}$ l'\'egalit\'e $e(T)\times e'(T)=1$ fournit $\vert\vert e \vert\vert\times\vert\vert e' \vert\vert=\vert\vert e' \vert\vert=1$, donc $e'(T)\in \mathcal{V}\{T\}$ et $e(T)$ est inversible dans $\mathcal{V}\{T\}$. Par le lemme (3.3.1.1) l'\'el\'ement $e(T)$ est de la forme $e(T)=u+\pi T v(T)$ avec $u\in \mathcal{V}\setminus\pi \mathcal{V}$ et $v(T)\in \mathcal{V}\{T\}$; or la r\'eduction $mod\ \pi$ de $g$ est le polyn\^ome unitaire $\overline{f}$, donc $e(T)\ mod\ \pi=1$ et $\overline{f}(T)=\omega (T)\ mod\ \pi$: en particulier $s=d$ et par la factorialit\'e de $\mathcal{V}\{T\}$  [Sal] $\omega(T)$ est irr\'eductible dans $\mathcal{V}\{T\}$. Or les g\'en\'erateurs de l'id\'eal principal $I=\omega (T)\mathcal{V}\{T\}$ sont tous de la forme $e_{1}(T)\times\omega (T)$ o\`u $e_{1}(T)$ est inversible dans $\mathcal{V}\{T\}$ donc de la forme $e_{1}(T)=u_{1}(1+\pi v_{1}(T))$ avec $u_{1}\in\mathcal{V}\setminus\pi\mathcal{V},\ v_{1}(T)\in\mathcal{V}\{T\}$; par suite\\
   \textbf{$\omega(T)$ est l'unique g\'en\'erateur irr\'eductible dans $\mathcal{V}\{T\}$ de I qui soit un polyn\^ome unitaire de degr\'e $d$=degr\'e de $\overline{f}$, relevant $\overline{f}$, o\`u $\overline{f}$ est l'unique g\'en\'erateur unitaire de $\overline{I}=I\ mod\ \pi$: $\omega(T)$ est m\^eme irr\'eductible dans $\mathcal{V}[T]$} d'apr\`es le lemme suivant:
 \vskip 3mm
\noindent \textbf{Lemme (3.3.1.2)}. \textit{Soit $\omega(T)\in\mathcal{V}[T]$, unitaire. Les propri\'et\'es suivantes sont \'equivalentes:}
\begin{enumerate}
\item[$(i)$] $\omega (T)$ \textit{est irr\'eductible dans $\mathcal{V} [T] $}.
\item[$(ii)$] $\omega (T)$ \textit{est irr\'eductible dans $K [T] $.
\item[$(iii)$] $\omega (T)$ est irr\'eductible dans $\mathcal{V}\{ T\}$.}
\end{enumerate}

 \vskip 3mm
\noindent\textit{D\'emonstration}. \\
L'\'equivalence de $(i)$ et $(ii)$ est classique puisque les coefficients de $\omega(T)$ sont premiers entre eux.\\ 
Montrons l'\'equivalence $(i)\iff (iii)$. Soit $I=\omega(T)\mathcal{V}\{T\},\ \widetilde{I}=\omega(T)\mathcal{V} [T]$; puisque $\omega(T)$ est unitaire la fl\`eche naturelle
$$
\widetilde{R}:=\mathcal{V}[T]/ \widetilde{I}\overset{\sim}{\longrightarrow}\mathcal{V}\{T\}/I:=R
\leqno{\bigstar}
$$
est un isomorphisme. Or la factorialit\'e de $\mathcal{V}[T]$ fournit l'\'equivalence
$$
\omega(T)\  \mbox{irr\'eductible dans}\ \mathcal{V}[T]\iff \widetilde{I}\ \mbox{est premier} \ ;
$$
de m\^eme la factorialit\'e de $\mathcal{V}\{T\}$  [Sal] fournit l'\'equivalence
$$
\omega(T)\  \mbox{irr\'eductible dans}\ \mathcal{V}\{T\}\iff I \ \mbox{est premier} \ .
$$
D'o\`u l'\'equivalence du lemme gr\^ace \`a l'isomorphisme $\bigstar$ ci-dessus. $\square$\\

Ainsi on a prouv\'e:
 \vskip 3mm
\noindent \textbf{Proposition (3.3.1.3)}. \textit{Les points de $Spm\ K\{T\}$ sont en bijection avec les polyn\^omes unitaires irr\'eductibles de $\mathcal{V}[T]$}.\\

Rappelons le th\'eor\`eme suivant de Robert:
 \vskip 3mm
\noindent \textbf{Th\'eor\`eme (3.3.1.4) [Ro, chap 6, \S2.2, theo 2]}. \textit{Tout $S(T)\in K\{T\}$ s'\'ecrit de mani\`ere unique sous la forme 
$$
S(T)=\pi^{r}\times u\times (1+\pi v(T))\times \omega_{1}(T)\times \omega_{2}(T)
$$
o\`u}
\begin{enumerate}
\item[$(i)$]$r\in \mathbb{Z},\ u\in \mathcal{V}^{\ast}=\mathcal{V}\setminus\pi\mathcal{V}, \ v(T)\in\mathcal{V}\{T\}\ ,$
\item[$(ii)$]\textit{$\omega_{1}(T)\in \mathcal{V}[T]$, unitaire, $\vert\omega_{1}(0)\vert=1\ ,$
\item[$(iii)$]$\omega_{2}(T)\in \mathcal{V}[T]$, unitaire, de degr\'e $\mu,\ \omega_{2}(T)\equiv T^{\mu}\ mod\ \pi$ .}
\end{enumerate}

On en d\'eduit le corollaire suivant:
\vskip 3mm
\noindent \textbf{Corollaire (3.3.1.5)}.
\begin{enumerate}
\item[$(i)$]\textit{Les \'el\'ements inversibles de $K\{T\}$ sont les \'el\'ements de la forme
$$
S(T)=\pi^{r}\times u\times (1+\pi v(T))
$$
o\`u
$$
r\in \mathbb{Z},\ u\in \mathcal{V}^{\ast}=\mathcal{V}\setminus\pi\mathcal{V}, \ v(T)\in\mathcal{V}\{T\}\ .
$$
\item[$(ii)$]Aux inversibles pr\`es, les irr\'eductibles de $K\{T\}$ sont les polyn\^omes $\omega(T)\in \mathcal{V}[T]$, unitaires et irr\'eductibles dans $\mathcal{V}[T]$.
\item[$(iii)$]Aux inversibles pr\`es, les irr\'eductibles de $\mathcal{V}\{T\}$ sont de deux sortes: l'\'el\'ement $\pi$ et les polyn\^omes $\omega(T)\in \mathcal{V}[T]$, unitaires et irr\'eductibles dans $\mathcal{V}[T]$.}
\end{enumerate}
\vskip 3mm
\noindent\textit{D\'emonstration}. Les polyn\^omes unitaires irr\'eductibles sont de degr\'e sup\'erieur ou \'egal \`a un et $\pi$ est irr\'eductible dans  $\mathcal{V}\{T\}$, car le quotient $\mathcal{V}\{T\}/\pi \mathcal{V}\{T\}= k[T]$ est int\`egre.\  $\square$\\

\vskip 3mm
Ce qui permet de donner une autre description des points de $Spm\ K\{T\}$:
\vskip 3mm
\noindent \textbf{Corollaire (3.3.1.6)}. \textit{Les points de $Spm\ K\{T\}$ sont en bijection avec un syst\`eme de repr\'esentants des \'el\'ements irr\'eductibles de $K\{T\}$}.\\

\subsubsection*{3.3.2. Sp\'ecialisation, (bonne) r\'eduction et Teichm\" uller dans $Spm\ K\{T\}$}

Soit $f(T)\in\mathcal{V}[T]$, unitaire et irr\'eductible dans $\mathcal{V}[T]$, qu'on identifiera dans la suite \`a un point $x\in Spm\ K\{T\}$; par le lemme de Hensel [Bour, AC III, \S4, \no3, th\'eo 1] son image canonique $\overline{f}(T)\in k[T]$ dans $k[T]$ n'a qu'un seul facteur irr\'eductible not\'e 
$$
\overline{g}(T)=\sum_{i=0}^{i=r}a_{i}T^{i}\ , \ a_{i}\in k\ ,\ a_{r}=1\ ;\ \overline{f}(T)=\overline{g}(T)^{e}\ \mbox{pour}\ e\in\mathbb{N}^{\ast}\ .
$$
Au d\'ebut du \S3.1 on a d\'efini
$$
Teich(\overline{g})(T)=\sum_{i=0}^{i=r}b_{i}T^{i}
$$ 
appel\'e rel\`evement de Teichm\" uller de $\overline{g}$; puisque $\overline{g}$ est irr\'eductible, $Teich(\overline{g})$ est irr\'eductible, de m\^eme que tout rel\`evement $g\in \mathcal{V}[T]$ de $\overline{g}$. En vertu de la commutativit\'e du diagramme $(3.1.4)$ on a
$$
Teich_{K\{T\}}(\overline{g}k[T])=Teich(\overline{g})K\{T\}\ .
\leqno{(3.3.2.1)}
$$
On donne l'interpr\'etation de la sp\'ecialisation, de la r\'eduction et du rel\`evement de Teichm\" uller en termes de polyn\^omes dans la proposition \'evidente suivante:
\vskip 3mm
\noindent \textbf{Proposition (3.3.2.2)}. \textit{Sous les hypoth\`eses et notations 3.3.2 pr\'ec\'edentes on a:}
\begin{enumerate}
\item[$(i)$]
	\begin{enumerate}
	\item[]$ \mathfrak{m}_{x}=fK\{T\},\ \overset{\circ}{\mathfrak{m}}_{x}=f\mathcal{V}\{T\}\ ,$
	\item[]$\mathfrak{m}_{R_{x}}$ \textit{est engendr\'e par la classe de} $\pi \ et \ g \ $ \textit{dans} $R_{x}$\ \textit{pour tout rel\`evement} $g\ de\ \overline{g}$,
	 \item[]$\mathfrak{m}_{sp(x)}=\pi\mathcal{V}\{T\} +g\mathcal{V}\{T\}=\pi\mathcal{V}\{T\} +Teich(\overline{g})\mathcal{V}\{T\}\ ,$
	 \item[]$\mathfrak{m}_{\overline{sp}(x)}=\overline{g}k[T], \  \widetilde{\mathfrak{m}}_{x}=\overline{f}k[T]$, 
	\end{enumerate}
\item[$(ii)$]$\overline{sp}(f)=\overline{g}\ ,$
\item[$(iii)$]$r\acute{e}d(f)=\overline{f}\ ,$
\item[$(iv)$]$\overline{sp}(f)=r\acute{e}d(f)=\overline{f}$\ \textit{( i.e. f a bonne r\'eduction) $\iff$ $\overline{f}$ est irr\'eductible},
\item[$(v)$]$Teich_{K\{T\}}(\overline{g})=Teich(\overline{g})\ ,$
\item[$(vi)$]$\overline{sp}(Teich(\overline{g}))=r\acute{e}d(Teich(\overline{g}))=\overline{g}$ ,
\item[$(vii)$]$deg(f)=e\times deg(\overline{g})$; \textit{lorsque $R_{x}=\mathcal{V}_{x}, \ e$ est l'indice de ramification de $\mathcal{V}_{x}$ sur $\mathcal{V}$ et $deg(\overline{g})$ est le degr\'e r\'esiduel de  $K_{x}$ sur $K$} .
\end{enumerate}
\vskip 5mm
Ce qui, pour les tubes, fournit le th\'eor\`eme:

\vskip 3mm
\noindent \textbf{Th\'eor\`eme (3.3.2.3)}. \textit{Soit $x_{0}\in \ Max\ k[T]$; on note $\mathfrak{m}_{x_{0}} =\overline{g}(T)k[T]$ l'id\'eal maximal correspondant o\`u $\overline{g}(T)$ est le g\'en\'erateur unitaire irr\'eductible de $\mathfrak{m}_{x_{0}}$. Alors on a:}
\begin{enumerate}
\item[$(i)$]$r\acute{e}d^{-1}(x_{0})=\{ f\in \mathcal{V}[T]/ \ f \ unitaire \ tel\ que\ \overline{f}=\overline{g}  \}\ .$
\item[$(ii)$]$]x_{0}[:=\overline{sp}^{-1}(\overline{g})=\{ f\in \mathcal{V}[T]$/ \textit{f unitaire irr\'eductible},  $\exists e\in\mathbb{N},\ \overline{f}=\overline{g}^{e}  \}\ .$
\item[$(iii)$]$]x_{0}[ \setminus Êr\acute{e}d^{-1}(x_{0})=\{ f \in \mathcal{V}[T]$/ \textit{f unitaire irr\'eductible}, $\exists e\in\mathbb{N}^{\geqslant 2},\overline{f}=\overline{g}^{e}  \}\ .$
\item[$(iv)$]$Teich(\overline{g})\in r\acute{e}d^{-1}(x_{0})\ .$
\item[$(v)$]$f\in]x_{0}[ \setminus Êr\acute{e}d^{-1}(x_{0})$ \textit{si $\overline{g}(T)=T,f$ polyn\^ome d'Eisenstein, $\overline{f}=\overline{g}^{e}, e\geqslant 2$.}
\end{enumerate}

 \subsection*{3.4. Exemples}
Les exemples que nous allons traiter se rapportent aux points de $Spm\ K\{T\}$ pour diff\'erents corps $K$, extensions finies de $\mathbb{Q}_{p}$: comme dans le \S3.3 nous identifierons les points de $Spm\ K\{T\}$ ou $Max\ k[T]=Spm\ k[T]$ aux polyn\^omes unitaires irr\'eductibles auxquels ils correspondent. Nous utiliserons aussi les notations du \S3.2, $R_{x},\mathcal{V}_{x}$ avec pour corps r\'esiduels $k(R_{x}), k(\mathcal{V}_{x})$ et corps des fractions $K_{x}$, etc...\\
Pour illustrer la diversit\'e des situations nous allons exhiber tour \`a tour des exemples dans lesquels $R_{x}$ est un anneau de valuation discr\`ete ramifi\'e ou non sur $\mathcal{V}$, des cas o\`u $R_{x}$ n'est pas un anneau de valuation discr\`ete avec $\mathcal{V}_{x}$ ramifi\'e ou non sur $\mathcal{V}$ ou totalement ramifi\'e sur $\mathcal{V}$, en donnant des points dans $r\acute{e}d^{-1}(x_{0})$ et dans $ ]x_{0}[ \setminus Êr\acute{e}d^{-1}(x_{0})$.
\subsubsection*{3.4.1. \textbf{1\ier} cas: $R_{x}$ anneau de valuation discr\`ete non ramifi\'e sur $\mathcal{V}$}
$\rondI$
 $R_{x}=\mathcal{V}_{x}=\mathbb{Z}_{3}=\mathbb{Z}_{3}[T]/(T)$.\\
 
  Ici $K=\mathbb{Q}_{3},\ x_{0}=T\in \mathbb{F}_{3}[T] ,\ k(x_{0})=\mathbb{F}_{3} $ et $x= T\in r\acute{e}d^{-1}(x_{0})$.\\
 
$\rondII$
 $R_{x}=\mathcal{V}_{x}=\mathbb{Z}_{3}[\sqrt{2}].$
 \begin{enumerate}
 \item[]Ici $K=\mathbb{Q}_{3},\ K(R_{x})=K(\mathcal{V}_{x})=\mathbb{Q}_{3}[\sqrt{2}]$. 
 \item[]$R_{x}=\mathbb{Z}_{3}[T]/(T^2-2)$ est non ramifi\'e de degr\'e deux sur $\mathbb{Z}_{3}.$ 
 \item[]$ x_{0}=T^2-2\in\mathbb{F}_{3}[T]$.
 \item[]$k(x_{0})=\mathbb{F}_{9}=\mathbb{F}_{3}[T]/(T^2-2)=\mathbb{F}_{3}[T]/(T^2+1)=k(R_{x})$.
 \item[]$x= T^2-2\in r\acute{e}d^{-1}(x_{0}), y= T^2+1=Teich(x_{0})\in r\acute{e}d^{-1}(x_{0})$.
   \item[]$K(R_{x})=K(\mathcal{V}_{x})$ est une extension non ramifi\'ee de degr\'e $2$ de $K=\mathbb{Q}_{3}$.
     \item[]$k(R_{x})= \mathbb{F}_{9}=k(\mathcal{V}_{x})$.
 \end{enumerate}
\subsubsection*{3.4.2. \textbf{2\ieme} cas: $R_{x}$ anneau de valuation discr\`ete ramifi\'e sur $\mathcal{V}$}
$\rondI$ 
$R_{x}=\mathcal{V}_{x}=\mathbb{Z}_{2}[i].$
\begin{enumerate}
 \item[]Ici $K=\mathbb{Q}_{2},\ K(R_{x})=K(\mathcal{V}_{x})=\mathbb{Q}_{2}[i].$ 
 \item[]$R_{x}=\mathbb{Z}_{2}[T]/(T^2+1)=\mathbb{Z}_{2}[Y]/((Y-1)^2+1)=\mathbb{Z}_{2}[Y]/((Y^2-2Y+2)$ est un anneau de valuation discr\`ete totalement ramifi\'e de degr\'e deux sur $\mathbb{Z}_{2},$ car $Y^2-2Y+2$ est un polyn\^ome d'Eisenstein. 
 \item[]$\mathfrak{m}_{R_{x}}$ est engendr\'e par la classe de $T+1$.
 \item[]$ x_{0}=T+1\in\mathbb{F}_{2}[T]$.
 \item[]$ k(x_{0})=\mathbb{F}_{2}=\mathbb{F}_{2}[T]/(T+1)=\mathbb{Z}_{2}[T]/(T^2+1,T+1)=k(R_{x})$.
 \item[]$x= T^2+1\in  ]x_{0}[ \setminus Êr\acute{e}d^{-1}(x_{0}), y=T+1=Teich(x_{0})\in r\acute{e}d^{-1}(x_{0}), z=T+3\in r\acute{e}d^{-1}(x_{0})$.
  \item[]$K(R_{x})=K(\mathcal{V}_{x})$ est une extension totalement ramifi\'ee de degr\'e $2$ de $K=\mathbb{Q}_{2}$.
    \item[]$k(R_{x})=\mathbb{F}_{2}=k(\mathcal{V}_{x})$.
   \end{enumerate}
 
$\rondII$
 $R_{x}=\mathcal{V}_{x}=\mathbb{Z}_{3}[\sqrt{2},\sqrt{3}].$
 \begin{enumerate}
 \item[]Ici $K=\mathbb{Q}_{3},\ K(R_{x})=K(\mathcal{V}_{x})=\mathbb{Q}_{3}[\sqrt{2},\sqrt{3}]$. 
 \item[]Le polyn\^ome minimal de $\sqrt{2}+\sqrt{3}$ sur $\mathbb{Z}_{3}$ est $P=(T-\sqrt{2}-\sqrt{3})(T-\sqrt{2}+\sqrt{3})(T+\sqrt{2}+\sqrt{3})(T+\sqrt{2}-\sqrt{3})=T^4-10T^2+1\equiv (T^2+1)^2 \ mod\ 3$, car on v\'erifie que chacune des quatre racines de $P$ est primitive, i.e. engendre $\mathbb{Q}_{3}[\sqrt{2},\sqrt{3}]$ qui est de degr\'e $4$ sur $\mathbb{Q}_{3}$, comme extension compos\'ee des deux extensions de degr\'e $2$ totalement ramifi\'ees $\mathbb{Q}_{3}[\sqrt{3}]$ et $\mathbb{Q}_{3}[\sqrt{6}]$ de $\mathbb{Q}_{3}$: faisant par exemple le raisonnement avec $\sqrt{2}+\sqrt{3}$, on a $1/(\sqrt{3}+\sqrt{2})=\sqrt{3}-\sqrt{2} , 2\sqrt{3}=(\sqrt{3}+\sqrt{2})+(\sqrt{3}-\sqrt{2}),2\sqrt{2}=(\sqrt{3}+\sqrt{2})-(\sqrt{3}-\sqrt{2})$
 , d'o\`u $\sqrt{2}\in\mathbb{Q}_{3}[\sqrt{2}+\sqrt{3}]$ et $\sqrt{3}\in\mathbb{Q}_{3}[\sqrt{2}+\sqrt{3}]$; ainsi $P$ est bien irr\'eductible sur $\mathbb{Q}_{3}$ et $\mathbb{Q}_{3}[\sqrt{2}+\sqrt{3}]=\mathbb{Q}_{3}[\sqrt{2},\sqrt{3}]=\mathbb{Q}_{3}[T]/(T^4-10T^2+1)$. De plus $1/(\sqrt{3}+\sqrt{2})\in\mathbb{Z}_{3}[\sqrt{2}+\sqrt{3}]$ car le coefficient constant de $P$ est inversible dans $\mathbb{Z}_{3}$.
  \item[]$R_{x}=\mathbb{Z}_{3}[\sqrt{2}+\sqrt{3}]=\mathbb{Z}_{3}[T]/(T^4-10T^2+1)$.
  \item[]$\mathfrak{m}_{R_{x}}$ est engendr\'e par la classe de $T^2+1$.
  \item[]$ x_{0}=T^2+1\in\mathbb{F}_{3}[T]$.
  \item[]$ k(x_{0})=\mathbb{F}_{9}=\mathbb{F}_{3}[T]/(T^2+1)=\mathbb{Z}_{3}[T]/(T^4-10T^2+1,T^2+1)=k(R_{x}).$
 \item[]$x= T^4-10T^2+1\in  ]x_{0}[ \setminus Êr\acute{e}d^{-1}(x_{0}),y= T^2+1=Teich(x_{0})\in r\acute{e}d^{-1}(x_{0}), z=T^2+3T+1\in r\acute{e}d^{-1}(x_{0})$.
 \item[]$K(R_{x})=K(\mathcal{V}_{x})$ est une extension ramifi\'ee de degr\'e $4$ de $K=\mathbb{Q}_{3}$, d'indice de ramification $2$ et de degr\'e r\'esiduel $2$: remarquons que $K(R_{x})$ n'est pas totalement ramifi\'ee sur $\mathbb{Q}_{3}$, bien qu'elle soit extension compos\'ee de deux extensions totalement ramifi\'ees sur $\mathbb{Q}_{3}$.
   \item[]$k(R_{x})= \mathbb{F}_{9}=k(\mathcal{V}_{x})$.
 \end{enumerate}
 
 $\rondIII$
  $R_{x}=\mathcal{V}_{x}=\mathbb{Z}_{3}[i,\sqrt{3}].$
 \begin{enumerate}
 \item[]Ici $K=\mathbb{Q}_{3},\ K(R_{x})=K(\mathcal{V}_{x})=\mathbb{Q}_{3}[i,\sqrt{3}].$ 
 \item[] Le polyn\^ome minimal de $i+\sqrt{3}$ sur $\mathbb{Z}_{3}$ est $P=(T-i-\sqrt{3})(T-i+\sqrt{3})(T+i+\sqrt{3})(T+i-\sqrt{3})=T^4-4T^2+16\equiv (T^2+1)^2 \ mod\ 3$,  car on v\'erifie que chacune des quatre racines de $P$ est primitive, i.e. engendre $\mathbb{Q}_{3}[i,\sqrt{3}]$ qui est de degr\'e $4$ sur $\mathbb{Q}_{3}$, comme extension compos\'ee des deux extensions de degr\'e $2$ totalement ramifi\'ees $\mathbb{Q}_{3}[\sqrt{3}]$ et $\mathbb{Q}_{3}[\sqrt{-3}]$ de $\mathbb{Q}_{3}$: faisant par exemple le raisonnement avec $i+\sqrt{3}$, on a $4/(i+\sqrt{3})= -i+\sqrt{3}, 2\sqrt{3}=(\sqrt{3}+i)+(\sqrt{3}-i), 2i=(i+\sqrt{3})-(-i+\sqrt{3})$, d'o\`u $\sqrt{3}\in\mathbb{Q}_{3}[i+\sqrt{3}]$ et $i\in\mathbb{Q}_{3}[i+\sqrt{3}]$; ainsi $P$ est bien irr\'eductible sur $\mathbb{Q}_{3}$ et $\mathbb{Q}_{3}[i+\sqrt{3}]=\mathbb{Q}_{3}[i,\sqrt{3}]=\mathbb{Q}_{3}[T]/(T^4-4T^2+16)$. De plus $1/(i+\sqrt{3})\in\mathbb{Z}_{3}[i+\sqrt{3}]$ car le coefficient constant de $P$ est inversible dans $\mathbb{Z}_{3}$.
 \item[]$R_{x}=\mathbb{Z}_{3}[i+\sqrt{3}]=\mathbb{Z}_{3}[T]/(T^4-4T^2+16)$.
  \item[]$\mathfrak{m}_{R_{x}}$ est engendr\'e par la classe de $T^2+1$.
  \item[]$ x_{0}=T^2+1\in\mathbb{F}_{3}[T]$.
  \item[]$ k(x_{0})=\mathbb{F}_{9}=\mathbb{F}_{3}[T]/(T^2+1)=\mathbb{Z}_{3}[T]/(T^4-4T^2+16,T^2+1)=k(R_{x}).$
 \item[]$x= T^4-4T^2+16\in  ]x_{0}[ \setminus Êr\acute{e}d^{-1}(x_{0}), y=T^2+1=Teich(x_{0})\in r\acute{e}d^{-1}(x_{0}), z=T^2+7\in r\acute{e}d^{-1}(x_{0})$.
 \item[]$K(R_{x})=K(\mathcal{V}_{x})$ est une extension ramifi\'ee de degr\'e $4$ de $K=\mathbb{Q}_{3}$, d'indice de ramification $2$ et de degr\'e r\'esiduel $2$: remarquons que $K(R_{x})$ n'est pas totalement ramifi\'ee sur $\mathbb{Q}_{3}$, bien qu'elle soit extension compos\'ee de deux extensions totalement ramifi\'ees sur $\mathbb{Q}_{3}$.
   \item[]$k(R_{x})= \mathbb{F}_{9}=k(\mathcal{V}_{x})$.
 \end{enumerate}
 
\subsubsection*{3.4.3.  \textbf{3\ieme} cas: $R_{x}$ n'est pas un anneau de valuation discr\`ete}
\subsubsection*{3.4.3.1. $\mathcal{V}_{x}$ non ramifi\'e sur $\mathcal{V}$}

$R_{x}=\mathbb{Z}_{3}[3\sqrt{2}]\underset{\neq}{\hookrightarrow}\mathcal{V}_{x}=\mathbb{Z}_{3}[\sqrt{2}].$

\begin{enumerate}
 \item[]Ici $K=\mathbb{Q}_{3},\ K(R_{x})=K(\mathcal{V}_{x})=\mathbb{Q}_{3}[\sqrt{2}].$ 
 \item[]$R_{x}=\mathbb{Z}_{3}[T]/(T^2-18)$: $T^2-18$ est bien irr\'eductible sur $\mathbb{Z}_{3}$ sinon il aurait une racine dans $\mathbb{Q}_{3}$, \`a savoir $3\sqrt{2}$ ou $-3\sqrt{2}$ et donc $\sqrt{2}\in\mathbb{Q}_{3}$ ce qui n'est pas.
 \item[]$R_{x}$ n'est pas un anneau de valuation discr\`ete car $\sqrt{2}\notin R_{x}$.
  \item[]$\mathfrak{m}_{R_{x}}$ est engendr\'e par $3$ et $T$.
  \item[]$ x_{0}=T\in\mathbb{F}_{3}[T]$.
  \item[]$k(x_{0})=\mathbb{F}_{3}=\mathbb{F}_{3}[T]/(T)=\mathbb{Z}_{3}[T]/(T^2-18,3,T)=k(R_{x})$.
 \item[]$x= T^2-18\in  ]x_{0}[ \setminus Êr\acute{e}d^{-1}(x_{0}), y=T=Teich(x_{0})\in r\acute{e}d^{-1}(x_{0}), z=T+3\in r\acute{e}d^{-1}(x_{0})$.
  \item[]$K(R_{x})=K(\mathcal{V}_{x})$ est une extension non ramifi\'ee de degr\'e $2$ de $K=\mathbb{Q}_{3}$. 
  \item[]$k(R_{x})=\mathbb{F}_{3}\subset \mathbb{F}_{9}=k(\mathcal{V}_{x})$.
 \end{enumerate}

\subsubsection*{3.4.3.2. $\mathcal{V}_{x}$ totalement ramifi\'e sur $\mathcal{V}$}

$R_{x}=\mathbb{Z}_{2}[2i]\underset{\neq}{\hookrightarrow}\mathcal{V}_{x}=\mathbb{Z}_{2}[i]$.

\begin{enumerate}
 \item[]Ici $K=\mathbb{Q}_{2},\ K(R_{x})=K(\mathcal{V}_{x})=\mathbb{Q}_{2}[i].$ 
 \item[]$R_{x}=\mathbb{Z}_{2}[T]/(T^2+4)$: $T^2+4$ est bien irr\'eductible sur $\mathbb{Z}_{2}$ sinon il aurait une racine dans $\mathbb{Q}_{2}$, \`a savoir $2i$ ou $-2i$ et donc $i\in\mathbb{Q}_{2}$ ce qui n'est pas.
 \item[]$R_{x}$ n'est pas un anneau de valuation discr\`ete car $i\notin R_{x}$.
  \item[]$\mathfrak{m}_{R_{x}}$ est engendr\'e par $2$ et $T$.
  \item[]$ x_{0}=T\in\mathbb{F}_{2}[T]$.
  \item[]$k(x_{0})=\mathbb{F}_{2}=\mathbb{F}_{2}[T]/(T)=\mathbb{Z}_{2}[T]/(T^2+4,2,T)=k(R_{x})$.
 \item[]$x= T^2+4\in  ]x_{0}[ \setminus Êr\acute{e}d^{-1}(x_{0}), y=T=Teich(x_{0})\in r\acute{e}d^{-1}(x_{0}), z=T+2\in r\acute{e}d^{-1}(x_{0})$.
  \item[]$K(R_{x})=K(\mathcal{V}_{x})$ est une extension totalement ramifi\'ee de degr\'e $2$ de $K=\mathbb{Q}_{2}$. 
  \item[]$k(R_{x})=\mathbb{F}_{2}=k(\mathcal{V}_{x})$.
  \item[]Le morphisme $R_{x}=\mathbb{Z}_{2}[2i]\hookrightarrow\mathcal{V}_{x}=\mathbb{Z}_{2}[i]$ est ramifi\'e [SGA 1, I, \S7, cor 7.5].
 \end{enumerate}

\subsubsection*{3.4.3.3. $\mathcal{V}_{x}$ ramifi\'e non totalement sur $\mathcal{V}$}

$\rondI$ $R_{x}=\mathbb{Z}_{3}[\sqrt{2},3\sqrt{3}]\underset{\neq}{\hookrightarrow}\mathcal{V}_{x}=\mathbb{Z}_{3}[\sqrt{2},\sqrt{3}]$.\\
\begin{enumerate}
 \item[]Ici $K=\mathbb{Q}_{3},\ K(R_{x})=K(\mathcal{V}_{x})=\mathbb{Q}_{3}[\sqrt{2},\sqrt{3}].$ 
 \item[]Le polyn\^ome minimal de $\sqrt{2}+3\sqrt{3}$ sur $\mathbb{Z}_{3}$ est $P=(T-\sqrt{2}-3\sqrt{3})(T-\sqrt{2}+3\sqrt{3})(T+\sqrt{2}+3\sqrt{3})(T+\sqrt{2}-3\sqrt{3})=T^4-58T^2+625\equiv (T^2+1)^2 \ mod\ 3$, car on v\'erifie que chacune des quatre racines de $P$ est primitive, i.e. engendre $\mathbb{Q}_{3}[\sqrt{2},\sqrt{3}]$ qui est de degr\'e $4$ sur $\mathbb{Q}_{3}$: faisant par exemple le raisonnement avec $\sqrt{2}+3\sqrt{3}$, on a $1/(3\sqrt{3}+\sqrt{2})=(3\sqrt{3}-\sqrt{2})/25$, d'o\`u $\sqrt{2}\in\mathbb{Q}_{3}[\sqrt{2}+3\sqrt{3}]$ et $\sqrt{3}\in\mathbb{Q}_{3}[\sqrt{2}+3\sqrt{3}]$; ainsi $P$ est bien irr\'eductible sur $\mathbb{Q}_{3}$ et $\mathbb{Q}_{3}[\sqrt{2}+3\sqrt{3}]=\mathbb{Q}_{3}[\sqrt{2},\sqrt{3}]=\mathbb{Q}_{3}[T]/(T^4-58T^2+625)$. De plus $1/(\sqrt{2}+3\sqrt{3})\in\mathbb{Z}_{3}[\sqrt{2}+\sqrt{3}]$ car le coefficient constant de $P$ est inversible dans $\mathbb{Z}_{3}$; et donc $\sqrt{2}\in\mathbb{Z}_{3}[\sqrt{2}+3\sqrt{3}]$ et $\sqrt{3}\in\mathbb{Z}_{3}[\sqrt{2}+3\sqrt{3}]$.
  \item[]$R_{x}=\mathbb{Z}_{3}[\sqrt{2}+3\sqrt{3}]=\mathbb{Z}_{3}[T]/(T^4-58T^2+625)$.
  \item[]$\mathfrak{m}_{R_{x}}$ est engendr\'e par $3$ et $T^2+1$.
  \item[]$ x_{0}=T^2+1\in\mathbb{F}_{3}[T]$.
  \item[]$ k(x_{0})=\mathbb{F}_{9}=\mathbb{F}_{3}[T]/(T^2+1)=\mathbb{Z}_{3}[T]/(T^4-58T^2+625,3,T^2+1)=k(R_{x}).$
 \item[]$x= T^4-58T^2+625\in  ]x_{0}[ \setminus Êr\acute{e}d^{-1}(x_{0}), y=T^2+1=Teich(x_{0})\in r\acute{e}d^{-1}(x_{0}), z=T^2+3T+4\in r\acute{e}d^{-1}(x_{0})$.
  \item[]$K(R_{x})=K(\mathcal{V}_{x})$ est une extension ramifi\'ee de degr\'e $4$ de $K=\mathbb{Q}_{3}$, d'indice de ramification $2$ et de degr\'e r\'esiduel $2$.
    \item[]$k(R_{x})= \mathbb{F}_{9}=k(\mathcal{V}_{x})$.
    \item[]Le morphisme $R_{x}=\mathbb{Z}_{3}[\sqrt{2}+3\sqrt{3}]\hookrightarrow\mathcal{V}_{x}=\mathbb{Z}_{3}[\sqrt{2},\sqrt{3}]$ est ramifi\'e [SGA 1, I, \S7, cor 7.5].
 \end{enumerate}

$\rondII$$R_{x}=\mathbb{Z}_{3}[3\sqrt{2}+\sqrt{3}]\underset{\neq}{\hookrightarrow}\mathbb{Z}_{3}[3\sqrt{2},\sqrt{3}]\underset{\neq}{\hookrightarrow}\mathcal{V}_{x}=\mathbb{Z}_{3}[\sqrt{2},\sqrt{3}]$.\\

\begin{enumerate}
 \item[]Ici $K=\mathbb{Q}_{3},\ K(R_{x})=K(\mathcal{V}_{x})=\mathbb{Q}_{3}[\sqrt{2},\sqrt{3}].$ 
 \item[]Le polyn\^ome minimal de $3\sqrt{2}+\sqrt{3}$ sur $\mathbb{Z}_{3}$ est $P=(T-3\sqrt{2}-\sqrt{3})(T-3\sqrt{2}+\sqrt{3})(T+3\sqrt{2}+\sqrt{3})(T+3\sqrt{2}-\sqrt{3})=T^4-42T^2+225\equiv T^4 \ mod\ 3$,
  car on v\'erifie que chacune des racines de $P$ engendre l'extension $\mathbb{Q}_{3}[\sqrt{2},\sqrt{3}]$ de degr\'e $4$ sur $\mathbb{Q}_{3}$; ainsi $P$ est bien irr\'eductible sur $\mathbb{Q}_{3}$ et $\mathbb{Q}_{3}[3\sqrt{2}+\sqrt{3}]=\mathbb{Q}_{3}[\sqrt{2},\sqrt{3}]=\mathbb{Q}_{3}[T]/(T^4-42T^2+225)$.
  \item[]$R_{x}=\mathbb{Z}_{3}[3\sqrt{2}+\sqrt{3}]=\mathbb{Z}_{3}[T]/(T^4-42T^2+225)$.
  \item[]$\mathfrak{m}_{R_{x}}$ est engendr\'e par $3$ et $T$.
  \item[]$ x_{0}=T\in\mathbb{F}_{3}[T]$.
  \item[]$ k(x_{0})=\mathbb{F}_{3}=\mathbb{F}_{3}[T]/(T)=\mathbb{Z}_{3}[T]/(T^4-42T^2+225,3,T)=k(R_{x}).$
 \item[]$x= T^4-42T^2+225\in  ]x_{0}[ \setminus Êr\acute{e}d^{-1}(x_{0}), y=T=Teich(x_{0})\in r\acute{e}d^{-1}(x_{0}), z=T+3\in r\acute{e}d^{-1}(x_{0})$.
  \item[]$K(R_{x})=K(\mathcal{V}_{x})$ est une extension ramifi\'ee de degr\'e $4$ de $K=\mathbb{Q}_{3}$, d'indice de ramification $2$ et de degr\'e r\'esiduel $2$.
    \item[]$k(R_{x})= \mathbb{F}_{3}\subset\mathbb{F}_{9}=k(\mathcal{V}_{x})$.
 \end{enumerate}

$\rondIII$$R_{x}=\mathbb{Z}_{3}[i,3\sqrt{3}]\underset{\neq}{\hookrightarrow}\mathcal{V}_{x}=\mathbb{Z}_{3}[i,\sqrt{3}]$.\\

\begin{enumerate}
 \item[]Ici $K=\mathbb{Q}_{3},\ K(R_{x})=K(\mathcal{V}_{x})=\mathbb{Q}_{3}[i,\sqrt{3}].$ 
 \item[]Le polyn\^ome minimal de $i+3\sqrt{3}$ sur $\mathbb{Z}_{3}$ est $P=(T-i-3\sqrt{3})(T-i+3\sqrt{3})(T+i+3\sqrt{3})(T+i-3\sqrt{3})=T^4-52T^2+784\equiv (T^2+1)^2 \ mod\ 3$,  car on v\'erifie que chacune des racines de $P$ engendre l'extension $\mathbb{Q}_{3}[i,\sqrt{3}]$ de degr\'e $4$ sur $\mathbb{Q}_{3}$; ainsi $P$ est bien irr\'eductible sur $\mathbb{Q}_{3}$ et $\mathbb{Q}_{3}[i+3\sqrt{3}]=\mathbb{Q}_{3}[i,\sqrt{3}]=\mathbb{Q}_{3}[T]/(T^4-52T^2+784)$. De plus $1/(i+3\sqrt{3})=(3\sqrt{3}-i)/29\in\mathbb{Z}_{3}[(i+3\sqrt{3}]$ car le coefficient constant de $P$ est inversible dans $\mathbb{Z}_{3}$; et donc $i\in\mathbb{Z}_{3}[i+3\sqrt{3}]$ et $3\sqrt{3}\in\mathbb{Z}_{3}[i+3\sqrt{3}]$. 
  \item[]$R_{x}=\mathbb{Z}_{3}[i+3\sqrt{3}]=\mathbb{Z}_{3}[T]/(T^4-52T^2+784)$.
  \item[]$\mathfrak{m}_{R_{x}}$ est engendr\'e par $3$ et $T^2+1$.
  \item[]$ x_{0}=T^2+1\in\mathbb{F}_{3}[T]$.
  \item[]$ k(x_{0})=\mathbb{F}_{9}=\mathbb{F}_{3}[T]/(T^2+1)=\mathbb{Z}_{3}[T]/(T^4-52T^2+784,3,T^2+1)=k(R_{x}).$
 \item[]$x= T^4-52T^2+784\in  ]x_{0}[ \setminus Êr\acute{e}d^{-1}(x_{0}), y=T^2+1=Teich(x_{0})\in r\acute{e}d^{-1}(x_{0}), z=T^2+4\in r\acute{e}d^{-1}(x_{0})$.
  \item[]$K(R_{x})=K(\mathcal{V}_{x})$ est une extension ramifi\'ee de degr\'e $4$ de $K=\mathbb{Q}_{3}$, d'indice de ramification $2$ et de degr\'e r\'esiduel $2$.
    \item[]$k(R_{x})=\mathbb{F}_{9}=k(\mathcal{V}_{x})$.
     \item[]Le morphisme $R_{x}=\mathbb{Z}_{3}[i+3\sqrt{3}]\hookrightarrow\mathcal{V}_{x}=\mathbb{Z}_{3}[i,\sqrt{3}]$ est ramifi\'e [SGA 1, I, \S7, cor 7.5].
 \end{enumerate}

$\rondIV$$R_{x}=\mathbb{Z}_{3}[3i+\sqrt{3}]\underset{\neq}{\hookrightarrow}\mathbb{Z}_{3}[3i,\sqrt{3}]\underset{\neq}{\hookrightarrow}\mathcal{V}_{x}=\mathbb{Z}_{3}[i,\sqrt{3}]$.\\

\begin{enumerate}
 \item[]Ici $K=\mathbb{Q}_{3},\ K(R_{x})=K(\mathcal{V}_{x})=\mathbb{Q}_{3}[i,\sqrt{3}].$ 
 \item[]Le polyn\^ome minimal de $3i+\sqrt{3}$ sur $\mathbb{Z}_{3}$ est $P=(T-3i-\sqrt{3})(T-3i+\sqrt{3})(T+3i+\sqrt{3})(T+3i-\sqrt{3})= T^4+12T^2+144\equiv T^4 \ mod\ 3$,  car on v\'erifie que chacune des racines de $P$ engendre l'extension $\mathbb{Q}_{3}[i,\sqrt{3}]$ de degr\'e $4$ sur $\mathbb{Q}_{3}$; ainsi $P$ est bien irr\'eductible sur $\mathbb{Q}_{3}$ et $\mathbb{Q}_{3}[3i+\sqrt{3}]=\mathbb{Q}_{3}[i,\sqrt{3}]=\mathbb{Q}_{3}[T]/(T^4+12T^2+144)$.
  \item[]$R_{x}=\mathbb{Z}_{3}[3i+\sqrt{3}]=\mathbb{Z}_{3}[T]/(T^4+12T^2+144)$.
  \item[]$\mathfrak{m}_{R_{x}}$ est engendr\'e par $3$ et $T$.
  \item[]$ x_{0}=T\in\mathbb{F}_{3}[T]$.
  \item[]$ k(x_{0})=\mathbb{F}_{3}=\mathbb{F}_{3}[T]/(T)=\mathbb{Z}_{3}[T]/(T^4+12T^2+144,3,T)=k(R_{x}).$
 \item[]$x= T^4+12T^2+144\in  ]x_{0}[ \setminus Êr\acute{e}d^{-1}(x_{0}), y=T=Teich(x_{0})\in r\acute{e}d^{-1}(x_{0}), z=T+3\in r\acute{e}d^{-1}(x_{0})$.
  \item[]$K(R_{x})=K(\mathcal{V}_{x})$ est une extension ramifi\'ee de degr\'e $4$ de $K=\mathbb{Q}_{3}$, d'indice de ramification $2$ et de degr\'e r\'esiduel $2$.
    \item[]$k(R_{x})=\mathbb{F}_{3}\subset\mathbb{F}_{9}=k(\mathcal{V}_{x})$.
 \end{enumerate}

\newpage
\section*{4. Images directes de $F$-isocristaux convergents}
\subsection*{4.0. }On reprend les hypoth\`eses de 3.0.
\subsection*{4.1. Convergence des images directes}
\subsubsection*{4.1.1.}  On  suppose de plus  en 4.1 que $K$ contient les racines $q^{i\grave{e}mes}$ de l'unit\'e. Soient $\mathcal{S},\mathcal{S}'$ deux $\mathcal{V}$-sch\'emas formels lisses; on consid\`ere un diagramme commutatif

 $$
\begin{array}{cccccc}
\xymatrix{
 &X\ar@{.>}[dd]^(.3){f} |\hole \ \ar@{^{(}->}[rr]^{i_{X}} & & \mathcal{X} \ar[dd]^{h}&& \\
 X' \ \ar@{^{(}->}[rr]^(.8){i_{X'}} \ar[dd]_(.5){f'} \ar[ur]^{\varphi_{X}} & & \mathcal{X}' \ar[dd]_(.7){h'} \ar[ur]^{\varphi_{\mathcal{X}}}&&& \\
 &S\ \ar@{^{(}.>}[rr]^(.8){i_{S}}  |\hole &&\mathcal{S}\ar[rr]^{\rho}&&Spf\ \mathcal{V}\\
 S' \ \ar@{^{(}->}[rr]_{i_{S'}} \ar@{.>}[ur]_{\varphi_{S}} && \mathcal{S'} \ar[ur]_{\varphi_{\mathcal{S}}}\ar[rr]_{\rho'}&\ &\ \ Spf\  \mathcal{V}\ar[ur]_{\theta} &
 }
\end{array}
\leqno{(4.1.1.1)}
  $$
 dans lequel $f,f',\varphi_{S},\varphi_{X}$ sont des morphismes de $k$-sch\'emas lisses, $h,h',\rho,\rho', \varphi_{\mathcal{S}}$,\\
 $\varphi_{\mathcal{X}}$ sont des morphismes de $\mathcal{V}$-sch\'emas formels s\'epar\'es lisses, $h$, $h'$ sont propres et lisses, les $i$ sont des immersions ferm\'ees et les carr\'es verticaux sont cart\'esiens.

\vskip 3mm
\noindent \textbf{Th\'eor\`eme (4.1.2)}. \textit{Supposons le corps $k$ parfait. Soient $S$ un $k$-sch\'ema lisse, $\mathcal{S}$ un $\mathcal{V}$-sch\'ema formel lisse relevant $S$ et $f : X \rightarrow S$ un $k$-morphisme propre et lisse relevable en un morphisme propre et lisse $h : \mathcal{X} \rightarrow \mathcal{S}$ de $\mathcal{V}$-sch\'emas formels. Alors}\\
\begin{enumerate}
\item[(4.1.2.1)] \textit{Pour tout entier $i \geqslant 0$, $f$ induit un foncteur} 
$$
R^{i} f_{conv^{\ast}} : F^{a}\mbox{-}Isoc(X/K) \longrightarrow F^{a}\mbox{-}Isoc(S/K)\ .
$$

\item[ (4.1.2.2)] \textit{Le foncteur pr\'ec\'edent est compatible au passage \`a la fibre en un point ferm\'e $s_{0}=Spec\ k(s_{0})$ de $S$, c'est-\`a-dire : pour tout carr\'e cart\'esien}
$$
\xymatrix{
X_{s_{0}} \ar[r]^{i_{X_{s_{0}}}} \ar[d] _{f_{s_{0}}} & X \ar[d]^{f} \\
s_{0} \ar[r]_{i_{s_{0}}} & S
}
$$
\textit{et tout $\mathcal{E} \in F^{a}\mbox{-}Isoc(X/K)$ on a un isomorphisme de changement de base
$$
i^{\ast}_{s_{0}} R^{i} f_{conv^{\ast}}(\mathcal{E)}\ \displaystyle \mathop{\longrightarrow}^{\sim}   R^{i} f_{s_{0},conv^{\ast}}(i_{X_{s_{0}}}^{\ast}(\mathcal{E}))=H^{i}_{rig}(X_{s_{0}}/K(s_{0}), \mathcal{E}_{X_{s_{0}}})
$$
compatible aux Frobenius, dans lequel le corps $K(s_{0})$ a \'et\'e introduit au d\'ebut de 3.1 et $\mathcal{E}_{X_{s_{0}}}:=i_{X_{s_{0}}}^{\ast}(\mathcal{E})$.}
\item[ (4.1.2.3)] \textit{Sous les hypoth\`eses 4.1.1 le foncteur $R^{i} f_{conv^{\ast}}$ de 4.1.2.1 commute au changement de base $\varphi_{S}:S' \longrightarrow S$, i.e. pour $\mathcal{E}\in F^{a}\mbox{-}Isoc(X/K), \mathcal{E}'=\varphi_{X}^{\ast}(\mathcal{E})\in F^{a}\mbox{-}Isoc(X'/K)$, on a un diagramme commutatif d'isomorphismes}
$$
\xymatrix{
F^{\ast}_{S'}R^{i}f'_{rig^{\ast}}(X'/\mathcal{S}'; \mathcal{E}')\ar[r]^{\overset{\phi'^{i}}{\sim}}&R^{i}f'_{rig^{\ast}}(X'/\mathcal{S}'; \mathcal{E}')\\
F^{\ast}_{S}R^{i}f_{rig^{\ast}}(X/\mathcal{S}; \mathcal{E})\ar[r]^{\overset{\phi^{i}}{\sim}}\ar[u]_{\varphi^{\ast}_{S}}&R^{i}f_{rig^{\ast}}(X/\mathcal{S}; \mathcal{E})\ar[u]_{\varphi^{\ast}_{S}}
}
$$
\textit{dans lequel les fl\`eches verticales sont les isomorphismes de changement de base et les fl\`eches horizontales les isomorphismes de Frobenius.}
\end{enumerate}
\vskip 3mm
\noindent\textit{D\'emonstration}.\\
\textit{Pour (4.1.2.1).} Puisque l'on dispose d'un carr\'e cart\'esien

$$
\xymatrix{
X \ar@{^{(}->}[r]^{i_{X}} \ar[d] _{f} & \mathcal{X} \ar[d]^{h} \\
S \ar@{^{(}->}[r]_{i_{S}} & \mathcal{S}
}
$$
satisfaisant aux conditions de [Et 6,(3.4.4.3)], on en d\'eduit que pour $\mathcal{E}\in F^{a}\mbox{-}Isoc(X/K)$, on a [loc. cit.]
$$
\mathcal{E}_{i}:=R^{i}f_{rig^{\ast}}(X/\mathcal{S}; \mathcal{E})\in Isoc(S/K)\ .
$$
En notant $F_{S}:S\rightarrow S$ (resp $F_{X}:X\rightarrow X$) le Frobenius de $S$ (resp de $X$) (\'el\'evation \`a puissance $q$) et 
$$
\phi_{\mathcal{E}}: F^{\ast}_{X}(\mathcal{E})\overset {\sim}{\longrightarrow}\mathcal{E}
$$
l'isomorphisme de Frobenius de $\mathcal{E}$, il reste \`a construire un isomorphisme de Frobenius
$$
\phi_{i}=\phi_{\mathcal{E}_{i}}: F^{\ast}_{S}(\mathcal{E}_{i})\overset {\sim}{\longrightarrow}\mathcal{E}_{i}\ .
$$

Quitte \`a d\'ecomposer $S$ en somme de ses composantes connexes il suffit de d\'efinir $\phi_{\mathcal{E}_{i}}$ sur chacune de ces composantes connexes. Soit $S_{\alpha}$ un ouvert affine d'une composante connexe $S_{0}$ de $S$: comme le foncteur naturel
$$
F^{a}\mbox{-}Isoc(S_{0}/K) \longrightarrow F^{a}\mbox{-}Isoc(S_{\alpha}/K)
$$
est pleinement fid\`ele [Et 3, th\'eo 4] il suffit de d\'efinir $\phi_{\mathcal{E}_{i}}$ sur $S_{\alpha}$.\\

On part donc d'une d\'ecomposition [cor (2.1.3)]
 $\mathcal{S}= \displaystyle \mathop{\cup}_{\alpha}\ \mathcal{S}_{\alpha}$ , $\mathcal{S}_{\alpha} = \mbox{Spf}\ \hat{A}_{\alpha}$ o\`u les $A_{\alpha}$ sont des $\mathcal{V}$-alg\`ebres lisses et int\`egres,
 $S= \displaystyle \mathop{\cup}_{\alpha}\ S_{\alpha}$, $S_{\alpha}=Spec\ A_{\alpha,0}$, o\`u  $A_{\alpha,0} := A_{\alpha} / \pi A_{\alpha}$
 et on note $j_{S_{\alpha}}:S_{\alpha}\hookrightarrow S $ l'immersion ouverte. Choisissons comme en (3.1) un rel\`evement (fini et plat) $F_{\mathcal{S}_{\alpha}}: \mathcal{S}_{\alpha}\rightarrow \mathcal{S}_{\alpha}$ du Frobenius $F_{S_{\alpha}}: S_{\alpha}\rightarrow S_{\alpha}$, au-dessus de $\sigma: Spf\ \mathcal{V}\rightarrow Spf\ \mathcal{V}$. Puisque le diagramme commutatif \`a carr\'e cart\'esien
 
  $$
\begin{array}{c}
\xymatrix{
 &X\ar@{.>}[dd]^(.3){f} |\hole \ \ar@{^{(}->}[rr]^{i_{X}} & & \mathcal{X} \ar[dd]^{h} \\
 X_{\alpha} \ \ar@{^{(}->}[rr]^(.8){i_{X_{\alpha}}} \ar[dd]_(.5){f_{\alpha}} \ar@{^{(}->}[ur]^{j_{X_{\alpha}}} & & \mathcal{X}_{\alpha} \ar[dd]_(.7){h_{\alpha}} \ar@{^{(}->}[ur]^{j_{\mathcal{X}_{\alpha}}} \\
 &S\ \ar@{^{(}.>}[rr]^(.8){i_{S}}  |\hole &&\mathcal{S}&&\\
 S_{\alpha} \ \ar@{^{(}->}[rr]_{i_{S_{\alpha}}} \ar@{^{(}.>}[ur]_{j_{S_{\alpha}}} & & \mathcal{S}_{\alpha} \ar@{^{(}->}[ur]_{j_{\mathcal{S}_{\alpha}}}&&  
 }
\end{array}
  $$
fournit un isomorphisme de changement de base [Et 6, th\'eo (3.4.4)]
$$
j^{\ast}_{S_{\alpha}}R^{i}f_{rig^{\ast}}(X/\mathcal{S}; \mathcal{E})\overset{\sim}{\longrightarrow}R^{i}f_{\alpha, rig^{\ast}}(X_{\alpha}/\mathcal{S}_{\alpha}; j^{\ast}_{X_{\alpha}}(\mathcal{E}))\ ,
$$
on est ramen\'e \`a construire un isomorphisme de Frobenius sur
$$
\mathcal{E}^{i}_{\alpha}:=R^{i}f_{\alpha, rig^{\ast}}(X_{\alpha}/\mathcal{S}_{\alpha}; j^{\ast}_{X_{\alpha}}(\mathcal{E}))\in Isoc\ (S_{\alpha}/K)
$$
o\`u $\mathcal{E}_{\alpha}:= j^{\ast}_{X_{\alpha}}(\mathcal{E})$ est muni de l'isomorphisme de Frobenius
$$
\phi_{\mathcal{E}_{\alpha}}=j^{\ast}_{X_{\alpha}}(\phi_{\mathcal{E}}):F_{X_{\alpha}}^{\ast}(\mathcal{E}_{\alpha})\overset{\sim}{\longrightarrow}\mathcal{E}_{\alpha}\ .
$$
Notons
$$
\xymatrix{
X_{\alpha}^{(q)} \ar@{^{(}->}[r]^{i_{X^{(q)}_{\alpha}}} \ar[d] _{f^{(q)}_{\alpha}} & \mathcal{X}^{(q)}_{\alpha} \ar[d]^{h^{(q)}_{\alpha}} \\
S_{\alpha} \ar@{^{(}->}[r]_{i_{S_{\alpha}}} & \mathcal{S}_{\alpha}
}
$$
l'image inverse du carr\'e cart\'esien
$$
\xymatrix{
X_{\alpha} \ar@{^{(}->}[r]^{i_{X_{\alpha}}} \ar[d] _{f_{\alpha}} & \mathcal{X}_{\alpha} \ar[d]^{h_{\alpha}} \\
S_{\alpha} \ar@{^{(}->}[r]_{i_{S_{\alpha}}} & \mathcal{S}_{\alpha}
}
$$
par le carr\'e cart\'esien
$$
\xymatrix{
S_{\alpha} \ar@{^{(}->}[r]^{i_{S_{\alpha}}} \ar[d] _{F_{S_{\alpha}}} & \mathcal{X}_{\alpha} \ar[d]^{F_{\mathcal{S}_{\alpha}}} \\
S_{\alpha} \ar@{^{(}->}[r]_{i_{S_{\alpha}}} & \mathcal{S}_{\alpha}\ ,
}
$$
et consid\'erons la factorisation usuelle du Frobenius $F_{X_{\alpha}}$ (avec $F_{X_{\alpha}}^{\ast}(x)=x^{q}$)

$$
\xymatrix{X_{\alpha} \ar@/^1pc/[rrd]^{F_{X_{\alpha}}} \ar@/_/[rdd]_{f_{\alpha}}  \ar@{.>}[rd]^{F_{X_{\alpha}/S_{\alpha}}} \\
&  X_{\alpha}^{(q)} \ar[d]^{f_{\alpha}^{(q)}} \ar[r]^{\pi_{X_{\alpha}/S_{\alpha}}} & X_{\alpha} \ar[d]^{f_{\alpha}}\\
& S_{\alpha} \ar[r] ^{F_{S_{\alpha}}} & S_{\alpha}
}
\leqno{(4.1.2.4)}
$$
\noindent o\`u le carr\'e est cart\'esien. \\
\noindent On dispose d'un isomorphisme de changement de base [Et 6, (3.4.4.1)(ii)]
$$
F_{S_{\alpha}}^{\ast}R^{i}f_{\alpha, rig^{\ast}}(X_{\alpha}/\mathcal{S}_{\alpha}; \mathcal{E}_{\alpha})\overset{\sim}{\longrightarrow}R^{i}f^{(q)}_{\alpha, rig^{\ast}}(X^{(q)}_{\alpha}/\mathcal{S}_{\alpha}; \pi^{\ast}_{X_{\alpha}/S_{\alpha}}(\mathcal{E}_{\alpha})),
\leqno{(4.1.2.5)}
$$

\noindent d'un morphisme fonctoriel [C-T, 10.5.2]
$$
\eta^{i} :\ R^{i}f^{(q)}_{\alpha, rig^{\ast}}(X^{(q)}_{\alpha}/\mathcal{S}_{\alpha}; \pi^{\ast}_{X_{\alpha}/S_{\alpha}}(\mathcal{E}_{\alpha}))\longrightarrow R^{i}f_{\alpha, rig^{\ast}}(X_{\alpha}/\mathcal{S}_{\alpha};F_{X_{\alpha}}^{\ast} (\mathcal{E}_{\alpha}))
\leqno{(4.1.2.6)}
$$
induit par l'identit\'e de $S_{\alpha}$ et $F_{X_{\alpha}/S_{\alpha}}$,\\

\noindent et de l'isomorphisme \\
$$
R^{i}f_{\alpha, rig^{\ast}}(X_{\alpha}/\mathcal{S}_{\alpha};\phi_{\mathcal{E}_{\alpha}}):R^{i}f_{\alpha, rig^{\ast}}(X_{\alpha}/\mathcal{S}_{\alpha};F_{X_{\alpha}}^{\ast} (\mathcal{E}_{\alpha}))\overset{\sim}{\longrightarrow}R^{i}f_{\alpha, rig^{\ast}}(X_{\alpha}/\mathcal{S}_{\alpha}; \mathcal{E}_{\alpha})
\leqno{(4.1.2.7)}
$$
induit par 
$$
\phi_{\mathcal{E}_{\alpha}}=j^{\ast}_{X_{\alpha}}(\phi_{\mathcal{E}}):F_{X_{\alpha}}^{\ast}(\mathcal{E}_{\alpha})\overset{\sim}{\longrightarrow}\mathcal{E}_{\alpha}\ .
$$
\noindent Par composition de ces trois morphismes on obtient le morphisme de Frobenius de $\mathcal{E}^{i}_{\alpha}:=R^{i}f_{\alpha, rig^{\ast}}(X_{\alpha}/\mathcal{S}_{\alpha}; \mathcal{E}_{\alpha})$
$$
\phi^{i}:F_{S_{\alpha}}^{\ast}R^{i}f_{\alpha, rig^{\ast}}(X_{\alpha}/\mathcal{S}_{\alpha}; \mathcal{E}_{\alpha})\longrightarrow R^{i}f_{\alpha, rig^{\ast}}(X_{\alpha}/\mathcal{S}_{\alpha}; \mathcal{E}_{\alpha})\ ,
\leqno{(4.1.2.8)}
$$
et il s'agit de prouver que $\phi^{i}$ est un isomorphisme: pour \c ca il suffit de prouver que c'est le cas pour $\eta^{i}$. On sait d\'ej\`a que $\eta^{i}$ est un morphisme d'isocristaux convergents: comme la source et le but de $\eta^{i}$ sont des $\mathcal{O}_{\mathcal{S}_{\alpha K}}$-modules coh\'erents [Et 6, (3.4.4.3)] il suffit, pour montrer que $\eta^{i}$ est un isomorphisme, de montrer que c'est un isomorphisme fibre \`a fibre aux points ferm\'es de $\mathcal{S}_{\alpha K}= Spm\ (\hat{A}_{\alpha K})$ [B-G-R, 9.4.2 cor 7].\\

Soit $s\in Spm\ (\hat{A}_{\alpha K})$ un point ferm\'e de $\mathcal{S}_{\alpha K}$, correspondant \`a un id\'eal maximal $\mathfrak{m}_{s}$ de $\hat{A}_{ K}:=\hat{A}_{\alpha K}$; on utilise alors les notations de (3.2), en particulier on pose $s_{0}=\overline{sp}(s)$:\\
$$
\xymatrix@C=1,2cm  
{
A_{0} \ar@{->>}[r]^{\widetilde{\psi}_{s}\ \ \ \  \ \ \ \ \ \ \ \ \ \ \ \ \ \ \ \ \ \ \   }\ar@{=}[d]&k(R_{s})=R_{s}/\mathfrak{m}_{R_{s}}=A_{0}/\mathfrak{m}_{\overline{sp}(s)}=\hat{A}/\mathfrak{m}_{sp(s)}\\
 A_{0}\ar@{->>}[r]^{\overline{\psi}_{s} \  \  \ \ \ \ \ \  \ \ \ \ \ \ \  \ \  \  \  \  }&R_{s}/\pi R_{s}=A_{0}/\widetilde{\mathfrak{m}}_{s}=\hat{A}/<\pi, \overset{\circ}{\mathfrak{m}}_{s}>\ar@{->>}[u]\\
 \hat{A} \ar@{->>}[r]^{{\psi}_{s}}\ar@{->>}[u]^{\rho}&R_{s}=\hat{A}/\overset{\circ}{\mathfrak{m}}_{s}\ar@{->>}[u] \\
\hat{A}_{K}\ar@{->>}[r]_{\theta_{s}}\ar@{<-^{)}}[u]^{j}&K_{s}=\hat{A}_{K}/\mathfrak{m}_{s}\ar@{<-^{)}}[u]^{j_{s}}\ .
}
\leqno{(4.1.2.9)}
$$
Notons
$$
\mathcal{Y}_{s}=Spf\ R_{s},\  \mathcal{Y}_{K_{s}}=Spm\ K_{s},\  Y_{s_{0}}=Spec\ k(R_{s})\ ,
$$
$$
\begin{array}{llcccrr}
\theta_{s}^{\#}:\mathcal{Y}_{K_{s}}=Spm\ K_{s}&\longrightarrow & Spm\ \hat{A}_{K}=\mathcal{S}_{\alpha K}\  , \\
\psi_{s}^{\#}:\mathcal{Y}_{s}=Spf\ R_{s}&\longrightarrow & Spf\ \hat{A}=\mathcal{S}_{\alpha}\  ,
\end{array}
$$
$$
i_{s_{0}}:=\widetilde{\psi}_{s}^{\#}:Y_{s_{0}}=Spec\ k(R_{s})=Spec\ k(s_{0}) \rightarrow  Spec\ A_{0}=S_{\alpha}\ .
$$
Soit $\mathcal{F}\in Isoc(S_{\alpha}/K)$ un isocristal convergent dont une r\'ealisation [B3, (2.3.2)] est un 
$\mathcal{O}_{\mathcal{S}_{\alpha K}}$-module coh\'erent $\mathcal{F}_{K}$ \`a connexion int\'egrable convergente: une r\'ealisation de $i_{s_{0}}^{\ast}(\mathcal{F})$, la fibre de $\mathcal{F}$ en $s_{0}$, est donn\'ee par $(\theta_{s}^{\#})^{\ast}(\mathcal{F}_{K})$ [B3, (2.3.2)(iv)] en prenant l'image inverse de $\mathcal{F}_{K}$ par le diagramme
$$
\xymatrix{
S_{\alpha}\ar@{^{(}->}[r]&\mathcal{S}_{\alpha}\\
Y_{s_{0}}\ar@{^{(}->}[r]\ar[u]^{\widetilde{\psi}_{s}^{\#}}&\mathcal{Y}_{s}\ar[u]_{\psi_{s}^{\#}}&,
}
\leqno{(4.1.2.10)}
$$
 et cette image inverse ne d\'epend, \`a isomorphisme canonique pr\`es, que de $\widetilde{\psi}_{s}^{\#}$, qui lui-m\^eme ne d\'epend que de la sp\'ecialisation $s_{0}=\overline{sp}(s)$ de $s$. En fait, la connexion dont est muni $\mathcal{F}_{K}$ fournit un syst\`eme compatible d'isomorphismes canoniques 
$$
s_{1}^{\ast}(\mathcal{F}_{K})\overset{\sim}{\rightarrow}s_{2}^{\ast}(\mathcal{F}_{K})
$$
pour tous points $s_{1},s_{2}\in Spm\ \hat{A}_{K}$ ayant m\^eme sp\'ecialisation $s_{0}$, i.e. pour des points $s_{i}\in]s_{0}[=\overline{sp}^{-1}(s_{0}), i=1,2$ [B3; (2.2.17),(2.3.2)(iv)]. L'image inverse $(\theta_{s}^{\#})^{\ast}(\mathcal{F}_{K})$ est canoniquement isomorphe \`a celle obtenue en rempla\c cant $\mathcal{Y}_{s}$ par un $\mathcal{V}$-sch\'ema formel fini \'etale $\mathcal{Y}_{s'}$ et  $\psi_{s}^{\#}$ par $\psi_{s'}^{\#}$ tels que le diagramme suivant 

$$
\xymatrix{
S_{\alpha}\ar@{^{(}->}[r]&\mathcal{S}_{\alpha}\\
Y_{s_{0}}\ar@{^{(}->}[r]\ar[u]^{i_{s_{0}}}&\mathcal{Y}_{s'}\ar[u]_{\psi_{s'}^{\#}}
}
\leqno{(4.1.2.11)}
$$
soit cart\'esien. En fait on choisit $s'=Teich_{\hat{A}_{K}}(s_{0}) $, donc $\psi_{s'}^{\#}$ est le rel\`evement de Teichm¬\"uller de $s_{0}$
$$
\tau_{\hat{A}}^{\#}(s_{0}): \mathcal{Y}_{s'}=Spf\ \mathcal{V}(s_{0})\rightarrow\mathcal{Y}_{\alpha}=Spf\ \hat{A}
$$
dont la commutation au Frobenius [(3.1.11)] va nous \^etre fort utile: $\tau_{\hat{A}_{K}}^{\ast}(s_{0})$
est une r\'ealisation de $i_{s_{0}}^{\ast}$ [B3; (2.3.6) et (2.3.2)(iv)], [LS, chap 7] qui commute aux Frobenius. Ainsi, gr\^ace \`a [B3, (2.1.10), (2.2.3), (2.2.10), (2.3.2)] on a d\'emontr\'e la proposition suivante:\\

\vskip 2mm
\noindent \textbf{Proposition (4.1.2.12)}. \textit{Avec les notations pr\'ec\'edentes, soit $\mathcal{F}(resp\ \mathcal{G})\in Isoc(S_{\alpha}/K)$ un isocristal convergent dont une r\'ealisation est un $\mathcal{O}_{\mathcal{S}_{\alpha K}}$-module coh\'erent $\mathcal{F}_{K} (resp\ \mathcal{G}_{K})$ \`a connexion int\'egrable convergente et $\psi:\mathcal{F}\rightarrow\mathcal{G}$ un morphisme de $Isoc(S_{\alpha}/K)$ dont une r\'ealisation est un morphisme $\psi_{K}:\mathcal{F}_{K}\rightarrow\mathcal{G}_{K}$.} 
\begin{enumerate}
\item[(i)]\textit{Pour tout point ferm\'e $s_{0}$ de $S_{\alpha}$, une r\'ealisation de la fibre $i_{s_{0}}^{\ast}(\mathcal{F})=:\mathcal{F}_{s_{0}}$ de $\mathcal{F}$ en $s_{0}$ est donn\'ee par la fibre $(\theta_{s}^{\#})^{\ast}(\mathcal{F}_{K})=:\mathcal{F}_{K,s}$ de $\mathcal{F}_{K}$ en un point quelconconque $s\in ]s_{0}[$. La connexion dont est muni $\mathcal{F}_{K}$ fournit un syst\`eme compatible d'isomorphismes canoniques 
$$
s_{1}^{\ast}(\mathcal{F}_{K})\overset{\sim}{\rightarrow}s_{2}^{\ast}(\mathcal{F}_{K})
$$
pour tous points $s_{1},s_{2}\in Spm\ \hat{A}_{K}$ ayant m\^eme sp\'ecialisation $s_{0}$, i.e. pour des points $s_{i}\in]s_{0}[=\overline{sp}^{-1}(s_{0}), i=1,2$.}
\item[(ii)]\textit{Pour tout point ferm\'e $s_{0}$ de $S_{\alpha}$, $\tau_{\hat{A}_{K}}^{\ast}(s_{0})$
est une r\'ealisation de $i_{s_{0}}^{\ast}$, appel\'ee r\'ealisation de Teichm\"uller de $i_{s_{0}}^{\ast}$, qui commute aux Frobenius, i.e. $\tau^{\ast}_{\hat{A}_{K}}(s_{0})\circ F_{\hat{A}_{K}}^{\ast}=\sigma^{\ast}_{K(s_{0})}\circ\tau^{\ast}_{\hat{A}_{K}}(s_{0})$; donc parmi les r\'ealisations du} (i) \textit{figure $\tau_{\hat{A}_{K}}^{\ast}(s_{0})(\mathcal{F}_{K})$, appel\'ee r\'ealisation de Teichm\"uller de $i_{s_{0}}^{\ast}(\mathcal{F})$.\\
Si de plus $\mathcal{F}\in F^{a}\mbox{-}Isoc(S_{\alpha}/K)$ a pour Frobenius $\phi_{K}: F^{\ast}_{\hat{A}_{K}}(\mathcal{F}_{K})\overset{\sim}{\rightarrow}\mathcal{F}_{K}$, alors $\mathcal{F}_{s_{0}}$ a pour Frobenius $\tau_{\hat{A}_{K}}^{\ast}(s_{0})(\phi_{\mathcal{F}_{K}})=\phi_{\mathcal{F}_{K,s_{0}}}:\sigma^{\ast}_{K(s_{0})}(\mathcal{F}_{s_{0}})\overset{\sim}{\rightarrow}\mathcal{F}_{s_{0}}$.}
\item[(iii)]\textit{Pour v\'erifier que $\psi$ est un monomorphisme (resp un \'epimorphisme, resp un isomorphisme) il suffit de le v\'erifier aux points ferm\'es de $S_{\alpha}$; en particulier il suffit de v\'erifier que la fibre de $\psi_{K}$ aux points de Teichm\"uller en est un, i.e. il suffit de v\'erifier que  $\tau_{\hat{A}_{K}}^{\ast}(s_{0})(\psi_{K})$ en est un pour tout point ferm\'e $s_{0}$ de $S_{\alpha}$.}
\end{enumerate}
\vskip 3mm
On consid\`ere le diagramme commutatif \`a carr\'es cart\'esiens
$$
\xymatrix{
  X_{s_{0}}^{(q)} \ar[d]_{f_{s_{0}}^{(q)}} \ar[r]&X_{s_{0}}\ar[r] \ar[d]^{f_{s_{0}}} & X_{\alpha} \ar[d]^{f_{\alpha}}\\
 Y \ar[r] _{\sigma^{\#}_{k(s_{0})}} & Y\ar[r]_{i_{s_{0}}}&S_{\alpha}&\overset{.}{,}}
 \leqno{(4.1.2.13)}
$$
$\mathcal{E}_{\alpha}=j^{\ast}_{X_{\alpha}}(\mathcal{E})\in Isoc(X_{\alpha}/K)$ a pour images inverses
$\mathcal{E}_{X_{s_{0}}}$ et $\mathcal{E}_{X^{(q)}_{s_{0}}}$ sur $X_{s_{0}}$ et  $X^{(q)}_{s_{0}}$ respectivement. Il s'agit de d\'emontrer que l'image inverse
$$
i_{s_{0}}^{\ast}(\eta^{i})=\tau_{\hat{A}_{K}}^{\ast}(s_{0})(\eta^{i})=\eta^{i}_{s_{0}}
$$
de $\eta^{i}$ est un isomorphisme; or cette image inverse s'ins\`ere dans le diagramme commutatif
$$
\xymatrix{
 \tau_{\hat{A}_{K}}^{\ast}(s_{0})R^{i}f^{(q)}_{\alpha, rig^{\ast}}(X^{(q)}_{\alpha}/\mathcal{S}_{\alpha}; \pi^{\ast}_{X_{\alpha}/S_{\alpha}}(\mathcal{E}_{\alpha}))\ar[r]^{\eta^{i}_{s_{0}}}&  \tau_{\hat{A}_{K}}^{\ast}(s_{0})R^{i}f_{\alpha, rig^{\ast}}(X_{\alpha}/\mathcal{S}_{\alpha};F_{X_{\alpha}}^{\ast} (\mathcal{E}_{\alpha}))\ar[dd]^{\simeq}_{[Et 6, (3.4.4)]}\\
 \tau_{\hat{A}_{K}}^{\ast}(s_{0})F_{S_{\alpha}}^{\ast}R^{i}f_{\alpha, rig^{\ast}}(X_{\alpha}/\mathcal{S}_{\alpha}; \mathcal{E}_{\alpha})\ar[u]_{\simeq}^{[Et 6, (3.4.4)]}&\\
 \sigma^{\ast}_{K(s_{0})}\tau_{\hat{A}_{K}}^{\ast}(s_{0})R^{i}f_{\alpha, rig^{\ast}}(X_{\alpha}/\mathcal{S}_{\alpha}; \mathcal{E}_{\alpha})\ar[u]_{\simeq}^{[(3.1.11)']}\ar[d]^{\simeq}_{[Et 6, (3.4.4)]}&R^{i}f_{s_{0}, rig^{\ast}}(X_{s_{0}}/K(s_{0});F_{X_{s_{0}}}^{\ast} (\mathcal{E}_{X_{s_{0}}}))\ar@{=}[dd]\\
\sigma^{\ast}_{K(s_{0})}R^{i}f_{s_{0}, rig^{\ast}}(X_{s_{0}}/K(s_{0}); \mathcal{E}_{X_{s_{0}}})\ar[d]^{\simeq}_{[Et 6, (3.4.4)]}&\\
R^{i}f^{(q)}_{s_{0}, rig^{\ast}}(X^{(q)}_{s_{0}}/K(s_{0}); \mathcal{E}_{X^{(q)}_{s_{0}}})\ar@{=}[d]&H^{i}_{rig,c}(X_{s_{0}}/K(s_{0});F_{X_{s_{0}}}^{\ast} (\mathcal{E}_{X_{s_{0}}}))\\
H^{i}_{rig,c}(X^{(q)}_{s_{0}}/K(s_{0}); \mathcal{E}_{X^{(q)}_{s_{0}}})\ar[ur]_{F^{\ast}_{X_{s_{0}}/k(s_{0})}}&.
}
\leqno{(4.1.2.14)}
$$
Ainsi $\eta^{i}_{s_{0}}$ s'identifie \`a
$$
F^{\ast}_{X_{s_{0}}/k(s_{0})}:H^{i}_{rig,c}(X^{(q)}_{s_{0}}/K(s_{0}); \mathcal{E}_{X^{(q)}_{s_{0}}})\rightarrow H^{i}_{rig,c}(X_{s_{0}}/K(s_{0});F_{X_{s_{0}}}^{\ast} (\mathcal{E}_{X_{s_{0}}}))\ ,
\leqno{(4.1.2.15)}
$$
et ce morphisme est un isomorphisme d'apr\`es la proposition (4.1.2.16) ci-dessous, ce qui ach\`evera la preuve de (4.1.2.1).\\

\textit{La preuve de (4.1.2.2)} a \'et\'e donn\'ee ci-dessus par la m\^eme occasion.\\

\textit{Pour (4.1.2.3)}, on se ram\`ene, comme pour la preuve de (4.1.2.1), au cas o\`u $\mathcal{S}=Spf\ \hat{A}_{\alpha}$ et  $\mathcal{S}'=Spf\ \hat{A}'_{\alpha}$ et la fonctorialit\'e des constructions nous assure de la commutativit\'e du diagramme de (4.1.2.3).\ $\square$\\

\noindent On a la g\'en\'eralisation suivante de  [E-LS 1, 2.1]:\\
\vskip 2mm
\noindent \textbf{Proposition (4.1.2.16)}. \textit{
Supposons le corps $k$ parfait. Soient $X$ un $k$-sch\'ema s\'epar\'e de type fini, $F_{X}$ l'it\'er\'e $a$-i\`eme du Frobenius absolu de $X$ ($F_{X}^{\ast}(x)=x^{q}$), $F_{X}=\pi_{X/k}\circ F_{X/k}$ sa factorisation }
$$
\xymatrix{
X\ar[r]^{F_{X/k}}&X'\ar[r]^{\pi_{X/k}}&X
}
$$
 \textit{et $\sigma: K\rightarrow K'=K$ le rel\`evement choisi de la puissance $q=p^{a}$ de $k$  .}
 \textit{Pour $E\in Isoc^{\dag}(X/K)$, on a:}
\begin{enumerate}

\item[(i)]\textit{ Pour tout entier $i\geqslant 0$, $F_{X/k}$ induit une bijection $K$-lin\'eaire
$$F^{\ast}_{X/k}: H^{i}_{rig,c}(X'/K',\pi^{\ast}_{X/k}(E))\rightarrow H^{i}_{rig,c}(X/K,F^{\ast}_{X}(E)).$$
 Si de plus $X$ est lisse sur $k$, la m\^eme assertion vaut pour la cohomologie rigide sans supports compacts.}
\item[(ii)] \textit{Pour tout entier $i\geqslant 0$, $F_{X}$ induit une bijection $\sigma$-lin\'eaire
$$
F_{X}^{\ast}: H^{i}_{rig,c}(X/K,E)\rightarrow H^{i}_{rig,c}(X/K,F^{\ast}_{X}(E))
$$
\noindent c'est-\`a dire un isomorphisme
$$
\sigma^{\ast}(H^{i}_{rig,c}(X/K,E))\displaystyle \mathop{\rightarrow}^{\sim} H^{i}_{rig,c}(X/K,F^{\ast}_{X}(E)).
$$
\noindent Si de plus $X$ est lisse sur $k$, la m\^eme assertion vaut pour la cohomologie rigide sans supports compacts.
}
\end{enumerate}

\noindent \textit{D\'emonstration de (4.1.2.16)}. La preuve suit celle de [E-LS 1, 2.1]. \\
\noindent \textit{Pour (i)}. On notera $H^{i}(X/K,E)$ la cohomologie rigide avec ou sans supports compacts.
En utilisant la suite exacte longue de localisation en cohomologie rigide \`a supports compacts (resp. la suite spectrale de localisation lorsque $X$ est lisse sur $k$) on se ram\`ene au cas o\`u $X$ est un sous-sch\'ema de $\mathbb{P}_{k}^{n}$ qui ne rencontre pas les hyperplans de coordonn\'ees; comme la cohomologie rigide commute aux extensions finies de $K$ on peut supposer que $K$ contient les racines $q$-i\`emes de l'unit\'e. On note $F_{\mathbb{P}}$ l'endomorphisme $\sigma$-lin\'eaire de $\mathbb{P}=\mathbb{P}_{\mathcal{V}}^{n}$ d\'efini par $F_{\mathbb{P}}^{^{\ast}}(T_{i})=T_{i}^{q}$ pour $i\in \llbracket0,n\rrbracket$; on a la factorisation usuelle de $F_{\mathbb{P}}$

$$
\xymatrix{\mathbb{P}_{\mathcal{V}}^{n} \ar@/^1pc/[rrd]^{F_{\mathbb{P}}} \ar@/_/[rdd]_{g}  \ar@{.>}[rd]^{F=F_{\mathbb{P}/\mathcal{V}}} \\
&  \mathbb{P}_{\mathcal{V}}^{n(\sigma)} \ar[d]^{g^{(\sigma)}} \ar[r]^{\pi_{\mathbb{P}/\mathcal{V}}} & \mathbb{P}_{\mathcal{V}}^{n} \ar[d]^g\\
& Spec(\mathcal{V}) \ar[r] ^{\sigma} & Spec(\mathcal{V})
}
$$

\noindent o\`u le carr\'e est cart\'esien: puisque $k$ est parfait, $\sigma$ est un isomorphisme, de m\^eme que $\pi_{\mathbb{P}/\mathcal{V}}$. En dehors des hyperplans de coordonn\'ees, le morphisme
$$
F_{K}=(F_{\mathbb{P}/\mathcal{V}})_{K}: \mathbb{P}_{K}^{n}\rightarrow  \mathbb{P}_{K}^{n(\sigma)}
$$ 
est un rev\^etement \'etale galoisien de groupe $\mu_{q}^{n}\simeq (\mathbb{Z}/q\mathbb{Z})^{n}$. Or, si $V$ d\'esigne un voisinage strict quasi-Stein suffisamment petit du tube de $X'$, alors $F_{K}$ induit un rev\^etement \'etale galoisien encore not\'e $F_{K}:W\rightarrow V$ de groupe $\mu_{q}^{n}\simeq (\mathbb{Z}/q\mathbb{Z})^{n}$. On peut supposer de plus que $E'= \pi_{X/k}^{\ast}(E)$ provient d'un module \`a connexion $\mathcal{M}$ sur $V$: tout $F$-automorphisme $\psi$ de $\mathbb{P}_{\mathcal{V}}^{n}$  induit un automorphisme $\psi^{\ast}$ de $F_{K}^{\ast}\mathcal{M}\otimes \Omega^{\bullet}$ et l'endomorphisme $F_{\ast}(\sum \psi^{\ast})$ de $F_{K^{\ast}}F_{K}^{\ast}\mathcal{M}\otimes \Omega^{\bullet}_{V}$ se factorise de mani\`ere unique par le morphisme de complexes
$$
F_{K}^{\ast}: \mathcal{M}\otimes \Omega^{\bullet}_{V} \rightarrow F_{K^{\ast}}F_{K}^{\ast}\mathcal{M}\otimes \Omega^{\bullet}_{V}
$$
pour donner l'application trace
$$
Tr:F_{K^{\ast}}F_{K}^{\ast}\mathcal{M}\otimes \Omega^{\bullet}_{V} \rightarrow \mathcal{M}\otimes \Omega^{\bullet}_{V}
$$
(cf [E-LS 1, 2.1] et [Mi, V, lemma 1.12]). Cette application induit des homomorphismes\\

$tr: H^{i}(W, j_{W}^{\dag}F_{K}^{\ast}\mathcal{M}\otimes \Omega^{\bullet}) \rightarrow H^{i}(V, j_{V}^{\dag}\mathcal{M}\otimes \Omega^{\bullet})$\\

et $ tr: H^{i}_{]X[}(W,F_{K}^{\ast}\mathcal{M}\otimes \Omega^{\bullet}) \rightarrow H^{i}_{]X'[}(V, \mathcal{M}\otimes \Omega^{\bullet}).$\\

Rappelons ici la d\'efinition de $ j^{\dag}_{W}$. Pour un voisinage strict $W$ (resp. un couple de voisinages stricts $W'\subset W$) du tube de $X$ dans $\mathbb{P}_{K}^{n\ an}$ on note $\alpha_{W}$ (resp.$\alpha_{W W'})$ l'immersion ouverte de $W$ dans  $\mathbb{P}_{K}^{n\ an}$ (resp. de $W'$ dans $W$). Si $\mathcal{A}$ est un faisceau d'anneaux sur $W$ et $\mathcal{N}$ un $\mathcal{A}$-module, on pose [B 3,(2.1.1.1)], [LS, chap 5]:\\
 
 \noindent  $\qquad\qquad j^{\dag}_{W}\mathcal{N}:=\displaystyle\mathop{\mbox{lim}}_{\longrightarrow\atop{W'\subset W}} \alpha^{}_{W W'^{\ast}} \alpha^{\ast}_{W W'} \mathcal{N}$ ,\\
 \noindent la limite \'etant prise sur les voisinages $W'\subset W$.\\

 Puisque $F$ prolonge $F_{X/k}$, l'application

$$F^{\ast}_{X/k}: H^{i}(X'/K',\pi^{\ast}_{X/k}(E))\rightarrow H^{i}(X/K,F^{\ast}_{X}(E))$$

\noindent est induite par le morphisme de complexes $F_{K}^{\ast}$, et comme $Tr\circ F_{K}^{\ast}=q^{n}$ sur $\mathcal{M}\otimes \Omega^{\bullet}_{V}$, on en d\'eduit que $tr\circ F_{X/k}^{\ast}=q^{n}$ sur $H^{i}(X'/K',\pi^{\ast}_{X/k}(E))$: remarquons que pour pouvoir calculer les groupes de cohomologie rigide de $X$ sans supports \`a l'aide de complexes de de Rham comme ci-dessus on a \'et\'e amen\'e \`a supposer $X$ lisse dans l'\'enonc\'e de (4.1.2.16); l'\'enonc\'e de [E-LS 1, 2.1], pour \^etre complet, doit lui-aussi comporter l'hypoth\`ese de lissit\'e de $X$ dans le cas de la cohomologie rigide sans supports. \\

 Par d\'efinition on a $F^{\ast}_{K}\circ Tr= F_{K\ast}(\sum \psi^{\ast})$ sur $F_{K^{\ast}}F_{K}^{\ast}\mathcal{M}\otimes \Omega^{\bullet}_{V}$. D'autre part, si l'on note $\psi_{0}$ la r\'eduction mod $\pi$ d'un $F$-automorphisme $\psi$ de $\mathbb{P}_{\mathcal{V}}^{n}$, alors, pour tout $i \in \llbracket 0,n\rrbracket$ et $\lambda_{i}\in k$, $\psi_{0}$ v\'erifie
$$
\psi_{0}^{\ast}(\lambda_{i}^{q}t_{i}^{q})=\psi_{0}^{\ast}F_{X/k}^{\ast}(\lambda_{i}\otimes t_{i})=F_{X/k}^{\ast}(\lambda_{i}\otimes t_{i})=F_{X/k}^{\ast}(1\otimes\lambda_{i}^{q} t_{i})=\lambda_{i}^{q} t_{i}^{q};
$$
puisque $k$ est parfait on en d\'eduit que $\psi_{0}$ est l'identit\'e de $X$. Ainsi on voit que
 $$
 F_{X/k}^{\ast}\circ tr= \sum Id^{\ast}=q^{n}
 $$
sur $H^{i}(X/K,F^{\ast}_{X}(E))$. On a donc montr\'e que $(1/q^{n})tr$ est un inverse pour $F_{X/k}^{\ast}$.\\

\noindent\textit{Pour (ii)}. Le morphisme du (ii) est le compos\'e
$$
\xymatrix{
\sigma^{\ast}(H^{i}(X/K, E))\ar[r]^{ u}&H^{i}(X'/K',\pi^{\ast}_{X/k}(E))\ar[r]^{F_{X/k}^{\ast}}&H^{i}(X/K,F^{\ast}_{X}(E))
}
$$
o\`u le morphisme $u$ est induit par le changement de base $F_{k}$ (puissance $q$ sur $k$) du carr\'e cart\'esien

$$
\xymatrix{
X'\ar[r]^{\pi_{X/k}}\ar[d]^{g'}&X\ar[d]^{g}&\\
Spec\ k\ar[r]_{F_{k}}&Spec\ k&.
}
$$
\noindent Puisque $k$ est parfait, $F_{k}$ et $\pi_{X/k}$ sont des isomorphismes, donc $u$ en est un aussi. Ceci ach\`eve la preuve de (4.1.2.16).  $\square$\\

\noindent\textbf {4.1.3.} On va pr\'eciser le th\'eor\`eme (4.1.2) en l'\'etendant.\\

Soient $S$ un $k$-sch\'ema lisse et s\'epar\'e et $f\ :  X\rightarrow S$ un $k$-morphisme projectif et lisse. Notons $S= \underset{\alpha}{\bigcup} \ S_{\alpha,0}$ une d\'ecomposition de $S$ en r\'eunion d'ouverts connexes affines $S_{\alpha,0}= Spec(A_{\alpha,0}), \ A_{\alpha}= \mathcal{V}[t_{1},...,t_{d_{\alpha}}]/J_{\alpha}$ une $\mathcal{V}$-alg\`ebre lisse relevant $A_{\alpha,0}$ dont on a fix\'e une pr\'esentation et $S_{\alpha}= Spec (A_{\alpha})$. On d\'esigne par $$
f_{\alpha}:X_{\alpha,0}=X\times_{S}S_{\alpha,0}\longrightarrow S_{\alpha,0}
$$
la restriction de $f$. Quitte \`a d\'ecomposer $X_{\alpha, 0}$ en somme disjointe de ses composantes connexes on peut supposer $X_{\alpha, 0}$ connexe. D'apr\`es [Et 5, (3.3.1)] il existe un rel\`evement projectif $h_{\alpha}$ de $f_{\alpha}$, $h_{\alpha}: X_{\alpha}\rightarrow S_{\alpha}$. Le compl\'et\'e formel de ce morphisme $h_{\alpha}$ est un morphisme projectif de $\mathcal{V}$-sch\'emas formels 
$$
\hat{h}_{\alpha}: \mathcal{X}_{\alpha}\rightarrow \mathcal{S}_{\alpha}\ .
$$
\noindent On suppose de plus que, pour tout $\alpha$, $X_{\alpha}$ est plat sur $\mathcal{V}$: ceci implique [Et 5, (3.3.4)] que, pour tout $\alpha$, $\hat{h}_{\alpha}$ est lisse.\\

\noindent\textbf{Th\'eor\`eme (4.1.3.1)}. \textit{Supposons $k$ parfait. Soient $S$ un $k$-sch\'ema lisse et supposons que $f : X \rightarrow S$ est un $k$-morphisme projectif et lisse satisfaisant aux hypoth\`eses 3.3.2 pr\'ec\'edentes ou que $f$ d\'efinit $X$ comme une intersection compl\`ete relativement \`a $S$ dans un espace projectif sur $S$ [Et 5, (3.2.5)]. Alors, pour tout entier $i \geqslant 0$, $f$ induit un foncteur}
$$
R^{i} f_{conv\ast}\  : F^{a}\mbox{-}\textrm{Isoc}(X/K) \rightarrow F^{a}\mbox{-}\textrm{Isoc}(S/K)
$$
\textit{qui commute \`a tout changement de base $S' \rightarrow S$ entre $k$-sch\'emas lisses.} 

\vskip 3mm
\noindent \textit{D\'emonstration}. Compte tenu de [Et 6, (3.4.8.2), (3.4.8.6)] il s'agit de v\'erifier que le Frobenius (qui est d\'efini par fonctorialit\'e [LS, \S8]) est un isomorphisme. D'apr\`es [B-G-R, (9.4.2/7)] et comme on l'a vu dans la preuve de (4.1.2) il suffit de v\'erifier l'isomorphisme sur les points ferm\'es de $S$ et le r\'esultat provient alors de [(4.1.2.15)]. $\square$

\newpage
\subsection*{4.2. Cas fini \'etale}

Dans le cas fini \'etale on n'a pas lieu de supposer le corps $k$ parfait ni que $K$ contient les racines $q^{iemes}$ de l'unit\'e:

\vskip 3mm
\noindent \textbf{Th\'eor\`eme (4.2.1)}. \textit{Soient $S$ un $k$-sch\'ema lisse et $f : X \rightarrow S$ un $k$-morphisme fini \'etale. Alors}\\

\noindent (4.2.1.1) \textit{Pour tout entier $i \geqslant 0$, on a des foncteurs }\\

\begin{itemize}
\item[(i)] $R^{i}  f_{\textrm{conv}^{\ast}} : \textrm{Isoc}(X/K) \rightarrow \textrm{Isoc}(S/K),$

\vskip 1mm
\item[(ii)] $R^{i}  f_{\textrm{conv}^{\ast}} : F^{a}\mbox{-}\textrm{Isoc}(X/K) \rightarrow F^{a}\mbox{-}\textrm{Isoc}(S/K),$

\vskip 1mm
\item[(iii)] \textit{Pour $\mathcal{E} \in \textrm{Isoc}(X/K)$ et $i \geqslant 1$ on a }
$$
R^{i} f_{\textrm{conv}^{\ast}} \mathcal{(E)} = 0.
$$
\end{itemize}
\noindent (4.2.1.2) \textit{Supposons de plus $f$ galoisien de groupe $G$. Pour $\mathcal{E} \in \textrm{Isoc}(X/K)$ on a des isomorphismes canoniques}\\
\begin{itemize}
\item[(i)] $\mathcal{E} \displaystyle \mathop{\longrightarrow}^{\sim}\ (f_{\textrm{conv}^{\ast}}\ f^{\ast}(\mathcal{E}))^G,$
\vskip 1mm
\item[(ii)] $H^i_{\textrm{conv}}(S/K, \mathcal{E}) \displaystyle \mathop{\longrightarrow}^{\sim}\ H^i_{\textrm{conv}}(X/K, f^{\ast}(\mathcal{E}))^G,$
\vskip 1mm
\item[(iii)] \textit{Si  $\mathcal{E} \in F^{a}\mbox{-}\textrm{Isoc}(S/K)$, ces isomorphismes sont compatibles aux Frobenius.}
\end{itemize}

\vskip 3mm
\noindent \textit{D\'emonstration.\\
Pour (4.2.1.1)}. Le (i) est l\`a pour m\'emoire, car prouv\'e en [Et 6, 3.4.8]. On a vu dans la d\'emonstration de [loc. cit.] que la d\'efinition de $R^{i}  f_{\textrm{conv}^{\ast}}\mathcal{(E)}$ est locale sur $S$ : on peut donc supposer $S = \textrm{Spec}\ A_{0}$ affine et lisse sur $k$.\\

Posons $\mathcal{S} = \textrm{Spf}\ \hat{A}$ o\`u $A$ est une $\mathcal{V}$-alg\`ebre lisse relevant $A_{0}$, et relevons $f : X \rightarrow S$ en un morphisme fini \'etale de $\mathcal{V}$-sch\'emas formels $h : \mathcal{X} \rightarrow \mathcal{S}$ [EGA IV, (18.3.2) ou (18.3.4)], et soit $F_{\mathcal{S}} : \mathcal{S} \rightarrow \mathcal{S}$ un rel\`evement du Frobenius de $S$. Puisque $f$ est \'etale, dans la d\'ecomposition classique du Frobenius $F_{X}$ de $X$

$$
\xymatrix{X \ar@/^1pc/[rrd]^{F_{X}} \ar@/_/[rdd]_{f}  \ar@{.>}[rd]^{F_{X/S}} \\
&  X^{(q)} \ar[d]^{f^{(q)}} \ar[r]^{\pi_{X/S}} & X \ar[d]^f\\
& S \ar[r] _{F_{S}} & S
}
$$

\noindent le morphisme $F_{X/S}$ est un isomorphisme et se rel\`eve de mani\`ere unique en un isomorphisme $F_{\mathcal{X}/\mathcal{S}}$ s'ins\'erant dans le diagramme commutatif \`a carr\'e cart\'esien

$$
\xymatrix{\mathcal{X}  \ar@/_/[rdd]_{h}  \ar[rd]^{F_{\mathcal{X}/\mathcal{S}}}  &\\
&  \mathcal{X}^{(q)} \ar[d]^{h^{(q)}} \ar[r]^{\pi_{\mathcal{X}/\mathcal{S}}} & \mathcal{X} \ar[d]^h &\\
& \mathcal{S} \ar[r] _{F_{\mathcal{S}}} & \mathcal{S} & .
}
$$\\
On pose $F_{\mathcal{X}} = \pi_{\mathcal{X}/\mathcal{S}} \circ F_{\mathcal{X}/\mathcal{S}}.$\\

Pour $\mathcal{E} \in \textrm{Isoc}(X/K)$, soit $\mathcal{E}_{\mathcal{X}}$ une r\'ealisation de $\mathcal{E}$ sur $\mathcal{X}_{K}$ ; par d\'efinition on a

$$
\xymatrix{
R^{i} f_{\textrm{conv}^{\ast}}(X/ \mathcal{S};\mathcal{E}) = H^{i}(\mathbb{R} h_{K^{\ast}}
(\mathcal{E}_{\mathcal{X}} \otimes \Omega^{\bullet}_{\mathcal{X}_{K}/\mathcal{S}_{K}}))\\
= R^{i} h_{K^{\ast}} (\mathcal{E}_{\mathcal{X}}),
}
$$

\noindent car $\mathcal{X}$ est \'etale sur $\mathcal{S}$. D'o\`u le (iii) par le th\'eor\`eme B de Kiehl car $h$ est affine. \\

Soit $\mathcal{E} \in F^{a}\mbox{-}\textrm{Isoc}(X/K)$ et $\phi : F^{\ast}_{\mathcal{X}}
(\mathcal{E}_{\mathcal{X}}) \displaystyle \mathop{\longrightarrow}^{\sim} \mathcal{E}_{\mathcal{X}}$ le Frobenius.\\
D'apr\`es [cor (1.2.3)] il suffit de construire un isomorphisme $\phi^i$ (de Frobenius) sur $R^{i} f_{\textrm{conv}^{\ast}}(\mathcal{E})$, compatible aux connexions. Comme $F_{\mathcal{S}}$ est plat, le morphisme de changement de base 
$$
F^{\ast}_{\mathcal{S}_{K}} R^{i} f_{\textrm{conv}^{\ast}}(X/\mathcal{S}, \mathcal{E}) \rightarrow R^{i} f_{\textrm{conv}^{\ast}}^{(q)}(X^{(q)}/\mathcal{S},\mathcal{E}) \simeq R^{i} h_{K^{\ast}}^{(q)}(\pi^{\ast}_{\mathcal{X}_{K}/\mathcal{S}_{K}}(\mathcal{E}_{\mathcal{X}}))
$$
\noindent est un isomorphisme [Et 6, (3.4.4)]; par composition avec les isomorphismes

$$
R^{i} h_{K^{\ast}}^{(q)}(\pi^{\ast}_{\mathcal{X}_{K}/\mathcal{S}_{K}}(\mathcal{E}_{\mathcal{X}})) \simeq R^{i} h_{K^{\ast}}(F^{\ast}_{\mathcal{X}_{K}}(\mathcal{E}_{\mathcal{X}}))
$$

 \noindent (puisque $F_{\mathcal{X}/\mathcal{S}}$ est un isomorphisme)

\noindent et

$$
R^{i}  h_{K^{\ast}}(\phi) : R^{i}  h_{K^{\ast}}(F^{\ast}_{\mathcal{X}_{K}}(\mathcal{E}_{\mathcal{X}})) \simeq R^{i}  h_{K^{\ast}}(\mathcal{E}_{\mathcal{X}}) ,
$$

\noindent on obtient le Frobenius $\phi^i$ cherch\'e

$$
\phi^i : F^{\ast}_{\mathcal{S}_{K}}  R^{i} f_{\textrm{conv}^{\ast}}(X/\mathcal{S}, \mathcal{E}) \simeq R^{i} f_{\textrm{conv}^{\ast}}(X/\mathcal{S}, \mathcal{E}).
$$

En reprenant la preuve de [Et 6, (3.4.4)] on v\'erifie que $\phi^i$ est compatible aux connexions, d'o\`u le (ii).\\

\textit{Pour (4.2.1.2)}. Soit $\mathcal{E}_{\mathcal{S}}$ une r\'ealisation de $\mathcal{E}$ sur $\mathcal{S}_{K}$. Par d\'efinition on a

$$
f_{\textrm{conv}^{\ast}}(X/\mathcal{S}, f^{\ast}(\mathcal{E})) = h_{K^{\ast}}
(h^{\ast}_{K}(\mathcal{E}_{\mathcal{S}})).
$$

\noindent Comme $h_{K}$ est fini \'etale galoisien de groupe $G$ [Et 6, (2.3.1)], la fl\`eche canonique
$$
\mathcal{E}_{\mathcal{S}} \rightarrow (h_{K^{\ast}} (h^{\ast}_{K}(\mathcal{E}_{\mathcal{S}})))^G
$$

\noindent est un isomorphisme, d'o\`u (i).\\

L'isomorphisme du (ii) est alors une cons\'equence classique du (i) [Et 1, III, 3.1.1].\\

La fonctorialit\'e des constructions pr\'ec\'edentes prouve le (iii).  $\square$

 \end{document}